\documentclass{m2an}
\usepackage{graphicx}%
\usepackage{multirow}%
\usepackage{amsmath,amssymb,amsfonts}%
\usepackage{amsthm}%
\usepackage{mathrsfs}%
\usepackage{appendix}%
\usepackage{xcolor}%
\usepackage{textcomp}%
\usepackage{manyfoot}%
\usepackage{booktabs}%
\usepackage{algorithm}%
\usepackage{algpseudocode}%
\usepackage{listings}%

\usepackage[T1]{fontenc}
\usepackage{lmodern}
\usepackage{subcaption}
\usepackage{mathtools}
\usepackage{here}
\usepackage{hyperref}
\usepackage{enumerate}
\usepackage{cleveref}
\usepackage{dsfont}
\usepackage{enumitem}
\usepackage{nicefrac}

\newcommand{\bbE}{{\mathbb{E}}}

\newcommand{\bbP}{{\mathbb{P}}}

\newcommand{\bbV}{{\mathbb{V}}}

\newcommand{\N}{{\mathbb{N}}} 
\newcommand{\R}{{\mathbb{R}}} 

\DeclareSymbolFont{bbold}{U}{bbold}{m}{n}
\DeclareSymbolFontAlphabet{\mathbbold}{bbold}

\newcommand{\calB}{{\mathcal{B}}}

\newcommand{\calF}{{\mathcal{F}}}

\newcommand{\calL}{{\mathcal{L}}}

\newcommand{\calO}{{\mathcal{O}}}

\newcommand{\calR}{{\mathcal{R}}}
\newcommand{\calS}{{\mathcal{S}}}

\newcommand{\calV}{{\mathcal{V}}}
\newcommand{\calW}{{\mathcal{W}}}

\newcommand{\abs}[1]{\left\lvert#1\right\rvert}

\DeclareMathOperator*{\argmin}{arg\,min}
\DeclareMathOperator*{\esssup}{ess\,sup}
\newcommand{\norm}[1]{\left\lVert#1\right\rVert}
\newcommand{\ceil}[1]{\left\lceil#1\right\rceil}
\newcommand{\floor}[1]{\left\lfloor#1\right\rfloor}
\newcommand{\inp}[2]{\left \langle #1, #2 \right \rangle}

\newcommand{\Exp}[1]{\mathbb{E}\left[#1\right]}

\newcommand{\algname}{SG-LSCV }
\newcommand{\algnamens}{SG-LSCV}
\newcommand{\galgname}{SG-MCV }
\newcommand{\galgnamens}{SG-MCV}

\crefname{thrm}{theorem}{theorems}
\Crefname{thrm}{Theorem}{Theorems}
\crefname{lmm}{lemma}{lemmas}
\Crefname{lmm}{Lemma}{Lemmas}
\crefname{crllr}{corollary}{corollaries}
\Crefname{crllr}{Corollary}{Corollaries}
\crefname{prpstn}{proposition}{propositions}
\Crefname{prpstn}{Proposition}{Propositions}
\crefname{crtrn}{criterion}{criteria}
\Crefname{crtrn}{Criterion}{Criteria}
\crefname{lgrthm}{algorithm}{algorithms}
\Crefname{lgrthm}{Algorithm}{Algorithms}
\crefname{dfntn}{definition}{definitions}
\Crefname{dfntn}{Definition}{Definitions}
\crefname{cnjctr}{conjecture}{conjectures}
\Crefname{cnjctr}{Conjecture}{Conjectures}
\crefname{xmpl}{example}{examples}
\Crefname{xmpl}{Example}{Examples}
\crefname{prblm}{problem}{problems}
\Crefname{prblm}{Problem}{Problems}
\crefname{rmrk}{remark}{remarks}
\Crefname{rmrk}{Remark}{Remarks}
\crefname{nt}{note}{notes}
\Crefname{nt}{Note}{Notes}
\crefname{clm}{claim}{claims}
\Crefname{clm}{Claim}{Claims}
\crefname{smmr}{summary}{summaries}
\Crefname{smmr}{Summary}{Summaries}
\crefname{cs}{case}{cases}
\Crefname{cs}{Case}{Cases}
\crefname{bsrvtn}{observation}{observations}
\Crefname{bsrvtn}{Observation}{Observations}
\crefname{asmptn}{assumption}{assumptions}
\Crefname{asmptn}{Assumption}{Assumptions}

\begin{document}
\title{Stochastic gradient with least-squares control variates}
\runningtitle{Stochastic gradient with least-squares control variates}
\runningauthors{F. Nobile, M. Raviola, N. Schaeffer}
%
\author{Fabio Nobile} \address{Institute of Mathematics, \'{E}cole Polytechnique F\'{e}d\'{e}rale de Lausanne, Station 8, CH-1015 Lausanne, Switzerland; \email{fabio.nobile@epfl.ch; matteo.raviola@epfl.ch; schaeffernathan@hotmail.com}}
\author{Matteo Raviola} \sameaddress{1}
\author{Nathan Schaeffer} \sameaddress{1}
\date{Accepted July 28, 2026}
\begin{abstract}
    The stochastic gradient descent (SGD) method is a widely used approach for solving stochastic optimization problems, but its convergence is typically slow.
    Existing variance reduction techniques, such as SAGA, improve convergence by leveraging stored gradient information; however, they are restricted to settings where the objective functional is a finite sum, and their performance degrades when the number of terms in the sum is large.
    In this work, we propose a novel approach which is well suited when the objective is given by an expectation over random variables with a continuous probability distribution.
    Our method constructs a control variate by fitting a linear model to past gradient evaluations using weighted discrete least-squares, effectively reducing variance while preserving computational efficiency.
    We establish theoretical sublinear convergence guarantees and demonstrate the method's effectiveness through numerical experiments on random PDE-constrained optimization problems.
\end{abstract}
\maketitle
\section{Introduction} \label{sec1}

{\color{black} In this paper we consider the optimization problem
    \begin{equation}\label{eq:stoch-opt-prob-intro}
        \min_{u \in U} J(u) := \mathbb{E}_{Y \sim \rho}\left[ g(u,Y) \right]
    \end{equation}
    over the Hilbert search space $U$, where $g$ is a smooth function measuring the \textit{cost} of using the design $u$ for a realization $Y \sim \rho$ of the random model parameters distributed according to the probability distribution $\rho$.
    We focus on problems where the distribution $\rho$ is known in the sense that it can be sampled from and its density can be evaluated.
    These are known in the literature as stochastic optimization problems \cite[Section 1.2]{ermoliev1988numerical} or stochastic programming problems \cite[Chapter 1]{shapiro2021lectures}.
    Furthermore, we assume the evaluation of $g$ and its gradient is computationally intensive.
    These problems, where the generative model is known yet expensive to query, arise in a wide range of scientific and engineering applications:
    examples include PDE-constrained optimization under uncertainty, maximum a posteriori (MAP) estimation for Bayesian inverse problems with complex likelihood models, and design under uncertainty in fields such as climate modeling, structural mechanics, and subsurface flow.

    We are particularly interested in the case where $J$ is strongly convex and $g$ has Lipschitz continuous gradients.
    In this setting, to solve problem \eqref{eq:stoch-opt-prob-intro}, one may consider a gradient based method such as steepest descent (GD), namely the iterative method
    \begin{equation*}
        u_{k+1} = u_k - \tau \nabla J(u_k), \quad k \geq 0,
    \end{equation*}
    with the initialization $u_0 \in U$.
    Under our assumptions and a suitable choice of the step size $\tau$, GD can be shown to converge exponentially fast in the number of iterations, i.e., $\norm{u_k  -  u^\ast}_U \leq  C\gamma^k$ for some $\gamma  \in  (0, 1)$ and a constant $C > 0$, where $u_k$ denotes the $k$-th iterate of the method and $u^\ast$ the (unique) minimizer of $J$.
    However, it is important to note that, for a given point $u$, evaluating the gradient $\nabla J(u)$  entails the computation of the full expectation with respect to the randomness in $Y$.}
The practical limitation of this approach is that $\Exp{\nabla g(u,Y)}$ cannot be computed exactly in general and any accurate discretization of it by, e.g., a Monte Carlo estimator may lead to an excessively high cost per iteration.
This issue is particularly pronounced in contexts where gradient evaluations are expensive, as for instance in PDE-constrained optimization, where evaluating the gradient of $g(\cdot,y)$ typically involves the numerical solution of one or more PDEs.

A popular technique to circumvent this limitation is the stochastic gradient descent (SGD) method, a foundational algorithm for stochastic optimization.
In its Robbins--Monro version \cite{10.1214/aoms/1177729586}, it reads
\begin{equation} \label{eq:robbins-monro-intro}
    u_{k+1} = u_k - \tau_k \widehat{\nabla J}_{\mathrm{SGD}}(u_k),
    \quad \widehat{\nabla J}_{\mathrm{SGD}}(u_k) := \nabla g(u_k,Y_k),
\end{equation}
so that only one gradient $\nabla g(u_k,Y_k)$ is evaluated at each iteration, for a randomly drawn $Y_k \sim Y$ independent of all previous draws, and the convergence is achieved by reducing the step size $\tau_k$ over the iterations.
This makes the cost of each iteration affordable.
It can be shown \cite{bottou2018optimization} that the root mean squared error $\Exp{\norm{u_k - u^\ast}_U^2}^{1/2}$ decays at the rate $1/\sqrt{k}$.
That is, in order to gain an order of magnitude in accuracy one must increase the number of iterations by a factor $100$.
The simple estimator $\widehat{\nabla J}_{\mathrm{SGD}}(u_k)$ in \eqref{eq:robbins-monro-intro} can be replaced by a mini-batch estimator
\begin{equation*}
    \widehat{\nabla J}_{\mathrm{SGD}}(u_k)
    := \frac{1}{s} \sum_{j=1}^s \nabla g(u_k,Y_{k,j}),
\end{equation*}
however this does not change the asymptotic convergence properties of the method if $s$ is kept small and fixed over the iterations.

\subsection{Variance reduction in stochastic optimization}
Over the past decade, several variance-reduced (VR) variants of SGD have been developed to address the slow convergence of classical stochastic methods.
Some of these approaches introduce control variates \cite{ripley2009stochastic}  -- auxiliary quantities with known expectations -- aimed at reducing the variance of gradient estimators with $\calO(1)$ gradient evaluations per iteration, thereby keeping the favorable cost per iteration of SGD.
Most prominently, this line of research has been pursued in the context of finite-sum optimization problems of the form
\begin{equation} \label{eq:stoch-opt-prob-finite-sum-intro}
    \min_{u \in U} J^s(u) := \frac{1}{s} \sum_{i=1}^s g_i(u),
\end{equation}
where $g_i : U \to \mathbb{R}$, $i=1,\dots,s$.
Examples include SAG \cite{schmidt2017minimizing}, which maintains a memory of the most recent gradient for each $g_i$ and uses their average to estimate the gradient of $J^s$.
SAGA \cite{defazio2014saga} improves on this by constructing an unbiased estimator via control variate correction: at each iteration, it draws one random index $i_k$ uniformly on $\{1,\dots,s\}$, it updates the design with
\begin{equation*}
    \widehat{\nabla J^s}_{\mathrm{SAGA}}(u_k) := \nabla g_{i_k}(u_k) - v_{i_k}^{(k)} + \frac{1}{s} \sum_{i=1}^s v_i^{(k)},
\end{equation*}
where $v_i^{(k)}$ stores the most recent evaluation of $\nabla g_i$, and subsequently updates the memory $v_{i_k}^{(k+1)} := \nabla g_{i_k}(u_k)$.
SVRG \cite{johnson2013accelerating}, instead, introduces an outer loop in which a full gradient of $J^s$ is computed at a reference point, which is used by inner-loop updates to obtain a variance-reduced estimate of the gradient.
Finally, SARAH \cite{nguyen2017sarah} employs a recursive gradient estimator initialized by a full gradient but corrected in a backward fashion.
Under strong convexity assumption of $J^s$ and Lipschitz continuous gradient of each of the $g_i$'s, these methods achieve mean squared error decay of the form $\mathbb{E}[\|u_k - u^\ast\|^2_U] \lesssim \gamma^k$ with $1 - \gamma \sim \mathcal{O}(s^{-1})$, that is they recover the exponential convergence of GD, yet with a rate which deteriorates in $s$.

Although powerful, these methods critically rely on the finite-sum structure of the objective in \eqref{eq:stoch-opt-prob-finite-sum-intro}.
When applied to expectations over continuous random variables, as in \eqref{eq:stoch-opt-prob-intro}, they require an a priori discretization of $J$ via a quadrature rule which, however, introduces a quadrature error that should be carefully controlled.
A small error requires many quadrature points, thus compromising the performance of these variance-reduced SGD methods.
On the other hand, a large error will quickly dominate the optimization error after a few iterations thus nullifying the advantage of variance reduction techniques.

These limitations call for the development of VR algorithms which can be directly applied to the general expectation setting in \eqref{eq:stoch-opt-prob-intro}, with random variables whose state space is excessively large or continuous.
Algorithms that address this challenge include the SCSG algorithm \cite{lei2017less} which replaces the full gradient in SVRG with a large mini-batch estimate, allowing application to cases where $s$ is not finite or too large.
It performs a geometrically distributed number of inner iterations using control variates and achieves $\mathbb{E}[\|u_k - u^\ast\|^2_U] \leq \epsilon$ with a sample complexity of $\mathcal{O}(\log(\epsilon^{-1})/\epsilon)$ in the strongly convex case.
SPIDER \cite{fang2018spider}, based on SARAH, avoids full gradient computations altogether and instead uses recursive updates initialized and corrected via mini-batches; although promising, its complexity has only been analyzed in the nonconvex setting, in which case it achieves optimal (algebraic) sample complexity.
STORM \cite{cutkosky2019momentum} avoids using giant batch sizes proposing a variant of momentum to achieve variance reduction in non-convex optimization, and, through adaptive learning rates, reduces the need for hyperparameter tuning.
To solve \eqref{eq:stoch-opt-prob-intro}, inspired by multilevel approximations, the authors of MICE \cite{carlon2020multi} proposed a telescoping sum idea where the gradient at the latest iterate is expressed in terms of gradient differences involving the design at previous iterations.
This approach, however, features the same asymptotic complexity of SGD.
This suggests that there remains room for developing variance-reduced algorithms tailored to strongly convex stochastic optimization problems over continuous random variables.
\textcolor{black}{All of these algorithms, however, were developed in the context of sample-based stochastic optimization, that is  they do not rely on the knowledge of the underlying distribution $\rho$ beyond sampling from it.
    We hence aim to design more efficient control variates when the density of $\rho$ is known.
}

Related variance reduction ideas have been explored in other domains as well.
For instance, control variates have been applied to stochastic gradient Langevin dynamics in Bayesian inference \cite{baker2019control}, and to discrete latent-variable models \cite{titsias2022double}.

\subsection{Summary of main contributions}
{\color{black}
    This work develops a variance-reduced stochastic gradient method for optimization problems of the form \eqref{eq:stoch-opt-prob-intro}, where the objective is an expectation over a continuous or high-dimensional random variable.
    The central idea is to use the known distribution of $Y$ and the regularity of the parameter-to-gradient map $y \mapsto \nabla g(u,y)$ to construct control variates directly in the continuous parameter space, without first replacing the expectation by a finite-sum quadrature approximation.
    The main contributions are as follows.
    \begin{enumerate}
        \item We introduce a memory-based control variate framework in which past sample-gradient pairs are used to fit a surrogate of $y \mapsto \nabla g(u_k,y)$.
              The resulting stochastic gradient estimator keeps the one-sample-per-iteration structure of SGD while reducing variance through an explicitly integrable surrogate.

        \item We specialize this framework to optimal weighted least-squares approximation.
              This yields the \algname algorithm, in fixed and variable approximation-space variants, and allows the control variate to exploit problem-dependent approximation spaces, such as polynomial spaces for parametric PDE-constrained optimization problems.

        \item We prove convergence guarantees in the strongly convex setting with Lipschitz-continuous gradients.
              For a fixed approximation space, the method converges exponentially up to an error floor determined by the best approximation of $\nabla g(u^\ast,\cdot)$ in that space.
              For growing approximation spaces, the error converges to zero at rates governed by the approximability of $\nabla g(u^\ast,\cdot)$: subexponential rates for exponential projection-error decay and algebraic rates for algebraic projection-error decay.

        \item We analyze the computational work of the variable-space method and show that the favorable stochastic cost structure is preserved asymptotically.
              In particular, the algorithm requires one gradient evaluation per iteration when the sampling measure is fixed, and only a controlled additional cost when resampling is used for changing optimal sampling measures.

        \item We demonstrate the method on random PDE-constrained optimization problems.
              The experiments illustrate the predicted tradeoff between approximation-space dimension, convergence speed, and error plateaus, and show that variable-space \algname can outperform SAGA-type finite-sum variance reduction without relying on high-dimensional tensorized quadrature discretizations of the objective.
    \end{enumerate}
}

\subsection{Other related approaches}
A popular approach to mitigate the cost of computing expectations in the context of scientific computing is the multilevel Monte Carlo (MLMC) method \cite{doi:10.1137/16M109870X}, which combines hierarchies of approximations and performs the bulk of the computation on coarse and cheap-to-query model discretizations.
The MLMC paradigm can be integrated with stochastic optimization methods \cite{martin2019multilevel,hu2024multi,frikha2016multi,dereich2019general,ganesh2023gradient}, and may be combined with our approach as well.
In contrast, the authors of \cite{guth2024one}, propose to reduce the cost of optimizing expensive to evaluate objectives by constructing a surrogate of the complex computational model that is learned during the optimization, which is however not used as a control variate.

Another idea to recover the GD convergence properties, instead of using variance reduction techniques, is to focus on adaptive or geometry-aware step size strategies.
In the interpolation setting for overparametrized ML models \cite{vaswani2019painless} show that a stochastic variant of the classical Armijo line search yields exponential convergence rates for SGD.
Instead, methods such as \cite{bollapragada2018adaptive,franchini2024stochastic} increase the mini-batch size geometrically to recover exponential convergence, hence asymptotically incurring the same cost as full gradient methods.

\subsection{Organization of the paper}
The structure of the paper is as follows.
In \Cref{sec:alg-desc}, we rigorously define the stochastic optimization problem we aim to solve and outline our core algorithmic idea involving control variates.
In \Cref{sec:least-squares-revision}, we review the theoretical foundations of optimal weighted least-squares approximation from \cite{cohen2016optimal}, employed to construct efficient control variates.
\Cref{sec:ls-mb-cv} provides detailed exposition of the resulting algorithm, \algname, and its implementation using least-squares. We subsequently establish theoretical convergence guarantees under the assumptions of strong convexity and Lipschitz-continuous gradients.
Finally, \Cref{sec:numerics} demonstrates the efficacy and computational benefits of our proposed methodology through numerical experiments on random PDE-constrained optimization problems.

\section{The \algname algorithm}

\subsection{Stochastic gradient with memory-based control variate}
\label{sec:alg-desc}

Let $(U, \inp{\cdot}{\cdot}_U)$ be a Hilbert space and recall that we are concerned with the solution of a stochastic optimization problem over $U$, that is we aim to minimize an objective function of the form
\begin{equation}\label{eq:stoch-opt-prob}
    J(u) = \mathbb{E}\left[ g(u,Y) \right],
\end{equation}
where $g(u,y)$ is the cost of using the design $u$ when the random model parameters are $Y=y$.
More specifically, we let $Y$ be a random vector taking values in $\Gamma \subseteq \R^l$ with joint probability measure $\rho$ on $(\Gamma,\calB(\Gamma))$, where $\calB(\Gamma)$ denotes the Borel $\sigma$-algebra.
For $1 \leq p \leq +\infty$, we define the weighted $L^p_\rho(\Gamma;U)$ spaces as
\begin{equation*}
    L^p_\rho(\Gamma; U)
    = \left\{ h: \Gamma \to U \; \mathrm{strongly \; measurable} \; : \; \norm{h}_{L^p_\rho(\Gamma;U)} := \left( \int_\Gamma  \norm{h(y)}_{U}^p d\rho(y) \right)^{1/p} < +\infty \right\}
\end{equation*}
if $1 \leq p < +\infty$ and
\begin{equation*}
    L^\infty_\rho(\Gamma; U)
    = \left\{ h: \Gamma \to U \; \mathrm{strongly \; measurable} \; : \; \norm{h}_{L^\infty_\rho(\Gamma;U)} := \esssup_{y \in \Gamma} \norm{h(y)}_{U} < +\infty \right\}.
\end{equation*}
Moreover, for any $y \in \Gamma$ fixed we let $u \mapsto g(u,y)$ be a $C^1$ function, namely we assume it is continuously Fr{\'e}chet differentiable with $\nabla g(u,y)$ defined as
\begin{equation*}
    \inp{\nabla g(u,y)}{v}_{U} := \lim_{\varepsilon \to 0} \frac{g(u+\varepsilon v,y) - g(u,y)}{\varepsilon}, \quad \forall \, v \in U.
\end{equation*}
Finally, we assume $J$ and $g$ to satisfy the following.

\begin{asmptn}[Strong convexity] \label{asmp:str-cvx}
    The function $J$ is $\alpha$-strongly convex for some $\alpha > 0$, that is
    \begin{equation*}
        J(u) \geq J(u') + \inp{\nabla J(u')}{u-u'}_U +\frac{\alpha}{2} \norm{u-u'}_U^2,
    \end{equation*}
    for all $u,u' \in U$.
\end{asmptn}

\begin{asmptn}[Lipschitz continuity]\label{asmp:lip-cont}
    The function $g(\cdot,y)$ has Lipschitz continuous gradients with respect to $u$ with a constant $\ell(y)$, $y \in \Gamma$, such that $L := \bbE_{Y\sim \rho}[\ell(Y)^2]^{\frac12} < \infty$.
\end{asmptn}
\noindent
It is easy to see that \Cref{asmp:str-cvx} implies that $J$ is bounded from below, coercive, weakly lower semicontinuous, and convex, which is enough to show that there exists a unique minimizer $u^\ast \in U$.

Recall the SGD method which reads, in its Robbins--Monro version,
\begin{equation*}
    u_{k+1} = u_k - \tau_k \widehat{\nabla J}_{\mathrm{SGD}}(u_k),
    \quad \widehat{\nabla J}_{\mathrm{SGD}}(u_k) := \nabla g(u_k,Y_k),
\end{equation*}
where $Y_k \sim \rho$ is drawn independently of all previous draws.
Borrowing ideas from works on variance-reduced SGD variants, such as SAGA \cite{defazio2014saga} -- which are able to recover an exponential convergence in the number of iterations by exploiting gradients computed at past iterations in order to build an approximation of $\nabla J$ -- we propose the following construction to solve problem \eqref{eq:stoch-opt-prob}.
At each iteration $k$, build a surrogate model $v_k:\Gamma \to U$ of $\nabla g(u_k,\cdot)$ based on prior gradient evaluations in order to reduce the variance of the SGD update rule.
More precisely, let us introduce nested spaces $\mathbb{V}_0 \subset \mathbb{V}_1 \subset \cdots \subset \mathbb{V}_{\infty} \subseteq L^2_\rho(\Gamma;U)$ of functions from $\Gamma$ to $U$ with $\mathrm{dim}(\bbV_m) = m$ and $\overline{\bbV}_{\infty} = L^2_\rho(\Gamma;U)$; then, at iteration $k$, we fit $v_k \in \bbV_{m_k}$ to $\nabla g(u_k,\cdot)$ using the history of the last $s_k \geq m_k$ iterate-gradient pairs $\mathcal{S}_k := \{ (u_j, \nabla g(u_j, Y_j)) \}_{j=k-s_k}^{k-1}$, which are kept in memory.
Then, the update rule is
\begin{equation*}
    u_{k+1} = u_k - \tau_k \widehat{\nabla J}_{\galgnamens}(u_k),
    \quad \widehat{\nabla J}_{\galgnamens}(u_k) := \nabla g(u_k,Y_k) - v_k(Y_k) + \bbE_{Y\sim\rho}\left[v_k(Y)\,\vert\, \calF_{k}\right],
\end{equation*}
where $\calF_{k}$ denotes the $\sigma$-algebra generated by the draws $\{Y_0,\dots,Y_{k-1}\}$.
Therefore, the design choices for the algorithm are a suitable sequence of model spaces $V_{m_k}$, the size of the memory $s_k$, and the step size $\tau_k$.
The procedure is summarized in \Cref{algo:general-algorithm}.

\begin{algorithm}[ht]
    \caption{Stochastic gradient with memory-based control variate (\galgnamens)}
    \label{algo:general-algorithm}
    \begin{algorithmic}[1]
        \State \textbf{Input:} spaces $\{ \mathbb{V}_{m_k} \}_{k \geq 0}$, step sizes $\{ \tau_k \}_{k\geq0}$, memory sizes $\{ s_k \}_{k\geq0}$, initial guess $u_0$, initial model $v_0$ and memory $\mathcal{S}_0$, number of iterations $K$.
        \State \textbf{Output:} Approximate solution $u_K$.
        \For{$k = 0, 1, \dots,K-1$}
        \State Sample $Y_k \sim \rho$ and compute $\nabla g(u_k, Y_k)$
        \State Fit the model $v_k \in \bbV_{m_k}$ to $\nabla g(u_k, \cdot)$ using the data in memory $\mathcal{S}_k$
        \State Compute $\bbE_{Y\sim\rho}\left[v_k(Y) \,\vert\, \calF_{k} \right]$
        \State Update $u_{k+1} = u_k - \tau_k \left( \nabla g(u_k,Y_k) - v_k(Y_k) + \bbE_{Y\sim\rho}\left[v_k(Y) \,\vert\, \calF_{k} \right] \right)$
        \State Update memory $\calS_k$
        \EndFor
    \end{algorithmic}
\end{algorithm}

\begin{rmrk}
    Let us remark that our method is closely related to SAGA \cite{defazio2014saga}.
    Consider an atomic measure $\rho = \frac{1}{s} \sum_{i=1}^{s} \delta_{y_i}$ with $\{y_i\}_{i=1}^s \subset \Gamma$, so that the functional to minimize is $J(u)=\frac{1}{s}\sum_{i=1}^s g_i(u)$ with $g_i(u):=g(u,y_i)$.
    SAGA exploits a variance reduction technique and requires introducing a memory term which stores all previously computed gradients in the sum and overwrites a term if the corresponding index is redrawn. Specifically, the SAGA algorithm reads
    \begin{equation*}
        u_{k+1} = u_k - \tau \widehat{\nabla J}_{\mathrm{SAGA}}(u_k),
        \quad \widehat{\nabla J}_{\mathrm{SAGA}}(u_k) := \nabla g_{i_k}(u_k) - \nabla g_{i_k}(z_{i_k}^k) + \frac{1}{s} \sum_{j=1}^{s} \nabla g_j(z_{j}^k),
    \end{equation*}
    where, at each iteration $k$, an index $i_k \in  \{ 1,\dots , s\}$  is selected uniformly at random on $\{ 1, \dots , s\}$, and only the memory term at the sampled position $i_k$ is updated, so that $z^{k+1}_j = u_k$ if $j = i_k$ and $z^{k+1}_j = z^{k}_j$ otherwise.
    The major improvement of SAGA with respect to SGD is its exponential convergence $\gamma^k$ for some $\gamma \in (0,1)$, provided that the objective functionals $g_i$ are all strongly convex and with Lipschitz continuous
    gradients with uniform bounds in $i$ \cite{defazio2014saga}.
    Note that SAGA can indeed be recovered from \Cref{algo:general-algorithm} considering a model space $\mathbb{V} = \{v :\{1,\dots,s\} \to U\}$ and a memory size equal to $s$.
    The limitation of SAGA, which is overcome by our method, is that it can only be applied to objective functionals which are sums of a finite number of terms and its performance and memory occupation degrades for $s$ large.
\end{rmrk}

\subsection{Review of optimal weighted least-squares} \label{sec:least-squares-revision}
{\color{black} Our \galgname method relies on constructing a model of the gradient of the functional $g(u_k, \cdot)$ using the data in memory $\calS_k$, and we propose to do so via least-squares fitting with inexact data and optimal sampling.}
In this section we hence review this technique following \cite{cohen2016optimal}.
For this purpose, fix $u \in U$ and denote $\psi := \nabla g(u, \cdot): \Gamma \to U$ to ease notation.
Then, let us fix some space $\mathbb{V}_m$ of dimension $m$ of functions defined everywhere on $\Gamma$ and assume that, for any $y\in \Gamma$, there exists $v \in \mathbb{V}_m$ such that $v(y) \neq 0$.
Let $\mu$ be a probability measure on $\Gamma$, such that $\rho$ is absolutely continuous with respect to $\mu$, and consider an i.i.d. random sample $\{ Y_j \}_{j=1}^s$ from $\mu$.
Define further the dataset of sample-noisy observation pairs $\calS := \{(Y_i,{\color{black} \hat{\psi}_i})\}_{i=1}^s$, where
\begin{equation*}
    {\color{black} \hat{\psi}_i} := \psi(Y_i) + {\color{black} \eta_i},
    \quad i=1,\dots,s,
\end{equation*}
{\color{black} and $\eta_i = h(Y_i)$ is an additive noise term modeled by a noise function $h \in L^2_\rho(\Gamma; U)$}.
Then, the weighted least-squares estimator based on data $\calS$ is given by
\begin{equation*}
    \hat{\Pi}^{\mathbb{V}_m}_{\calS} [\psi] = \operatornamewithlimits{argmin}\limits_{v \in \mathbb{V}_m} \frac{1}{s} \sum_{i=1}^s w(Y_i) \norm{v(Y_i)- {\color{black} \hat{\psi}_i}}_U^2,
\end{equation*}
where $w := \frac{d\rho}{d\mu}$.
The solution to this least-squares problem is provided by the well-known normal equations, which are described next.
We first introduce the empirical semi-norm for $U$-valued functions $v: \Gamma \to U$
\begin{equation*}
    \norm{v}_{s,U}^2 = \frac{1}{s} \sum_{i=1}^s w(Y_i) \norm{v(Y_i)}_U^2,
\end{equation*}
which is a Monte Carlo approximation of the $L^2_\rho(\Gamma,U)$-norm.
Take now an orthonormal basis $\{ \phi_1,\dots,\phi_m \}$ of $\mathbb{V}_m$ and define the pseudo-Vandermonde matrix $V \in U^{s \times m}$, the weight matrix $W \in \R^{s \times s}$, and the noise-corrupted data vector $\hat{\Psi} \in U^{s}$ as
\begin{equation*}
    V_{ij} := \phi_j(Y_i),
    \quad
    W := \mathrm{diag}(w(Y_1),\dots,w(Y_s)),
    \quad
    \hat{\Psi}_i := {\color{black} \hat{\psi}_i}.
\end{equation*}
as well as the Gram matrix $G \in \R^{m \times m}$ and the vector $J \in \R^m$
\begin{equation*}
    G := \frac{1}{s} \inp{V}{W V}_{U},
    \quad J := \frac{1}{s} \inp{V}{W \hat{\Psi}}_{U},
\end{equation*}
where for $A \in U^{s \times m}$ and $B \in U^{s \times p}$ we defined the inner product
\begin{equation*}
    \left(\inp{A}{B}_{U}\right)_{ij} := \sum_{k=1}^{s} \inp{A_{ki}}{B_{kj}}_U.
\end{equation*}
Then, the normal equations read
\begin{equation*}
    G c = J
\end{equation*}
and the least-squares estimator reads $\hat{\Pi}^{\bbV_m}_{\calS} [\psi] = \sum_{i=1}^m c_i \phi_i$.
It turns out that, whenever $G$ is well conditioned, the reconstruction error is quasi-optimal, namely it is bounded by the best approximation error
\begin{equation} \label{eq:best-approx}
    e_m(\psi) = \norm{\psi-\Pi^{\mathbb{V}_m} [\psi]}_{L^2_\rho(\Gamma;U)}
    \quad \mathrm{with} \quad
    \Pi^{\mathbb{V}_m} [\psi] = \argmin_{v \in \bbV_m} \norm{\psi-v}_{L^2_\rho(\Gamma;U)},
\end{equation}
plus noise-dependent terms, up to a constant term.
This motivates the definition of the conditioned weighted least-squares estimator
\begin{equation} \label{eq:cond-least-squares-estim}
    \tilde{\Pi}^{\mathbb{V}_m}_{\calS} [\psi] :=
    \begin{cases}
        \hat{\Pi}^{\mathbb{V}_m}_{\calS} [\psi], \quad & \text{if } \norm{G-I}_2 \leq \delta < 1 \\
        0, \quad                                       & \text{otherwise},
    \end{cases}
\end{equation}
where $\norm{\cdot}_2$ denotes the matrix spectral norm.
The idea now is to choose the sample size $s$ big enough so that $\norm{G-I} < \delta$ with high probability.
Let us introduce the function
\begin{equation*}
    \Gamma \ni y \mapsto k_{m,w}(y) = \sum_{j=1}^m w(y) \norm{\phi_j(y)}_U^2,
\end{equation*}
whose reciprocal is known as the Christoffel function \cite{Neval1986GezaFO}, and define
\begin{equation*}
    K_{m,w}:= \norm{k_{m,w}}_{L_\rho^{\infty}(\Gamma)},
\end{equation*}
which, assuming $\rho \ll \mu$, trivially satisfies $K_{m,w} \geq m$ since $\int_\Gamma k_{m,w} d \mu = m$.
We can now state the following theorem.
\begin{thrm}{(From \cite[Theorem 4.1]{cohen2016optimal})}
    Let $\delta = \frac{1}{2}$ in \eqref{eq:cond-least-squares-estim}.
    For any $r>0$, if $m$ and $s$ are such that
    \begin{equation} \label{eq:samp-ineq}
        K_{m,w} \leq \kappa \frac{s}{\log(s)}
        \quad \text{with} \quad
        \kappa =\frac{1- \log(2)}{2+2r},
    \end{equation}
    where $w = \frac{d \rho}{d \mu}$, then the conditioned weighted least-squares estimator satisfies
    \begin{equation*}
        \Exp{\norm{\psi-\tilde{\Pi}^{\mathbb{V}_m}_{\calS} [\psi]}_{L^2_\rho(\Gamma;U)}^2} \leq (1 + 2\epsilon(s)) e_m(\psi)^2 + (8+2 \epsilon(s)) {\color{black} \norm{h}_{L^2_\rho(\Gamma;U)}^2} + 2 \norm{\psi}_{L^2_\rho(\Gamma;U)}^2 s^{-r},
    \end{equation*}
    with $\epsilon(s):= \frac{4 \kappa}{\log(s)}$.
    \label{thm:theorem-approx-lstsq}
\end{thrm}

From the result of the previous theorem, we have interest in choosing the sampling distribution $\mu$ to minimize $K_{m,w}$.
In this sense, the optimal sampling measure $\mu^*$ is given by
\begin{equation} \label{eq:opt-samp-meas}
    d \mu^\ast (y) := \frac{1}{m} \sum_{j=1}^m \norm{\phi_j(y)}_U^2 d \rho (y),
\end{equation}
which yields the optimal $K_{m,w^*} = m$.
Observe further that, if $\mu$ is any sampling measure satisfying the quasi optimality relation
\begin{equation}\label{eq:samp-meas-quasi-opt}
    c_\mu \frac{d \mu}{d \rho} \leq \frac{d \mu^*}{d \rho} \leq C_\mu \frac{d \mu}{d \rho}
\end{equation}
for some constants $C_\mu\geq1\geq c_\mu>0$, then $K_{m,w} \leq C_\mu m$.

\begin{xmpl}[Polynomial spaces with uniform measure] \label{exmp:legendre}
    Let $\Gamma=[-1,1]$ with Lebesgue density $d\rho(y) = \frac{d y}{2}$, $U = \R$, and consider the space of polynomials of degree smaller than $m$, that is $\mathbb{V}_m = \mathbb{P}_{m-1}$.
    Then we may take $\phi_j = L_{j-1}$, where $L_j$ is the Legendre polynomial of degree $j$.
    When the sampling measure is taken as $d\mu = d\rho$, hence $w(y)=1$, it holds $K_{m,w} = \calO(m^2)$, $m$ being the dimension of $\mathbb{V}_m$, which implies that the sample size must be at least of order $m^2\log(m)$.
    This can be reduced to the optimal sample size of $\calO(m \log(m))$ by taking the optimal sampling measure $\frac{1}{m} \sum_{j=1}^m \phi_j(y)^2$.
    Furthermore, in this case, it holds that the arcsine distribution on $\Gamma$
    \begin{equation*}
        d \mu = \frac{d y}{\pi \sqrt{1-y^2}}
    \end{equation*}
    satisfies a quasi optimality relation of the type \eqref{eq:samp-meas-quasi-opt} with constants $c_\mu =1$, $C_\mu = \frac{4}{\pi}$, hence the sampling can be done from $\mu$, independently of the size $m$ of the polynomial space.
\end{xmpl}

\begin{xmpl}[Polynomial spaces with Gaussian measure] \label{exmp:hermite}
    If $\Gamma=\R$ with Gaussian density $d \rho(y) = \frac{1}{\sqrt{2 \pi}} e^{-\frac{y^2}{2}}dy$, we may take $\phi_j = H_{j-1}$ instead, $H_j$ being the normalized Hermite polynomial of degree $j$.
    When the sampling measure is taken as $d\mu = d\rho$, hence $w(y)=1$, $K_{m,w}$ is unbounded \cite{migliorati2013polynomial}.
    However, if we consider the optimal sampling measure $d \mu^*(y) = \frac{1}{m} \sum_{j=0}^{m-1} H_j(y)^2 d\rho(y)$, drawing $\calO(m \log(m))$ samples is sufficient to achieve quasi-optimal reconstruction error.
\end{xmpl}

\subsubsection{Least-squares in tensor-product spaces}
\label{sec:tensor-product-spaces}
Polynomial spaces like those introduced in \Cref{exmp:legendre,exmp:hermite} {\color{black} can be used} as building blocks in order to construct spaces of functions from $\Gamma \subseteq \R^d$ to $U$ by tensorization.
In this section, we will briefly detail how this can be done as it is relevant to the PDE-constrained optimization problems we consider in the numerical experiments, see \Cref{sec:numerics}.

First, given a subspace $\calV_m$ of $L^2_\rho(\Gamma;\R)$ of dimension $m$ we can extend the output space from $\R$ to $U$ by considering the tensor space $\mathbb{V}_m = \calV_m \otimes U$.
In this case, the dimension of $\bbV_m$ is $m \cdot \mathrm{dim}(U)$, which is infinite if $\mathrm{dim}(U) = \infty$.
Letting $\phi_1,\dots,\phi_m$ be an orthonormal basis of $\calV_m$, we can redefine the Vandermonde matrix $V\in \R^{s\times m}$ and the matrices $G \in \R^{m\times m}$ and $J \in \R^m \otimes U$ for the $\R$-valued orthonormal basis functions $\{\phi_j\}_{j=1}^m$ spanning $\calV_m$, so that the normal equations read
\begin{equation*}
    (G \otimes I) c = J, \quad c \in \R^m \otimes U,
\end{equation*}
where $I$ denotes the identity operator on $U$.

In applications related to PDE-constrained optimization, we could consider for instance $U = L^2(D)$ for some $n$-dimensional physical domain $D\subset \R^n$ or, in practice, a discretization of it (e.g., a finite element space).
Moreover, in this setting one often works with tensorized parameter domains of the form $\Gamma = \prod_{i=1}^d \Gamma_i \subseteq \mathbb{R}^d$,
where each $\Gamma_i \subseteq \mathbb{R}$ is an interval and $d \in \mathbb{N} \cup \{\infty\}$.
It is then common to equip $\Gamma$ with a product measure $\rho = \bigotimes_{i=1}^d \rho_i$, where $\rho_i$ is a Borel probability measure supported on $\Gamma_i$ with finite moments of all orders and such that polynomials are dense in $L^2_{\rho_i}(\Gamma_i)$ (e.g., uniform or Gaussian, as in \Cref{exmp:legendre,exmp:hermite}).
When $d > 1$, a central tool in the approximation of smooth functions $\psi : \Gamma \to U$ is the use of polynomial expansions in terms of an orthonormal tensorized basis.
To this end, we define the set of \emph{multi-indices}
\begin{equation*}
    \mathcal{F} :=
    \begin{cases}
        \mathbb{N}^d                                                                           & \text{if } d < \infty, \\
        \left\{ \nu \in \mathbb{N}^{\mathbb{N}} : |\operatorname{supp}(\nu)| < \infty \right\} & \text{if } d = \infty,
    \end{cases}
\end{equation*}
where $\operatorname{supp}(\nu) := \{ i \in \mathbb{N} : \nu_i \neq 0 \}$ denotes the (finite) support of $\nu$ in the infinite-dimensional case.
For each $i \in \mathbb{N}$, let $\{\phi_{j}^{(i)}\}_{j=0}^\infty$ be a univariate orthonormal polynomial basis in $L^2_{\rho_i}(\Gamma_i)$. Then, for any $\nu \in \mathcal{F}$, we define the multivariate orthonormal polynomial
\begin{equation*}
    \phi_\nu(y) := \prod_{i \in \operatorname{supp}(\nu)} \phi_{\nu_i}^{(i)}(y_i),
\end{equation*}
which defines a complete orthonormal basis in $L^2_\rho(\Gamma)$.
A multivariate function $\psi \in L^2_\rho(\Gamma; U)$ can thus be expanded as
\begin{equation*}
    \psi(y) = \sum_{\nu \in \mathcal{F}} \psi_\nu \phi_\nu(y),
\end{equation*}
with coefficients $\psi_\nu \in U$.
See \cite{cohen2018multivariate} for more details.
In practical approximations, the infinite expansion is truncated to a finite sum over a \emph{multi-index set} $\Lambda \subset \mathcal{F}$.
A widely used class of multi-index sets is that of \emph{downward closed} or \emph{lower sets}, defined by
\begin{equation*}
    \text{$\Lambda$ is downward closed} \quad \Longleftrightarrow \quad \nu \in \Lambda, \, \tilde{\nu} \leq \nu \Rightarrow \tilde{\nu} \in \Lambda,
\end{equation*}
where the inequality $\tilde{\nu} \leq \nu$ should be understood componentwise.
Downward closed sets are natural in sparse polynomial approximation, adaptive algorithms \cite{cohen2018multivariate}, and compressed sensing \cite{adcock2022sparse}, as they ensure closedness of polynomial spaces under differentiation and allow efficient construction via greedy or adaptive selection.
In the convergence analysis of \Cref{sec:conv-analyis} it will be important to link the supremum of the weight $w = \frac{d \rho}{d \mu}$ to the dimension $m$ of the approximation space $\bbV_m$.

\begin{prpstn} \label{prop:sup-weight}
    Let $\Lambda \subset \calF$ be a multi-index set such that $\abs{\Lambda} = m$.
    Define the weight
    \begin{equation*}
        w_\Lambda(y) := \frac{m}{\sum_{\nu \in \Lambda} \phi_\nu^2(y)}
    \end{equation*}
    and its supremum
    \begin{equation*}
        q_\Lambda := \sup_{y\in\Gamma} w_\Lambda(y).
    \end{equation*}
    Assume the constant function is in $\calV_\Lambda$, then $q_\Lambda \leq m$.
    Furthermore, when $d=1$, $\rho$ is Gaussian, and $\calV_\Lambda = \bbP_{m-1}$, there exists $C>0$ such that $q_\Lambda \leq C \sqrt{m}$.
\end{prpstn}
\begin{proof}
    See \Cref{sec:appendix-prop-sup-weight}.
\end{proof}

\subsection{Least-squares memory-based control variate}
\label{sec:ls-mb-cv}
In this section we present the \galgname algorithm where the model fitting is done via the optimal weighted least-squares method introduced in the previous section.
We call this algorithm stochastic gradient with least-squares control variates (\algnamens).
\subsubsection{Fixed approximation space}
We start by giving a first version of the algorithm where the approximation space does not depend on the iteration, that is $m_k = m$ for all $k$.
Owing to \Cref{thm:theorem-approx-lstsq}, this entails that, in order to have quasi-optimal reconstruction, it is sufficient to fix the sampling measure to $\mu$ such that $\rho \ll \mu$ and the memory size to some sufficiently large integer $s > 0$.
We further assume that computing $\bbE_{Y \sim \rho}\left[v(Y)\right]$ is cheap for any (deterministic) function $v \in \mathbb{V}_m$.
The idea of the algorithm is to use least-squares approximations using past gradient evaluations to create a control variate for the current gradient.
As we have seen in \Cref{sec:least-squares-revision}, to numerically approximate a continuous function by least-squares, one needs to solve a linear system which involves the Gram matrix $G_k$ and the noise-corrupted evaluations of the target function which are stored in the memory $\mathcal{S}_k = \{(Y_j,\nabla g(u_j,Y_j))\}_{j=k-s}^{k-1}$.
Hereafter, we denote $\calS_k(1)$ the collection $\{Y_j\}_{j=k-s}^{k-1}$ and $\calS_k(2)$ the collection $\{\nabla g(u_j,Y_j)\}_{j=k-s}^{k-1}$, the weighted least-squares approximation of the gradient at $u_k$ with data $\mathcal{S}_k$ as $\hat{\Pi}^{\mathbb{V}_m}_{{\mathcal{S}_k}}  [\nabla g]$, and, similarly, the conditioned weighted least-squares approximation as follows
\begin{equation*}
    \tilde{\Pi}^{\mathbb{V}_m}_{{\mathcal{S}_k}}  [\nabla g] =
    \begin{cases}
         & \hat{\Pi}^{\mathbb{V}_m}_{{\mathcal{S}_k}}  [\nabla g] \quad \text{if } \norm{G_k-I}_2 \leq \frac{1}{2}, \\
         & 0 \quad \text{ otherwise}.
    \end{cases}
\end{equation*}
We remark that $\tilde{\Pi}_{\calS_k}^{\bbV_m}$ aims at approximating $\nabla g(u_k,\cdot)$, however using inconsistent data since $\calS_k(2)$ contains gradient evaluations for previous designs $u_j \neq u_k$.

Recall that, in practice, one does not solve the normal equations in order to obtain the solution to a least-squares problem.
Indeed, a more stable method makes use of the QR factorization of the weighted Vandermonde matrix $Q_k R_k = \sqrt{W_k} V_k$, which we assume to be implemented in a routine to solve least-squares problems $\texttt{SolveLSTSQ}(Q,R,\sqrt{W_k} \hat{\Psi})$, where $\hat{\Psi}$ denotes a noise-corrupted data vector, which returns the coefficient vector $c \in U^m$.
Furthermore, in the context of our algorithm, we must solve a least-squares problem at each iteration, hence we assume that at each iteration $k$ the QR factorization at the previous step $Q_{k-1} R_{k-1} = \sqrt{W_{k-1}} V_{k-1}$ is already computed and available (i.e., stored in memory).
Since $\sqrt{W_{k-1}} V_{k-1}$ and $\sqrt{W_k} V_k$ differ by two rows -- more precisely, {\color{black} we obtain $V_k$ by removing the first row of $V_{k-1}$} and appending a new row at the bottom -- we can use standard routines to cheaply update $Q_{k-1}$ and $R_{k-1}$ to obtain $Q_k$ and $R_k$, hence avoiding recomputing the QR factorization from scratch.
We then consider a generic method $\texttt{UpdateQR}(A,B,Q,R)$ taking matrices $A,B,Q,R$ such that $QR=A$ which computes efficiently the QR decomposition of $B$ if it differs from $A$ by a few rank-$1$ updates.

With this in mind, the different steps of the method at iteration $k$ are as follows.
\begin{enumerate}[label=(\alph*)]
    \item Sample the random variable $Y_k$ from a chosen distribution $\mu$.
    \item Approximate $\nabla g(u_k,\cdot)$ in the space $\mathbb{V}_m$ using memory $\calS_k$.
          The approximation is done via the weighted discrete least-squares method described in \Cref{sec:least-squares-revision}, i.e. solving the least-squares problem
          \begin{equation*}
              c_k = \texttt{SolveLSTSQ}(Q_k,R_k,\sqrt{W_k} \hat{\Psi}_k),
          \end{equation*}
          where $Q_k$, $R_k$, and $\hat{\Psi}_k = \calS_k(2)$ are respectively the QR factors of the weighted Vandermonde matrix $\sqrt{W_k} V_k$, and the noise-corrupted data vector $(\nabla g (u_j,Y_j))_{j=k-s}^{k-1} \in U^s$ at iteration $k$.
          This allows us to construct the approximation of $\nabla g(u_k,\cdot)$
          \begin{equation*}
              \tilde{\Pi}^{\mathbb{V}_m}_{{\mathcal{S}_k}}  [\nabla g] =
              \begin{cases}
                   & \sum_{i=1}^{m} (c_k)_i \phi_i \quad \text{if } \norm{\frac{1}{s} \inp{V_k}{W_k V_k}_{U}-I}_2 \leq \frac{1}{2}, \\
                   & 0 \quad \text{ otherwise}.
              \end{cases}
          \end{equation*}
          Finally, we compute
          \begin{equation*}
              \bbE_{Y\sim\rho}\left[\tilde{\Pi}^{\mathbb{V}_m}_{{\mathcal{S}_k}}[\nabla g](Y)\,\vert\,\calF_k\right],
          \end{equation*}
          where $\calF_k = \sigma\{Y_0,\cdots, Y_{k-1} \}$ is the $\sigma$-algebra generated by the random variables $Y_0,\cdots,Y_{k-1}$.
    \item Update the parameter $u_k$ using the constructed approximation of the gradient
          \begin{equation*}
              u_{k+1}= u_k - \tau_k \left( w(Y_k)\left( \nabla g(u_k,Y_{k}) - \tilde{\Pi}^{\mathbb{V}_m}_{{\mathcal{S}_k}}[\nabla g](Y_{k})\right) + \bbE_{Y\sim\rho}\left[\tilde{\Pi}^{\mathbb{V}_m}_{{\mathcal{S}_k}}[\nabla g](Y)\,\middle\vert\,\calF_k\right] \right),
          \end{equation*}
          where $\tau_k$ is the step size.
    \item \label{step:update-S-V-W} Update the data points by removing the oldest data point, $(Y_{k-s}, \nabla g(u_{k-s},Y_{k-s}))$, and adding the new data point $(Y_k,\nabla g(u_k,Y_k))$.
          Update the Vandermonde and weight matrices by deleting the first row and appending the new row corresponding to $Y_k$ at the end.
          Finally, update the QR factorization of the weighted Vandermonde matrix $\sqrt{W_{k+1}} V_{k+1}$ via
          \begin{equation*}
              Q_{k+1}, R_{k+1} = \texttt{UpdateQR}(\sqrt{W_{k}} V_{k},\sqrt{W_{k+1}} V_{{k+1}},Q_{k},R_{k}).
          \end{equation*}
\end{enumerate}
The full procedure is given in \Cref{alg:algo-fixed}.

\begin{algorithm}[ht]
    \caption{\algnamens: least-squares on fixed approximation space}
    \label{alg:algo-fixed}
    \begin{algorithmic}[1]
        \State \textbf{Input:} approximation space $\bbV_m$, sampling measure $\mu$, step sizes $\{\tau_k\}_{k\geq 0}$, $u_{-s},\dots,u_0$, initial memory $\mathcal{S}_0$, QR factors $Q_0$, $R_0$ of $\sqrt{W_0}V_0$, maximum number of iterations $K$
        \State \textbf{Output:} Approximated solution $u_K$
        \For{$k = 0, 1, \dots, K-1$}
        \State Sample $Y_k \sim \mu$ and compute $\nabla g(u_k, Y_k)$
        \State Compute the weighted least-squares approximation by solving
        \begin{equation*}
            c_k = \texttt{SolveLSTSQ}(Q_k,R_k,\sqrt{W_k} \hat{\Psi}_k)
            \quad \mathrm{with} \quad \hat{\Psi}_k = \calS_k(2)
        \end{equation*}
        \State Create the least-squares approximation $\tilde{\Pi}^{\mathbb{V}_m}_{{\mathcal{S}_k}} [\nabla g]$ using $c_k$
        \State Update $u_{k+1} = u_k - \tau_k \left( w(Y_k)\left(\nabla g(u_k,Y_{k}) - \tilde{\Pi}^{\mathbb{V}_m}_{{\mathcal{S}_k}} [\nabla g](Y_{k})\right) + \mathbb{E}_{Y\sim \rho}\left[\tilde{\Pi}^{\mathbb{V}_m}_{{\mathcal{S}_k}} [\nabla g](Y) \,\middle\vert\, \calF_k \right] \right)$
        \State Update $\mathcal{S}_k$, $V_k$, $W_k$ as explained in step \ref{step:update-S-V-W} above and the QR factors as
        \begin{equation*}
            Q_{k+1}, R_{k+1} = \texttt{UpdateQR}(\sqrt{W_{k}} V_{k},\sqrt{W_{k+1}} V_{{k+1}},Q_{k},R_{k}).
        \end{equation*}
        \EndFor
    \end{algorithmic}
\end{algorithm}

\begin{rmrk}\label{rem:var-no-to-zero}
    Let us point out what we expect at this level in terms of convergence and complexity from this approach.
    Recall that the asymptotic variance of $\widehat{\nabla J}_{\mathrm{SGD}}$ is of the order $\norm{\nabla g (u^\ast,\cdot)}_{L^2_\rho}$.
    On the other hand, provided everything works correctly, we expect the asymptotic variance of $\widehat{\nabla J}_{\algname}$ to scale like $\norm{(I-\Pi^{\bbV_m}) \nabla g (u^\ast,\cdot)}_{L^2_\rho}$, which does not vanish unless $\nabla g (u^\ast,\cdot) \in \bbV_m$.
    Hence, we expect that the asymptotic complexity of \Cref{alg:algo-fixed} matches the one exhibited by SGD, but hopefully with a much smaller constant.
\end{rmrk}

\begin{rmrk}
    In \Cref{alg:algo-fixed} we need to provide initial iterates $u_{-s},\dots,u_0$ and an initial memory $\mathcal{S}_0$.
    The easiest way to do so is to choose $u_{0} \in U$, set $u_k = u_0$ for all $k=-s,\dots,-1$, sample $\{Y_k\}_{k=-s}^0 \stackrel{\mathrm{iid}}{\sim}\mu$, and compute gradients $\nabla g (u_0,Y_k)$, $k=-s,\dots,0$.
    Another simple way to achieve this is to choose $u_{-s} \in U$ and run $s$ iterations of a gradient-based optimization method like SGD, keeping all iterates and gradient evaluations in memory.
\end{rmrk}

\subsubsection{Variable approximation space}
\label{sec:var-space-algo}
The considerations made in \Cref{rem:var-no-to-zero} shed light on how to modify \Cref{alg:algo-fixed} in order to beat the SGD asymptotic complexity.
Indeed, we shall decrease the variance of $\widehat{\nabla J}_{\algname}$ progressively through the iterations sending $m$ to infinity, that is choosing an index $m_k$ such that $m_k \to \infty$.
A variable approximation space also addresses the following issue: if we choose $\bbV_m$ large in \Cref{alg:algo-fixed} in order to keep $\norm{(I-\Pi^{\bbV_m}) \nabla g (u^\ast,\cdot)}_{L^2_\rho}$ sufficiently small, the least-squares fit of the control variate during the first iterations, where $u_k$ is rapidly changing hence the norm of the noise term in the least-squares is high, may incur a large approximation error making the control variate less effective.
Hence, it would be natural to regularize the approximation at the beginning of the optimization process keeping the approximation space size small, and progressively enlarging it as $u_k$ converges to $u^\ast$, until a desired precision is reached.
Note that, increasing the size of the approximation space entails, owing to \Cref{thm:theorem-approx-lstsq}, increasing the number of samples, that is the size of the memory needed for the least-squares approximation.

Let $\{\bbV_{m_p}\}_{p \geq 0}$ be a sequence of nested approximation spaces, $\mathrm{dim}(\bbV_{m_p}) = m_p$, and let $\bar{s}_{p}$ be the minimum sample size needed to perform a stable approximation of the gradient in $\mathbb{V}_{m_p}$ by the weighted least-squares method with sampling measure $\mu_p$, as prescribed by \Cref{thm:theorem-approx-lstsq}.
We choose the memory to have exactly this size.
Let further $\sigma : \mathbb{N}_0 \to \mathbb{N}$ be a scheduling function such that $\sigma(0)=0$ and $\sigma(p+1) - \sigma(p)$ is the number of iterations to perform with approximation in $\bbV_{m_p}$ before switching to $\bbV_{m_{p+1}}$.
Clearly, we need $\sigma(p) + \abs{\calS_0} \geq \bar{s}_p$, that is the total available data when we switch to $\bbV_{m_p}$ must be enough to perform stable approximation.
With such scheduling function at iteration $\sigma(p+1)$ we switch from $\bbV_{m_p}$ to $\bbV_{m_{p+1}}$ hence the memory requirement switches from $\abs{\calS_{\sigma(p+1)-1}} = \bar{s}_{p}$ to $\abs{\calS_{\sigma(p+1)}} = \bar{s}_{p+1}$.
Then, $(\bar{s}_{p+1} - \bar{s}_p)$ iterations prior to iteration $\sigma(p+1)$, we must begin storing new observations.
We formalize this by introducing an extended memory $\calS^e_k$ which contains these extra observations.
More precisely, the size of the extended memory is defined as follows
\begin{equation*}
    \abs{\calS^e_k} =
    \begin{cases}
        \bar{s}_p   & \text{if } \sigma(p) \leq k < \sigma(p+1) - (\bar{s}_{p+1} - \bar{s}_p),   \\
        s_{k-1} + 1 & \text{if } \sigma(p+1) - (\bar{s}_{p+1} - \bar{s}_p) \leq k < \sigma(p+1).
    \end{cases}
\end{equation*}
Then, setting $\calS_{\sigma(p)} = \calS^e_{\sigma(p)}$ ensures that the size of the memory satisfies
\begin{equation*}
    s_k = \abs{\calS_k} = \bar{s}_p, \quad \sigma(p) \leq k < \sigma(p+1).
\end{equation*}

We henceforth assume that either $\mu_p = \mu$ is some fixed probability measure $\mu$,
or $\mu_p$ is the optimal sampling measure \eqref{eq:opt-samp-meas} associated to the approximation space $\bbV_{m_p}$.
In the latter case we must make sure that, whenever we switch from $\mu_p$ to $\mu_{p+1}$, the set $\{Y_{i}\}_{i=\sigma(p+1)-\bar{s}_{p+1}}^{\sigma(p+1)-1}$ is an independent and identically distributed sample drawn from $\mu_{p+1}$.
To this end, we propose a resampling procedure inspired by optimal sampling techniques for nested approximation spaces \cite{arras2019sequential}.
Indeed, since we assume that $\bbV_{m_p} \subset \bbV_{m_{p+1}}$, the following relationship between optimal sampling measures of two consecutive approximation spaces holds
\begin{equation*}
    d \mu_{p+1} = \frac{1}{m_{p+1}} \sum_{j=1}^{m_{p+1}} \norm{\phi_j}^2_U d \rho = \frac{m_p}{m_{p+1}} \mu_p + \frac{1}{m_{p+1}} \sum_{j=m_p+1}^{m_{p+1}} \norm{\phi_j}^2_U d \rho,
\end{equation*}
that is $\mu_{p+1}$ is a mixture of the old optimal probability measure $\mu_p$ and the measure $\hat{\mu}_{p+1}$ defined as
\begin{equation*}
    d \hat{\mu}_{p+1} := \frac{1}{m_{p+1}-m_p}\sum_{j=m_p+1}^{m_{p+1}} \norm{\phi_j}^2_U d \rho.
\end{equation*}
Hence, assuming $\{Y_{j}\}_{j=\sigma(p+1)-\bar{s}_{p+1}}^{\sigma(p+1)-1} \sim \mu_p$, we shall add a resampling step of the extended memory at iteration $\sigma(p+1)$, where, for each $j = \sigma(p+1)-\bar{s}_{p+1},\dots,\sigma(p+1)-1$, we draw an independent Bernoulli variable $b_j \sim \mathrm{Bernoulli}(\frac{m_p}{m_{p+1}})$ and, if $b_j = 0$, sample $\tilde{Y}_{j} \sim \hat{\mu}_{p+1}$ and compute $\nabla g(u_j, \tilde{Y}_{j})$ to substitute the old data pair $(Y_{j}, \nabla g(u_j, Y_{j}))$ with the new pair $(\tilde{Y}_{j}, \nabla g(u_j, \tilde{Y}_{j}))$.
In other words, we set
\begin{equation*}
    \calS_{\sigma(p+1)} = \calR_p(\calS^e_{\sigma(p+1)}),
\end{equation*}
where the resampling operator $\calR_p(\calS)$ is a random operator defined by applying \Cref{algo:resampling} with input $\calS, \mu_p,\hat{\mu}_{p+1},m_p$, and $m_{p+1}$.
By construction, then, at any iteration $k \geq 0$ the points in memory at which the gradient is evaluated are sampled according to $\mu_{\sigma^{-1}(k)}$, where we use the generalized inverse
\begin{equation*}
    \sigma^{-1}(k) := \inf \{p : \sigma(p+1) > k\}, \quad k\geq0.
\end{equation*}
We shall also redefine the filtration $\calF_k$ to be the filtration generated by any draw and re-draw up to iteration $k$.

\begin{algorithm}[ht]
    \caption{Resampling routine}
    \label{algo:resampling}
    \begin{algorithmic}[1]
        \State \textbf{Input:} set (memory) $\calS = \{(Y_j, \nabla g(u_j, Y_j))\}_{j=1}^{\abs{\calS}}$, probability measures $\mu_0$, $\hat{\mu}_1$, approximation-space dimensions $m_0$,$m_1$.
        \State \textbf{Output:} Resampled set $\tilde{\calS}$.
        \State Initialize $\tilde{\calS} \gets \varnothing$
        \For{$j = 1, 2, \dots,\abs{\calS}$}
        \State Sample $b_j \sim \mathrm{Bernoulli}(\frac{m_0}{m_{1}})$
        \If{$b_j=1$}
        \State Update memory $\tilde{\calS} \gets \tilde{\calS} \cup \{(Y_j, {\color{black} \nabla g(u_j, Y_j))}\}$
        \Else
        \State Sample $\tilde{Y}_j \sim \hat{\mu}_{1}$ and evaluate gradient $\nabla g(u_j, \tilde{Y}_j)$
        \State Update memory $\tilde{\calS} \gets \tilde{\calS} \cup \{(\tilde{Y}_j, \nabla g(u_j, \tilde{Y}_j))\}$
        \EndIf
        \EndFor
    \end{algorithmic}
\end{algorithm}

Finally, notice that, at the implementation level, switching from $\bbV_{m_p}$ to $\bbV_{m_{p+1}}$, from $\mu_p$ to $\mu_{p+1}$, and setting $\calS_{\sigma(p+1)} = \calR_p(\calS^e_{\sigma(p+1)})$ entails appending $m_{p+1} - m_{p}$ new columns and $\bar{s}_{p+1} - \bar{s}_p$ new rows to the Vandermonde matrix $V_k$, and reconstructing the weight matrix (all elements must be recomputed since they depend on the sampling measure).
We assume that all of these operations are handled by the routine $\texttt{UpdateQR}$.
The pseudocode of the variable approximation space algorithm is given in \Cref{algo:algo-variable}.

\begin{algorithm}[ht]
    \caption{\algnamens: variable approximation space}
    \label{algo:algo-variable}
    \begin{algorithmic}[1]
        \State \textbf{Input:} approximation spaces $\{\bbV_{m_p}\}_{p\geq0}$, sampling measures $\{\mu_p\}_{p\geq0}$, step sizes $\{\tau_k\}_{k\geq0}$, scheduling function $\sigma$, initial iterates $u_{-s_0},\dots,u_0$, initial memory $\mathcal{S}_0$, QR factors of $Q_0$, $R_0$ of $\sqrt{W_0}V_0$, maximum number of iterations $K$.
        \State \textbf{Output:} Approximated solution $u_K$
        \State Initialize iteration counter $k=0$
        \For{$p = 0, 1, \dots, \sigma^{-1}(K)$}
        \For{$j = 1, \dots, \sigma(p+1)-\sigma(p)$}
        \State Sample $Y_k \sim \mu_p$ and compute $\nabla g(u_k, Y_k)$
        \State Compute the weighted least-squares approximation by solving
        \begin{equation*}
            c_k = \texttt{SolveLSTSQ}(Q_k,R_k,\sqrt{W_k} \hat{\Psi}_k)
            \quad \mathrm{where} \quad \hat{\Psi}_k = \calS_k(2)
        \end{equation*}
        \State Create the least-squares approximation $\tilde{\Pi}^{\mathbb{V}_{m_p}}_{{\mathcal{S}_k}} \nabla g(u_k,\cdot)$ using $c_k$
        \State Update $u_{k+1} = u_k - \tau_k \left( w_k(Y_k) \left(\nabla g(u_k,Y_{k}) - \tilde{\Pi}^{\mathbb{V}_{m_p}}_{{\mathcal{S}_k}} [\nabla g](Y_{k})\right) + \bbE_{Y\sim\rho}\left[\tilde{\Pi}^{\mathbb{V}_{m_p}}_{{\mathcal{S}_k}} [\nabla g] (Y) \,\middle\vert\, \calF_k\right]  \right)$
        \State Update $\mathcal{S}_k$, $V_k$, $W_k$ {\color{black} as explained in step \ref{step:update-S-V-W} on \cpageref{step:update-S-V-W}} and the QR factors as
        \begin{equation*}
            Q_{k+1}, R_{k+1} = \texttt{UpdateQR}(\sqrt{W_{k}} V_{k},\sqrt{W_{k+1}} V_{{k+1}},Q_{k},R_{k}).
        \end{equation*}
        \If{$j < \sigma(p+1) - \sigma(p) - (\bar{s}_{p+1} - \bar{s}_p)$}
        \State Set $\mathcal{S}^e_k = \mathcal{S}_k$
        \Else
        \State Update extended memory adding data point $(Y_k, \nabla g (u_k,Y_{k}))$
        \EndIf
        \State Update iteration counter $k = k+1$
        \EndFor
        \State Resample extended memory $\tilde{\calS}^e_{k} = \calR_p(\calS^e_{k})$
        \State Set $\calS_k = \tilde{\calS}^e_{k}$
        \State {Update the weight matrix $W_k$ recomputing all weights according to $\mu_{p+1}$, the Vandermonde matrix appending $m_{p+1} - m_{p}$ new columns and $\bar{s}_{p+1} - \bar{s}_p$ new rows as explained in \Cref{sec:var-space-algo}, and update the QR decomposition as
            \begin{equation*}
                Q_{k}, R_{k} = \texttt{UpdateQR}(\sqrt{W_{k-1}} V_{k-1},\sqrt{W_{k}} V_{{k}},Q_{k-1},R_{k-1}).
            \end{equation*}}
        \EndFor
    \end{algorithmic}
\end{algorithm}

\section{Convergence analysis} \label{sec:conv-analyis}
In this section we present the convergence analysis of the \Cref{alg:algo-fixed,algo:algo-variable} when using weighted least-squares approximations.
For ease of exposition, we will abuse notation and express all parameters involved as functions of the iteration number $k$ (e.g., the approximation space size at iteration $k$, namely $m_{\sigma^{-1}(k)}$, will be simply denoted by $m_k$).
Furthermore, we also abuse notation and denote $Y_j$ (as opposed to $\tilde{Y}_j$) any re-draw in the memory.

\subsection{Error recursion}
The first goal is to derive a recurrence relation for the error of \Cref{algo:algo-variable}.
Notice that this will apply to \Cref{alg:algo-fixed} as a particular case where the approximation space and sampling measures are fixed.
Let us denote $\tilde{\Pi}^{\mathbb{V}_m}_{\mathcal{S}_k}  [\nabla g]$ the conditioned weighted least-squares approximation.
Denote the $k$-th iteration of the variable approximation space algorithm as
\begin{equation*}
    u_{k+1} = u_k - \tau_k \widehat{\nabla J}(u_{k},Y_{k};\calS_k),
\end{equation*}
where
\begin{equation}
    \widehat{\nabla J}(u_{k},Y_{k};\calS_k)
    = w_k(Y_k)\left(\nabla g(u_k,Y_k) -  \tilde{\Pi}^{\mathbb{V}_{m_k}}_{\mathcal{S}_k}  [\nabla g](Y_k) \right)  + \mathbb{E}_{Y \sim \rho}\left[\tilde{\Pi}^{\mathbb{V}_{m_k}}_{\mathcal{S}_k} [\nabla g](Y) \,\middle\vert\,\calF_k\right].
    \label{eq:hat-grad-J}
\end{equation}
{\color{black} Let us recall that $Y_k \sim \mu_k$ throughout this section} and that the least-squares projection uses the memory $\calS_k = \{ (u_j, \nabla g(u_j, Y_j)) \}_{j=k-s_k}^{k-1}$.
Moreover, for the purpose of the analysis, let us introduce the idealized memory at step $k$ containing the observations of the gradient at the optimum, that is $\mathcal{S}^*_k := \{ (u^*, \nabla g(u^*, Y_j)) \}_{j=k-s_k}^{k-1}$, so that
\begin{equation*}
    \widehat{\nabla J}(u^*,Y_{k};\calS^*_k)
    = w_k(Y_k)\left(\nabla g(u^*,Y_k) - \tilde{\Pi}^{\mathbb{V}_{m_k}}_{\mathcal{S}^*_k}  [\nabla g](Y_k)\right)  + \bbE_{Y\sim \rho}\left[ \tilde{\Pi}^{\mathbb{V}_{m_k}}_{\mathcal{S}^*_k} [\nabla g](Y) \,\middle\vert\,\calF_k\right].
\end{equation*}


\begin{lmm} \label{lem:first}
    {\color{black} Under \Cref{asmp:str-cvx}}, let $u_k$ be the $k$-th iterate of \Cref{algo:algo-variable}, then
    \begin{equation*}
        \begin{split}
            \Exp{\norm{u_{k+1}-u^*}_U^2\,\middle\vert\, \calF_k} \leq & (1-2 \alpha \tau_k )\norm{u_k - u^*}_U^2                                                                                                          \\
                                                                   & +2\tau_k^2 \Exp{\norm{ \widehat{\nabla J}(u_{k},Y_{k};\calS_k)-     \widehat{\nabla J}(u^*,Y_{k};\calS^*_k)}_U^2\,\middle\vert\, \calF_k}         \\
                                                                   & + 2 \tau_k^2\Exp{w_k(Y_k)^2 \norm{\nabla g(u^*,Y_k) -  \tilde{\Pi}^{\bbV_{m_k}}_{\mathcal{S}^*_k}  [\nabla g](Y_k)}_U^2\,\middle\vert\, \calF_k}.
        \end{split}
    \end{equation*}
\end{lmm}

\begin{proof}
    Recall that $u^*=\mathrm{argmin}_u J(u)$.
    Then, $u^*$ is a first order critical point for $J$, that is
    \begin{equation}
        \nabla J(u^*) = \mathbb{E}[w_k(Y_k)\nabla g(u^*,Y_k)]= 0,
        \label{eq:exp-sol}
    \end{equation}
    where $Y_k \sim \mu_k$.
    Using \eqref{eq:exp-sol}, one gets
    \begin{equation*}
        u_{k+1}-u^* = u_k - u^* - \tau_k \left( \widehat{\nabla J}(u_{k},Y_{k};\calS_k)- \nabla J(u^*)\right).
    \end{equation*}
    Taking the squared $U$-norm of this expression and expanding the right-hand side yields
    \begin{equation*}
        \begin{split}
            \norm{u_{k+1}-u^*}_U^2 = & \norm{u_k - u^*}_U^2                                                                          \\
                                     & - 2\tau_k \langle u_k -u^* , \widehat{\nabla J}(u_{k},Y_{k};\calS_k)- \nabla J(u^*) \rangle_U \\
                                     & + \tau_k^2\norm{\widehat{\nabla J}(u_{k},Y_{k};\calS_k)- \nabla J(u^*)}_U^2.
        \end{split}
    \end{equation*}
    Recall now that $\calF_k$ denotes the $\sigma$-algebra generated by the random variables $\{ Y_i \}_{i=0}^{k-1}$ (including eventual resampling).
    Then,
    \begin{equation*}
        \begin{split}
            \Exp{\norm{u_{k+1}-u^*}_U^2\,\middle\vert\, \calF_k} = & \norm{u_k - u^*}_U^2                                                                                              \\
                                                                   & - 2\tau_k \inp{u_k -u^*}{\Exp{\widehat{\nabla J}(u_{k},Y_{k};\calS_k)\,\middle\vert\, \calF_k} - \nabla J(u^*)}_U \\
                                                                   & + \tau_k^2 \Exp{\norm{  \widehat{\nabla J}(u_{k},Y_{k};\calS_k)- \nabla J(u^*)}_U^2\,\middle\vert\, \calF_k}.
        \end{split}
    \end{equation*}
    Noticing that, conditioned on $\calF_k$, $\widehat{\nabla J}(u_{k},Y_{k};\calS_k)$ is an unbiased estimator of $\nabla J(u_k)$, indeed,
    \begin{equation*}
        \begin{split}
            \Exp{\widehat{\nabla J}(u_{k},Y_{k};\calS_k)\,\vert\, \calF_k}
             & = \Exp{w_k(Y_k) \left(\nabla g(u_k,Y_k) - \tilde{\Pi}^{\bbV_{m_k}}_{\mathcal{S}_k}  [\nabla g](Y_k)\right)  + \Exp{w_k(Y_k)\tilde{\Pi}^{\bbV_{m_k}}_{\mathcal{S}_k} [\nabla g](Y_k) \,\middle\vert\,\calF_k}\,\middle\vert\, \calF_k} \\
             & = \Exp{w_k(Y_k)\nabla g(u_k,Y_k)\,\vert\,\calF_k} - \Exp{w_k(Y_k)\tilde{\Pi}^{\bbV_{m_k}}_{\mathcal{S}_k}  [\nabla g](Y_k)\,\middle\vert\, \calF_k}
            + \Exp{w_k(Y_k)\tilde{\Pi}^{\bbV_{m_k}}_{\mathcal{S}_k} [\nabla g](Y_k) \,\middle\vert\,\calF_k}                                                                                                                                         \\
             & = \Exp{w_k(Y_k)\nabla g(u_k,Y_k) \,\vert\,\calF_k}                                                                                                                                                                                    \\
             & = \nabla J(u_k),
        \end{split}
    \end{equation*}
    we obtain
    \begin{equation} \label{eq:A-and-B}
        \begin{split}
            \Exp{\norm{u_{k+1}-u^*}_U^2\,\middle\vert\, \calF_k} = & \norm{u_k - u^*}_U^2                                                                                                           \\
                                                                   & \underbrace{- 2\tau_k \inp{ u_k -u^*}{\nabla J(u_k)-\nabla J(u^*)}_U }_{:=A}                                                   \\
                                                                   & + \underbrace{\tau_k^2 \Exp{\norm{\widehat{\nabla J}(u_{k},Y_{k};\calS_k)- \nabla J(u^*)}_U^2\,\middle\vert\, \calF_k}}_{:=B}.
        \end{split}
    \end{equation}

    The goal is now to derive upper bounds for the terms $A$ and $B$ in \eqref{eq:A-and-B}.
    To treat $A$, we use the $\alpha$-strong convexity assumption on $J$ to obtain
    \begin{equation*}
        \begin{split}
            A
             & = - 2\tau_k \inp{u_k -u^*}{\nabla J(u_k)-\nabla J(u^*)}_U \\
             & \leq -2 \tau_k \alpha \norm{u_k-u^*}_U^2.
        \end{split}
    \end{equation*}
    To treat $B$, we observe the following
    \begin{equation} \label{eq:D-term-origin}
        \begin{split}
                 & \norm{\widehat{\nabla J}(u_{k},Y_{k};\calS_k)- \nabla J(u^*)}_U^2                              \\
            \leq & 2 \norm{\widehat{\nabla J}(u_{k},Y_{k};\calS_k)- \widehat{\nabla J}(u^*,Y_{k};\calS^*_k)}_U^2
            + 2 \norm{\widehat{\nabla J}(u^*,Y_{k};\calS^*_k) - \nabla J(u^*)}_U^2                                \\
            =    & 2\norm{ \widehat{\nabla J}(u_{k},Y_{k};\calS_k) - \widehat{\nabla J}(u^*,Y_{k};\calS^*_k)}_U^2 \\
                 & + 2 \underbrace{\norm{w_k(Y_k)\left(\nabla g(u^*,Y_k)
            -  \tilde{\Pi}^{\bbV_{m_k}}_{\mathcal{S}^*_k}  [\nabla g](Y_k)\right) + \Exp{w_k(Y_k) \tilde{\Pi}^{\bbV_{m_k}}_{\mathcal{S}^*_k} [ \nabla g](Y_k)\,\middle\vert\, \calF_k}
            - \nabla J(u^*)}_U^2}_{:=C}
        \end{split}
    \end{equation}
    To bound $C$, recall the identity
    \begin{equation}
        \mathbb{E}[\norm{X}_U^2\,\vert\, \calF_k] = \mathbb{E}[\norm{X -\mathbb{E}[X\,\vert\, \calF_k]}_U^2 \,\vert\, \calF_k ] + \norm{\mathbb{E}[X\,\vert\, \calF_k]}_U^2
        \label{eq:ineq-expec}
    \end{equation}
    for any $U$-valued square integrable random variable $X$.
    Hence, taking the conditional expectation with respect to $\calF_k$, yields
    \begin{equation*}
        \begin{split}
            \mathbb{E}[C \mid \calF_k]
            \leq & \Exp{w_k(Y_k)^2 \norm{\nabla g(u^*,Y_k) -  \tilde{\Pi}^{\bbV_{m_k}}_{\mathcal{S}^*_k}  [\nabla g](Y_k)}_U^2\,\middle\vert\, \calF_k},
        \end{split}
    \end{equation*}
    hence the thesis.
\end{proof}
To treat the term
\begin{equation}\label{eq:D-term}
    D
    := \norm{ \widehat{\nabla J}(u_{k},Y_{k};\calS_k)-
        \widehat{\nabla J}(u^*,Y_{k};\calS^*_k)}_U^2
\end{equation}
in \eqref{eq:D-term-origin} we need the following Lemmas.
\begin{lmm}
    {\color{black} Under \Cref{asmp:lip-cont}}, let $\tilde{\Pi}^{\bbV_{m_k}}_{\mathcal{S}_k}, \tilde{\Pi}^{\bbV_{m_k}}_{\mathcal{S}^*_k}$ denote the weighted discrete least-squares approximations \eqref{eq:cond-least-squares-estim} with memory $\calS_k$ and $\calS_k^*$, respectively.
    Denote moreover $G_k \in \R^{m_k \times m_k}$ the Gram matrix
    \begin{equation} \label{eq:gram-mat}
        \left(G_k\right)_{ij} = \frac{1}{s_k} \sum_{l=1}^{s_k} w_k(Y_{k-l}) \inp{\phi_i(Y_{k-l})}{\phi_j(Y_{k-l})}_U,
    \end{equation}
    and let $q_k := \sup_{y \in \Gamma} w_k(y)$, then
    \begin{equation*}
        \begin{split}
                 & \Exp{w_k(Y_k)^2 \norm{\tilde{\Pi}^{\bbV_{m_k}}_{\mathcal{S}_k}  [\nabla g](Y_k)-\tilde{\Pi}^{\bbV_{m_k}}_{\mathcal{S}^*_k}  [\nabla g](Y_k)}^2_U \,\middle\vert\, \calF_k} \\
            \leq & \norm{G_k^{-1}}_2^2 \norm{G_k}_2 \frac{q_k}{s_k} \sum_{i=1}^{s_k} w_{k}(Y_{k-i})\ell(Y_{k-i})^2 \norm{u_{k-i} - u^*}_U^2\mathds{1}_{\{\norm{G_k-I}_2\leq \frac{1}{2}\}},
        \end{split}
    \end{equation*}
    where $s_k$ denotes the sample (memory) size and $\ell(y)$ is the Lipschitz constant in \Cref{asmp:lip-cont}.
    \label{lem:lemma-intermediate}
\end{lmm}

\begin{proof}
    We recall that the gradient approximation is computed via weighted discrete least-squares using an orthonormal basis of $\bbV_{m_k}$.
    This means that we can express the projections as follows
    \begin{equation*}
        \begin{split}
            \tilde{\Pi}^{\bbV_{m_k}}_{\mathcal{S}_k}  [\nabla g] = \sum_{i=1}^{m_k} c_{k,i} \phi_i,
            \quad
            \tilde{\Pi}^{\bbV_{m_k}}_{\mathcal{S}^*_k}  [\nabla g] = \sum_{i=1}^{m_k} c_i^* \phi_i,
        \end{split}
    \end{equation*}
    whenever the least-squares problem is well conditioned, that is $\{\norm{G_k-I}_2\leq \frac{1}{2}\}$.
    Moreover, the vectors $c_k= (c_{k,1}, \dots , c_{k,m_k})$ and $c^*_k= (c^*_1, \dots , c^*_{m_k})$ are obtained by solving the following systems
    \begin{equation*}
        \begin{split}
            G_k c_k = \frac{1}{s_k} J_k,
            \qquad G_k c^* = \frac{1}{s_k} J^*_k,
        \end{split}
    \end{equation*}
    where we defined
    \begin{equation*}
        J_k := \inp{\sqrt{W_k} V_k}{\sqrt{W_k}\hat{\Psi}_k}_U,
        \quad
        J^*_k := \inp{\sqrt{W_k} V_k}{\sqrt{W_k}\Psi^*_k}_U,
    \end{equation*}
    with $V_k$ the Vandermonde matrix, $W_k$ the weight matrix, and
    \begin{equation*}
        \hat{\Psi}_k = \calS_k(2) =
        \begin{pmatrix}
            \nabla g(u_{k-1},Y_{k-1}) \\
            \vdots                    \\
            \nabla g(u_{k-s_k},Y_{k-s_k})
        \end{pmatrix}
        , \quad \Psi^*_k = \calS_k^\ast(2) =
        \begin{pmatrix}
            \nabla g(u^*,Y_{k-1}) \\
            \vdots                \\
            \nabla g(u^*,Y_{k-s_k})
        \end{pmatrix}
        .
    \end{equation*}
    Hence, we can write
    \begin{equation*}
        \begin{cases}
            c_k = \frac{1}{s_k} G_k^{-1} J_k , \quad c^* = \frac{1}{s_k} G_k^{-1} J^*_k & \text{if } \norm{G_k-I}_2 \leq \frac{1}{2}, \\
            c_k=c^*=0                                                                   & \text{otherwise},
        \end{cases}
    \end{equation*}
    so that
    \begin{equation*}
        \norm{c_k -c^*_k}_2^2 \leq
        \begin{cases}
            \norm{G_k^{-1}}_2^2 \norm{\frac{1}{s_k} (\sqrt{W_k}V_k)^T}_{{\calL(U^{s_k},\R^{m_k})}}^2 \norm{\sqrt{W_k}(\Psi_k-\Psi^*_k)}_{U^{s_k}}^2 \quad & \norm{G_k-I}_2 \leq \frac{1}{2}, \\
            0 \quad                                                                                                                                       & \text{otherwise}.
        \end{cases}
    \end{equation*}
    The fact that the Gram matrix $G_k$ is symmetric positive definite whenever $\norm{G_k-I}_2 \leq \frac{1}{2}$ allows us to write its spectral norm as
    \begin{equation*}
        \norm{G_k}_2 = \sup_{\substack{v \neq 0 \\ v\in \R^{m_k}}} \frac{v^TG_kv}{\norm{v}_2^2}.
    \end{equation*}
    Furthermore,
    \begin{equation*}
        \begin{split}
            \norm{\frac{1}{s_k}(\sqrt{W_k}V_k)^T}_{{\calL(U^{s_k},\R^{m_k})}}^2 = \norm{\frac{1}{s_k} \sqrt{W_k}V_k}_{{{\calL(\R^{m_k},U^{s_k})}}}^2
             & = \sup_{\substack{v \neq 0               \\ v\in \R^{m_k}}} \frac{\norm{\frac{1}{s_k}\sqrt{W_k}V_k v}_{{{\calL(\R^{m_k},U^{s_k})}}}^2}{\norm{v}_2^2} \\
             & = \frac{1}{s_k} \sup_{\substack{v \neq 0 \\ v\in \R^{m_k}}} \frac{v^TG_kv}{\norm{v}_2^2}          \\
             & = \frac{\norm{G_k}_2}{s_k}.
        \end{split}
    \end{equation*}
    Then, taking the definition of the spectral norm and recalling that $g(\cdot,y)$ has $\ell(y)$-Lipschitz gradient yields
    \begin{equation*}
        \begin{split}
            \norm{c_k -c^*_k}_2^2 \mathds{1}_{\{\norm{G_k-I}_2 \leq \frac{1}{2}\}}
             & \leq \norm{G_k^{-1}}^2 \frac{\norm{G_k}}{s_k} \sum_{i=1}^{s_k} \norm{\sqrt{W_k}((\hat{\Psi}_{k})_i-(\Psi_{k}^*)_i)}_{U}^2 \mathds{1}_{\{\norm{G_k-I}_2 \leq \frac{1}{2}\}}          \\
             & = \norm{G_k^{-1}}^2 \frac{\norm{G_k}}{s_k} \sum_{i=1}^{s_k} w_k(Y_{k-i})\norm{\nabla g(u_{k-i},Y_{k-i})-\nabla g(u^*,Y_{k-i})}_U^2 \mathds{1}_{\{\norm{G_k-I}_2 \leq \frac{1}{2}\}} \\
             & \leq \norm{G_k^{-1}}^2 \frac{\norm{G_k}}{s_k} \sum_{i=1}^{s_k} w_k(Y_{k-i}) \ell(Y_{k-i})^2 \norm{u_{k-i} - u^*}_U^2 \mathds{1}_{\{\norm{G_k-I}_2 \leq \frac{1}{2}\}}.
        \end{split}
    \end{equation*}
    Finally, the thesis follows from the fact that for any $v = \sum_{i=1}^{m_k} c_i \phi_i \in \bbV_{m_k}$ with $\{c_i\}_{i=1}^{m_k}$ $\calF_k$-measurable
    \begin{equation*}
        \Exp{w_k(Y_k)^2 \norm{v(Y_k)}^2_U \,\middle\vert\, \calF_k}
        \leq q_k \bbE_{Y\sim\rho}\left[ \norm{v(Y)}^2_U \,\middle\vert\, \calF_k \right]
        = q_k \norm{c}_2^2.
    \end{equation*}
\end{proof}

We can now give an upper bound of the quantity $D$ in \eqref{eq:D-term}.

\begin{lmm}
    {\color{black} Under \Cref{asmp:str-cvx,asmp:lip-cont}}, let $\widehat{\nabla J}(u_{k},Y_{k};\calS_k)$ be as in \eqref{eq:hat-grad-J}, $D$ as in \eqref{eq:D-term}, $u^*=\mathrm{argmin}_u J(u)$, and $u_k$ the $k$-th iterate of \Cref{alg:algo-fixed} or \ref{algo:algo-variable}.
    Let further $G_k$ be the Gram matrix at iteration $k$ as in \eqref{eq:gram-mat}, $q_k = \sup_{y\in\Gamma} w_k(y)$, and $s_k$ be the memory size to perform the gradient approximation.
    Then,
    \begin{equation*}
        \begin{split}
            \Exp{D \,\vert\, \calF_k}
            =    & \Exp{\norm{\widehat{\nabla J}(u_{k},Y_{k};\calS_k)- \widehat{\nabla J}(u^*,Y_{k};\calS_k^*)}_U^2\,\middle\vert\, \calF_k}                                                                               \\
            \leq & 2q_k L^2 \norm{u_k-u^*}_U^2 + 2 \norm{G_k^{-1}}_2^2  \norm{G_k}_2 \frac{q_k}{s_k} \sum_{i=1}^{s_k} w_{k}(Y_{k-i})\ell(Y_{k-i})^2 \norm{u_{k-i} - u^*}_U^2\mathds{1}_{\{\norm{G_k-I}\leq \frac{1}{2}\}}.
        \end{split}
    \end{equation*}
    \label{lemma_hard}
\end{lmm}

\begin{proof}
    Let us start by expanding $\widehat{\nabla J}(u_{k},Y_{k};\calS_k)$ and using the triangular inequality:
    \begin{equation*}
        \begin{split}
            D =  & \norm{\widehat{\nabla J}(u_{k},Y_{k};\calS_k)- \widehat{\nabla J}(u^*,Y_{k};\calS_k^*)}_U^2                                                                                                                                                                                                                                            \\
            =    & \left\|w_k(Y_k)\left(\nabla g (u_k,Y_k) - \tilde{\Pi}^{\bbV_{m_k}}_{\mathcal{S}_k}  [\nabla g](Y_k)\right) + \Exp{w_k(Y_k) \tilde{\Pi}^{\bbV_{m_k}}_{\mathcal{S}_k}  [\nabla g](Y_k)\,\middle\vert\, \calF_k} \right.                                                                                                                    \\
                 & \left.-w_k(Y_k)\left(\nabla g(u^*,Y_k) -  \tilde{\Pi}^{\bbV_{m_k}}_{\mathcal{S}^*_k}  [\nabla g](Y_k)\right) - \Exp{w_k(Y_k) \tilde{\Pi}^{\bbV_{m_k}}_{\mathcal{S}^*_k}  [\nabla g](Y_k)\,\middle\vert\, \calF_k}\right\|^2_U                                                                                                          \\
            \leq & 2 w_k(Y_k)^2\norm{\nabla g (u_k,Y_k)- \nabla g(u^*,Y_k)}_U^2                                                                                                                                                                                                                                                                           \\
                 & + 2 \norm{\Exp{w_k(Y_k)\left(\tilde{\Pi}^{\bbV_{m_k}}_{\mathcal{S}_k}  [\nabla g](Y_k)-  \tilde{\Pi}^{\bbV_{m_k}}_{\mathcal{S}^*_k}  [\nabla g](Y_k)\right) \,\middle\vert\,\calF_k} - w_k(Y_k)\left(\tilde{\Pi}^{\bbV_{m_k}}_{\mathcal{S}_k}  [\nabla g](Y_k)-\tilde{\Pi}^{\bbV_{m_k}}_{\mathcal{S}^*_k}  [\nabla g](Y_k)\right)}^2_U
        \end{split}
    \end{equation*}
    We now take the conditional expectation with respect to $\calF_k$ of the above inequality and use \eqref{eq:ineq-expec}, the fact that $g(\cdot,y)$ has Lipschitz gradient with constant $\ell(y)$, and $w_k(y) \leq q_k$ to obtain
    \begin{equation*}
        \begin{split}
            \mathbb{E}[D \mid \calF_k]\leq & 2 \Exp{w_k(Y_k)^2\ell(Y_k)^2} \norm{u_k-u^*}_U^2 \notag                                                                                                                                           \\
                                           & + 2 \Exp{w_k(Y_k)^2 \norm{\tilde{\Pi}^{\bbV_{m_k}}_{\mathcal{S}_k}  [\nabla g](Y_k)-\tilde{\Pi}^{\bbV_{m_k}}_{\mathcal{S}^*_k}  [\nabla g](Y_k)}_U^2\,\middle\vert\, \calF_k}, \label{eq:eq-hard}
        \end{split}
    \end{equation*}
    and, from \Cref{lem:lemma-intermediate},
    \begin{equation*}
        \begin{split}
                 & \Exp{D\,\vert\,\calF_k}                                                                                                                                                                               \\
            \leq & 2q_k L^2 \norm{u_k-u^*}_U^2 + 2 \norm{G_k^{-1}}_2^2  \norm{G_k}_2 \frac{q_k}{s_k} \sum_{i=1}^{s_k} w_k(Y_{k-i})\ell(Y_{k-i})^2 \norm{u_{k-i} - u^*}_U^2\mathds{1}_{\{\norm{G_k-I}\leq \frac{1}{2}\}}.
        \end{split}
    \end{equation*}
\end{proof}

We can now prove the following error recursion.
\begin{thrm}
    {\color{black} Under \Cref{asmp:str-cvx,asmp:lip-cont}}, let $u^*=\mathrm{argmin}_u J(u)$, and $u_k$ the $k$-th iterate of \Cref{alg:algo-fixed} or \ref{algo:algo-variable}.
    Let further $\tau_k$ be the step size, $s_k$ be the memory size to perform the gradient approximation, and $q_k = \sup_{y\in\Gamma} w_k(y)$.
    Then,
    \begin{equation}
        \begin{split}
            \Exp{\norm{u_{k+1}-u^*}_U^2} \leq & \left( 1 - 2 \tau_k \alpha + 4 q_k \tau_k^2 L^2\right) \mathbb{E}[\norm{u_k-u^*}_U^2]                                                \\
            +                                 & 2 q_k \tau_k^2 \bbE_{Y\sim\rho}\left[\norm{ \nabla g(u^*,Y) -  \tilde{\Pi}^{\bbV_{m_k}}_{\mathcal{S}^*_k}  [\nabla g](Y)}_U^2\right] \\
            +                                 & 24  \tau_k^2 \frac{q_k}{s_k} \sum_{i=1}^{s_k} L^2 \Exp{\norm{u_{k-i} - u^*}_U^2}.
        \end{split}
        \label{eq:theorem-rec}
    \end{equation}
    \label{thm:theorem-rec}
\end{thrm}

\begin{proof}
    Putting together the results from \Cref{lem:first} and \Cref{lemma_hard} we have
    \begin{equation*}
        \begin{split}
            \Exp{\norm{u_{k+1}-u^*}_U^2\,\vert\, \calF_k} \leq & (1-2 \tau_k \alpha)\norm{u_k-u^*}_U^2                                                                                                                                                                                \\
                                                               & + 2 q_k \tau_k^2 \Exp{w_k(Y_k) \norm{ \nabla g(u^*,Y_k) -  \tilde{\Pi}^{\bbV_{m_k}}_{\mathcal{S}^*_k}  [\nabla g](Y_k)}_U^2\,\middle\vert\, \calF_k}                                                                 \\
                                                               & + 4q_k \tau_k^2 L^2 \norm{u_k-u^*}_U^2                                                                                                                                                                               \\
                                                               & + 4 \tau_k^2 \frac{q_k}{s_k} \Exp{\norm{G_k^{-1}}_2^2  \norm{G_k}_2 \sum_{i=1}^{s_k} w_{k}(Y_{k-i})\ell(Y_{k-i})^2 \norm{u_{k-i} - u^*}_U^2\mathds{1}_{\{\norm{G_k-I}_2\leq \frac{1}{2}\}}\,\middle\vert\, \calF_k}.
        \end{split}
    \end{equation*}
    We now take the expectation of the whole expression and, using the law of total expectation and the fact that $Y_{k-i}$ is independent of $u_{k-i}$,
    \begin{equation*}
        \begin{split}
            \Exp{\norm{u_{k+1}-u^*}_U^2} \leq & \left( 1 - 2 \tau_k \alpha + 4 q_k \tau_k^2 L^2\right) \Exp{\norm{u_k-u^*}_U^2}                                         \\
            +                                 & 2 \tau_k^2 q_k \Exp{w_k(Y_k)\norm{ \nabla g(u^*,Y_k) - \tilde{\Pi}^{\bbV_{m_k}}_{\mathcal{S}^*_k} [\nabla g](Y_k)}_U^2} \\
            +                                 & 24 \tau_k^2 \frac{q_k}{s_k} \sum_{i=1}^{s_k} L^2 \mathbb{E}\left[ \norm{u_{k-i} - u^*}_U^2 \right],
        \end{split}
    \end{equation*}
    where we used the independence of $Y_{k-i}$ and $u_{k-i}$, and the fact that $\norm{G_k-I}_2\leq \frac{1}{2}$ implies $\norm{G_k}_2\leq \frac{3}{2}$ and $\norm{G_k^{-1}}_2\leq 2$.
\end{proof}

\begin{crllr} \label{thm:rec-simplified}
    Assume the same hypothesis of \Cref{thm:theorem-rec} are satisfied and that the memory size $s_k$ satisfies the sampling inequality \eqref{eq:samp-ineq}.
    Introduce the notation
    \begin{equation} \label{eq:error-notation}
        e_k := \Exp{\norm{u_{k}-u^*}_U^2}
    \end{equation}
    and define
    \begin{equation*}
        D_{m_k}
        := (1 + 2 \epsilon(s_k)) e_{m_k}(\nabla g(u^*,\cdot))^2 + 2 \norm{\nabla g(u^\ast,\cdot)}_{L^2_\rho(\Gamma;U)}^2 s_{k}^{-r},
    \end{equation*}
    where $e_{m_k}(\nabla g(u^*,\cdot))$ denotes the best approximation error of $\nabla g(u^*,\cdot)$ in $\bbV_{m_k}$ (see \eqref{eq:best-approx}) and $\epsilon(s) = \frac{4 \kappa}{\log(s)}$ with $\kappa =\frac{1- \log(2)}{2+2r}$.
    Then, for all $k \geq 0$,
    \begin{equation}
        \begin{split}
            e_{k+1} \leq & a_k e_k
            + 2 q_k \tau_k^2 D_{m_k}
            +b_k \sum_{i=1}^{s_k} e_{k-i},
        \end{split}
        \label{eq:theorem-rec-simplified}
    \end{equation}
    where $a_k =  1 - 2 \tau_k \alpha + 4 q_k \tau_k^2 L^2$ and $b_k= 24 \tau_k^2 L^2 \frac{q_k}{s_k}$.
\end{crllr}
\begin{proof}
    Since the sampling inequality \eqref{eq:samp-ineq} is satisfied, per \Cref{thm:theorem-approx-lstsq} ({\color{black} which we apply to the least-squares approximation with idealized memory, hence noiseless data}) the second term in the right-hand side of \eqref{eq:theorem-rec} can be bounded as
    \begin{equation*}
        \bbE_{Y\sim\rho}\left[\norm{ \nabla g(u^*,Y) - \tilde{\Pi}^{\bbV_{m_k}}_{\mathcal{S}^*_k} [\nabla g](Y)}_U^2\right]
        \leq (1 + 2 \epsilon(s_k)) e_{m_k}(\nabla g(u^*,\cdot))^2 + 2 \norm{\nabla g(u^\ast,\cdot)}_{L^2_\rho(\Gamma;U)}^2 s_{k}^{-r} = D_{m_k}.
    \end{equation*}
\end{proof}

We now prove a useful lemma that is needed for the convergence analysis in the next two sections.

\begin{lmm} \label{thm:rec_X}
    Under the same hypotheses and definitions of \Cref{thm:rec-simplified},
    let $\bar{k} \geq 0$ and assume further that there exists a sequence $\{X_k\}_{k\geq \bar{k}-s_{\bar{k}}} \subset \R_+$ such that
    \begin{equation} \label{eq:X_k_rec_rel}
        X_{k+1} \geq a_k X_k + b_k \sum_{i=1}^{s_k} X_{k-i}, \quad k\geq \bar{k},
    \end{equation}
    and, if $\bar{k}-s_{\bar{k}} < 0$, $X_k \geq 1$ for all $k = \bar{k}-s_{\bar{k}},\dots,0$.
    Then, the error satisfies the bound
    \begin{equation*}
        e_k \leq X_k \left( \tilde{e}_{\bar{k}} + 2 \sum_{j=\bar{k}}^{k-1} \frac{q_j \tau_j^2 D_{m_j}}{X_{j+1}} \right), \quad k\geq \bar{k}+1,
    \end{equation*}
    where
    \begin{equation}\label{eq:init-err-notation}
        \tilde{e}_{\bar{k}} :=
        \begin{cases}
            \max\{e_{\bar{k}-s_{\bar{k}}},\dots,e_{0}\} & \text{if $\bar{k}-s_{\bar{k}} < 0$}, \\
            e_{\bar{k}}                                 & \text{otherwise}.
        \end{cases}
    \end{equation}
\end{lmm}
\begin{proof}
    We prove the case $\bar{k} = 0$, then the result for $\bar{k} \geq 1$ follows by a translation argument.
    Let us make the following ansatz.
    For any $k \geq 0$, there exist $Y_{j,k} > 0$, $j \in \{0,\cdots,k-1 \}$, such that
    \begin{equation}
        e_k \leq X_k \tilde{e}_0 + \sum_{j=0}^{k-1}Y_{j,k} q_j \tau_j^2 D_{m_j}.
        \label{eq:err_bound_with_X_and_Y}
    \end{equation}
    We then plug \eqref{eq:err_bound_with_X_and_Y} into \eqref{eq:theorem-rec-simplified} to obtain
    \begin{equation*}
        \begin{split}
            e_{k+1} \leq & a_k \left[ X_k\tilde{e}_0 + \sum_{j=0}^{k-1} Y_{j,k} q_j \tau_j^2 D_{m_j} \right]      + 2 q_k \tau_k^2 D_{m_k}                                                                                                \\
                         & + b_k \sum_{i=1}^{\min\{s_k,k-1\}}\left[ X_{k-i}\tilde{e}_0 + \sum_{j=0}^{k-i-1} Y_{j,k-i} q_j \tau_j^2 D_{m_j} \right]        + b_k \sum_{i=1}^{\max\{0,s_k-k+1\}} e_{1-i}                                    \\
            \leq         & \tilde{e}_0 \left[ a_k X_k + b_k \sum_{i=1}^{\min\{s_k,k-1\}} X_{k-i} \right]         + 2 q_k \tau_k^2 D_{m_k}                                                                                                 \\
                         & + \sum_{j=0}^{k-1}  q_j \tau_j^2 D_{m_j} \left[ a_k Y_{j,k} + b_k \sum_{i=1}^{\min\{s_k,k-j-1\}} Y_{j,k-i} \right]       + b_k \sum_{i=1}^{\max\{0,s_k-k+1\}} X_{1-i} \tilde{e}_0,                             \\
            =            & \tilde{e}_0 \left[ a_k X_k + b_k \sum_{i=1}^{s_k} X_{k-i} \right] + 2 q_k \tau_k^2 D_{m_k} + \sum_{j=0}^{k-1}  q_j \tau_j^2 D_{m_j} \left[ a_k Y_{j,k} + b_k \sum_{i=1}^{\min\{s_k,k-j-1\}} Y_{j,k-i} \right],
        \end{split}
    \end{equation*}
    where we used $X_i \geq 1$ whenever $i \leq 0$.
    Notice that the above can be further bounded by
    \begin{equation*}
        e_{k+1} \leq X_{k+1}\tilde{e}_0 + \sum_{j=0}^{k} Y_{j,k+1} q_j \tau_j^2 D_{m_j}
    \end{equation*}
    owing to \eqref{eq:X_k_rec_rel} and if we can exhibit a sequence $\{Y_{j,k}\}$ such that for any $k \geq 0$
    \begin{align}
         & Y_{j,k+1} \geq
        \begin{cases}
            a_k Y_{j,k} + b_k \sum_{i=1}^{\min\{s_k,k-j-1\}} Y_{j,k-i}, & \quad \text{if } j \in \{0,\cdots,k-1\} \\
            2,                                                          & \quad \text{if } j=k
            \label{eq:Y_k_fixed}
        \end{cases}.
    \end{align}
    Such a sequence can be constructed using the sequence $\{X_k\}_k$.
    Indeed, setting $Y_{j,k+1} = 2\frac{X_{k+1}}{X_{j+1}}$ gives, for $j=k$, $Y_{k,k+1} = 2 \geq 2$ and, for $j<k$,
    \begin{equation*}
        \begin{split}
            Y_{j,k+1}
             & = 2\frac{X_{k+1}}{X_{j+1}}                                                                        \\
             & \geq a_{k} 2 \frac{X_{k}}{X_{j+1}}  + b_k \sum_{i=1}^{s_k} 2 \frac{X_{k-i}}{X_{j+1}}              \\
             & \geq a_{k} 2 \frac{X_{k}}{X_{j+1}}  + b_k \sum_{i=1}^{\min\{s_k,k-j-1\}} 2\frac{X_{k-i}}{X_{j+1}} \\
             & = a_{k} Y_{j,k} + b_k \sum_{i=1}^{\min\{s_k,k-j-1\}} Y_{j,k-i},
        \end{split}
    \end{equation*}
    where we used \eqref{eq:X_k_rec_rel} in the second inequality, so that \eqref{eq:Y_k_fixed} is satisfied.
    This proves the claim.
\end{proof}

\subsection{Convergence of the fixed approximation space algorithm}
Let us now prove the convergence of the fixed approximation space \Cref{alg:algo-fixed} with a constant step size.
In this version of the algorithm, one may be limited by the accuracy of the approximation in $\bbV_m$.
    {\color{black} We demonstrate that this algorithm exhibits an exponential convergence rate up to an additive constant term which depends on the accuracy of the best approximation in $\bbV_m$ of $\nabla g(u^*,\cdot)$.
        In this section, to ease notation, we shall drop the subscript $k$ for all the quantities which are kept constant throughout the iterations.}

\begin{thrm} \label{thm:conv-fixed}
    {\color{black} Fix the approximation space $\bbV_m$.
        Under \Cref{asmp:str-cvx,asmp:lip-cont}, let $e_k$ be the squared error defined in \eqref{eq:error-notation}, $\tau > 0$ be the constant step size, $s$ be the constant memory size satisfying the sampling inequality \eqref{eq:samp-ineq}, $\tilde{e}_0$ be the initial error defined as in \eqref{eq:init-err-notation}, and $q = \sup_{y\in\Gamma} w(y)$.}
    If ${\color{black} 0 < }\tau < \frac{\alpha}{14 q L^2}$, then there exists ${\color{black} x = x(\tau)} \in (0,1)$ such that
    \begin{equation}\label{eq:thm-conv-fixed}
        e_k \leq x^k \tilde{e}_0 + \frac{2 q \tau^2 D_m}{1-x}.
    \end{equation}
    In particular, if
    \begin{equation*}
        {\color{black} 0 < }
        \tau
        \leq \frac{\alpha}{4qL^2}\min\left\{ \frac{1 - \left( \frac{1}{1+s \delta} \right)^{\frac{1}{s}}}{1+ 6\delta - \left( \frac{1}{1+s \delta} \right)^{\frac{1}{s}}}, \frac{\alpha^2 - 24 q L^2 \delta}{\alpha^2} \right\},
        \qquad \delta < \frac{\alpha^2}{24qL^2},
    \end{equation*}
    then \eqref{eq:thm-conv-fixed} is true for
    \begin{equation}
        \label{eq:x-of-tau}
        {\color{black} x = \frac{\alpha - 4q \tau L^2(1+6\delta)}{\alpha - 4q \tau L^2}.}
    \end{equation}
\end{thrm}

\begin{proof}
    As the step size and memory size are fixed, we shall use the simplified notation $a$ and $b$ for the terms $a_k$ and $b_k$ in \Cref{thm:rec-simplified}.
    Let us first show that \eqref{eq:X_k_rec_rel} is true for $X_k = x^k$ for some $x \in (0,1)$, $k \geq -s$.
    Substituting $x^k$ for $X_k$ yields
    \begin{equation*}
        \begin{split}
            x^{k+1} \geq a x^k + b \sum_{i=1}^{s} x^{k-i}
            = ax^k + b \frac{x^{k-s}-x^{k}}{1-x}
        \end{split}
    \end{equation*}
    which, dividing by $x^k>0$, reads
    \begin{equation*}
        x \geq a + b \frac{x^{-s}-1}{1-x}.
    \end{equation*}
    The above is in turn equivalent to
    \begin{equation}
        x^{s}(-x^2 +(1+a)x+b-a)-b \geq 0.
        \label{ineq:degree_s_polynomial}
    \end{equation}
    We now aim to show that, for a suitable choice of the step size, the inequality above is satisfied for some $x \in (0,1)$.
    We propose two alternative proofs of this result, with different choices of step size.
    The second one, in particular, will be useful to analyze the variable subspace algorithm discussed in the next subsection.

    The first approach yields a step size working for any $m$.
    Let us define the function corresponding to the left-hand-side of \eqref{ineq:degree_s_polynomial}
    \begin{equation*}
        f(x) = x^{s}(-x^2 +(1+a)x+b-a)-b,
    \end{equation*}
    and notice that $f$ is $C^1(\R)$ since it is a polynomial.
    Moreover, observe that in $x=0$ we have $f(0)=-b<0$ and
    in $x=1$ we have $f(1)=0$.
    Therefore, if $f'(1)<0$, there exists $x \in (0,1)$ such that $f(x)>0$.
    In particular, this entails that $x^{*}:=\inf\{x\in (0,1) \,:\,f(x)>0\}$ is well-defined.
    The derivative of $f$ is
    \begin{equation*}
        f'(x) = sx^{s-1}(-x^2 + (1+a)x+b-a) + x^{s}(-2x + 1 + a),
    \end{equation*}
    which, evaluated at $x=1$, yields
    \begin{equation*}
        f'(1)= a + sb -1 = 28 q \tau^2 L^2 - 2 \tau \alpha.
    \end{equation*}
    Hence, setting $f'(1)<0$ and solving for $\tau>0$ gives
    \begin{equation*}
        \tau <  \frac{\alpha}{14q L^2}.
    \end{equation*}

    The second approach yields a step size which depends on $s$, namely the memory size, which, in turn, depends on the chosen $m$ via the sampling inequality \eqref{eq:samp-ineq}.
    Let us first rewrite \eqref{ineq:degree_s_polynomial} as
    \begin{equation*}
        x^{s}(-x^2+(1+a)x+b-a) \geq b.
    \end{equation*}
    The above states that the product of two polynomials must be at least as great as $b$.
    Since it is in general hard to solve the equation for a generic $s$, we assume that the inequality (up to some constant $1+s \delta$, $\delta > 0$) is satisfied by the degree-two polynomial, and then we solve for the degree-$s$ polynomial.
    More specifically, we impose first
    \begin{equation}
        -x^2+(1+a)x+b-a \geq (1+s \delta)b, \quad \delta>0,
        \label{ineq:second_degree_second_approach}
    \end{equation}
    and, expanding the terms $a$ and $b$ in the left-hand side of \eqref{ineq:second_degree_second_approach}, yields the following second degree polynomial in $\tau$
    \begin{equation}\label{eq:second-deg-approach-in-tau}
        4q \tau^2L^2 \left( x -1 - 6\delta \right) + 2\tau\alpha \left( 1-x \right) - \left( 1-x \right)^2 \geq 0.
    \end{equation}
    We compute the discriminant
    \begin{equation*}
        \begin{split}
            \Delta & = 4\alpha^2\left( 1-x \right)^2 + 16qL^2\left( x -1 - 6\delta  \right) \left( 1-x \right)^2 \\
                   & = \left( 1-x \right)^2 \left(4\alpha^2 + 16qL^2 \left(x -1 - 6\delta\right) \right)
        \end{split}
    \end{equation*}
    For the equation to admit solutions, we need $\Delta \geq 0$.
    As $\left( 1-x \right)^2  \geq 0$ for any $x\in (0,1)$, the condition becomes
    \begin{equation*}
        \alpha^2 + 4qL^2\left(x -1 - 6 \delta \right) \geq 0.
    \end{equation*}
    Solving for $x$ the above yields
    \begin{equation*}
        \begin{split}
            x \geq 1 - \frac{\alpha^2 - 24qL^2\delta}{4 q L^2}.
        \end{split}
    \end{equation*}
    Since we are looking for $x\in (0,1)$, we need to ensure that the right-hand side of the above inequality is strictly less than $1$, that is
    \begin{align*}
        \alpha^2 - 24qL^2\delta > 0 \quad
        \Leftrightarrow \quad \delta < \frac{\alpha^2}{24qL^2}.
    \end{align*}
    Then note that, as $x-1-6\delta < 0$, the sign of the coefficient in front of the second degree term of \eqref{eq:second-deg-approach-in-tau} is negative.
    Hence, the set of solutions of \eqref{eq:second-deg-approach-in-tau} coincides with the points that lie between the two roots of the polynomial.
    In particular, we can choose
    \begin{equation}\label{eq:tau-x-second-approach}
        \tau = \frac{\alpha\left( 1-x \right)}{4q L^2\left( 1+ 6\delta - x\right)},
    \end{equation}
    {\color{black} which leads to \eqref{eq:x-of-tau}.}
    With this choice of $\tau$, we may write
    \begin{equation*}
        x^{s}(-x^2+(1+a)x+b-a) \geq x^{s}(1+s \delta)b,
    \end{equation*}
    hence, to guarantee that the right-hand side is greater or equal to $b$, we need to impose the further condition
    \begin{equation*}
        x \geq \left( \frac{1}{1+s \delta} \right)^{\frac{1}{s}}.
    \end{equation*}
    Ultimately, the condition
    \begin{equation}\label{eq:x-second-approach-inproof}
        1 > x \geq \max \left\{ \left( \frac{1}{1+s \delta} \right)^{\frac{1}{s}}, 1 - \frac{\alpha^2 - 24qL^2\delta }{4 qL^2} \right\},
        \qquad \delta < \frac{\alpha^2}{24qL^2},
    \end{equation}
    leads to \eqref{ineq:degree_s_polynomial}.
    Then, per \eqref{eq:tau-x-second-approach} and the fact that $x \mapsto \frac{1-x}{1 + 6\delta - x}$ is monotonically decreasing, the above implies that the step size should satisfy the condition
    \begin{equation*}
        {\color{black} 0 <}
        \tau
        \leq \frac{\alpha}{4qL^2}\min\left\{ \frac{1 - \left( \frac{1}{1+s \delta} \right)^{\frac{1}{s}}}{1+ 6 \delta - \left( \frac{1}{1+s \delta} \right)^{\frac{1}{s}}}, \frac{\alpha^2 - 24q L^2 \delta}{\alpha^2} \right\},
        \qquad \delta < \frac{\alpha^2}{24qL^2}.
    \end{equation*}
    The conclusion follows from \Cref{thm:rec_X}.
\end{proof}

\begin{rmrk}\label{rem:rate-large-s}
    Note that the condition $\tau < \frac{\alpha}{14 q L^2}$ in the \Cref{thm:conv-fixed} does not depend on the size of the memory $s$ and, hence, on the size $m$ of the approximation space (provided we assume no implicit dependence of $q$ on $s$).
    We interpret this as a minimal condition to have the initial error which decays to zero, however without an explicit rate $x$.
    The second condition is more involved, yet it quantifies the speed at which the initial error decay happens.
    In particular, it is interesting to consider the case $s \gg 1$, where the first term in the maximum in \eqref{eq:x-second-approach-inproof} dominates.
    Expanding up to the first order yields
    \begin{equation*}
        x \sim 1 - \frac{\log(s \delta)}{s},
    \end{equation*}
    which reveals how the convergence rate deteriorates as the memory size increases.
        {\color{black} Furthermore, owing to the sampling inequality \eqref{eq:samp-ineq}, if the sampling measure $\mu$ is quasi-optimal then $m = \calO(s/\log(s))$, hence the convergence rate deteriorates linearly in $m^{-1}$ up to logarithmic terms.
            Since $m$ is the dimension of the approximation space $\bbV_m$, namely the number of degrees of freedom in our model, this result mirrors the deterioration of the convergence rate of SAGA \cite{defazio2014saga}, whose analysis leads to a rate of the form $x = 1 - \frac{C}{s}$ for some constant $C>0$, and where the number of degrees of freedom coincides with the dataset size $s$.
            On the other hand, in the case of a finite-sum optimization problem, our algorithm could select a much smaller approximation space than the number of terms in the sum (e.g., because of smoothness properties which allow for extrapolation between points in the dataset) and, hence, one could get away with a better convergence rate by needing to visit fewer points.
        }
\end{rmrk}

\subsection{Convergence of the variable approximation space algorithm}
We shall now turn to the proof of a similar error bound for the case of a variable approximation space.
In order to give an intuition regarding the bound, let us first note that in \Cref{thm:conv-fixed} we assumed $X_k = x^k$ for some $x \in (0,1)$, $X_k$ being the cumulative error contraction factor at iteration $k$ from \Cref{thm:rec_X}, namely we assumed that at each iteration the error contracts by a factor $x$ -- independently of the iteration.
This is justified by the fact that the approximation space and the step size in \Cref{alg:algo-fixed} are fixed over the iterations.
Furthermore, owing to \eqref{eq:x-of-tau}, it holds
$x = 1 - \gamma \tau + \mathcal{O}(\tau^2)$ as $\tau \to 0$ for some $\gamma >0$.
Combining this with \Cref{rem:rate-large-s} reveals that, up to logarithmic terms, at iteration $k$ we should pick the step size $\tau$ inversely proportional to the memory size $s$.
This result then suggests that, in the setting of \Cref{algo:algo-variable} where the approximation space and, owing to the sampling inequality \eqref{eq:samp-ineq}, the memory size $s_k$ grow larger over the iterations, we shall pick an iteration-dependent step size $\tau_k \propto s_k^{-1}$ to obtain a contraction of $x_k = 1 - \gamma \tau_k$ for some $\gamma >0$ at iteration $k$.
This, in turn, entails that a more reasonable model for the cumulative contraction factor is $X_k = \prod_{j=0}^{k} x_j$.
Indeed, we can prove the following {\color{black} lemma}.
\begin{lmm} \label{thm:conv_var1}
    {\color{black} Under the same assumptions and definitions of \Cref{thm:rec-simplified}}, let, for all $k \geq 0$,
    \begin{equation*}
        \tau_k = \tau (k+1)^{t-1},
        \quad s_k = \ceil{s (k+1)^{1-t}},
    \end{equation*}
    for some constants $t \in [0,1)$ and $\tau,s>0$ with $s<1$ and $2 \alpha \tau \neq 1$ if $t=0$.
    Assume further that $q_k = o(\tau_k^{-1})$.
    Then, there exists $\bar{k} \in \N$ such that
    \begin{equation}\label{eq:error_bound_thm_var}
        e_k \leq X_k \left(\tilde{e}_{\bar{k}} + 2 \sum_{j=\bar{k}}^{k-1}  \frac{q_j \tau_j^2 D_{m_j}}{X_{j+1}}\right)
        \quad \forall \, k\geq \bar{k},
    \end{equation}
    with
    \begin{equation*}
        X_k = \prod_{j=\bar{k}}^k x_j, \quad x_j :=  1 - \alpha \tau_j,
    \end{equation*}
    {\color{black} and the initial error $\tilde{e}_{\bar{k}}$ defined as in \eqref{eq:init-err-notation}}.
    \label{thm:theorem_rate_variable}
\end{lmm}

\begin{proof}
    Let us show that \eqref{eq:X_k_rec_rel} holds with
    \begin{equation}\label{eq:prod:ansatz}
        X_k = \prod_{j=\bar{k}}^k x_j, \quad x_j = \abs{1 - \gamma \tau_j},
    \end{equation}
    for some $\gamma > 0$ {\color{black} which will be specified later} and for $\bar{k} \geq 0$ sufficiently large.
    Then, the result will follow from \Cref{thm:rec_X}.

    Substituting the ansatz \eqref{eq:prod:ansatz} in \eqref{eq:X_k_rec_rel} and dividing by $X_k > 0$ gives
    \begin{equation}
        x_{k+1} \geq a_k + b_k \sum_{i=1}^{s_{k}} \prod_{j=k-i+1}^k x_j^{-1}.
        \label{eq:variable_proof_to_prove_simplified}
    \end{equation}
    Let now $\bar{k}$ be such that
    \begin{equation*}
        \bar{k} + 1 - s_{\bar{k}} > (\gamma \tau)^{\frac{1}{1-t}},
    \end{equation*}
    which is always possible since $s_k$ grows slower than $k$.
    This condition is enough to ensure $\gamma \tau_{k-i+1} < 1$ for all $k \geq \bar{k}$ and $i \leq s_k$.
    With such choice of $\bar{k}$, we can then estimate
    \begin{equation}\label{eq:log-X-bound}
        \log \left( \prod_{j=k-i+1}^k x_j^{-1} \right) = - \sum_{j=k-i+1}^k \log(1- \gamma \tau_j) \leq -\frac{\log(1-\gamma \tau_{k-i+1})}{\gamma \tau_{k-i+1}} \sum_{j=k-i+1}^k \gamma \tau_j.
    \end{equation}
    Let us further assume $\bar{k}$ large enough to ensure that, say,
    \begin{equation} \label{eq:second-condition-on-bark}
        -\frac{\log(1-\gamma \tau_{k-i+1})}{\gamma \tau_{k-i+1}} \leq 2
        \qquad \forall \, k \geq \bar{k}, \, i\leq s_k,
    \end{equation}
    which is true when
    \begin{equation*}
        \bar{k} + 1 - s_{\bar{k}} > (0.79 \gamma \tau)^{\frac{1}{1-t}},
    \end{equation*}
    where the factor $0.79$ was found numerically.
    Then, defining $T_k := \sum_{j=0}^k \tau_j$, \eqref{eq:second-condition-on-bark} yields
    \begin{equation*}
        \prod_{j=k-i+1}^k x_j^{-1} \leq e^{2\gamma (T_k-T_{k-i})},
        \quad \forall \, k \geq \bar{k}, \; i \leq s_k,
    \end{equation*}
    which, in turn, leads to the following inequality
    \begin{equation*}
        \begin{split}
            \sum_{i=1}^{s_k} \prod_{j=k-i+1}^k x_j^{-1} & \leq e^{2\gamma T_k} \sum_{i=1}^{s_k} e^{-2 \gamma T_{k-i}}       \\
                                                        & = e^{2\gamma T_k} \sum_{i=k-s_k}^{k-1} e^{-2 \gamma T_{i}}        \\
                                                        & \leq e^{2\gamma T_k} \sum_{i=k-s_k}^{\infty} e^{-2 \gamma T_{i}}.
        \end{split}
    \end{equation*}
    We aim to bound the terms above exploiting the choices of ${\color{black} \tau_k} = \tau(k+1)^{t-1}$ with $t \in [0,1)$ and $s_k = \ceil{s(k+1)^{1-t}}$.
    In order to do this, we will repeatedly use the fact that, for a decreasing function $f:\R\to\R$, it holds for any two integers $a,b \in \mathbb{Z}$
    \begin{equation*}
        \int_{a}^{b+1}f(x)dx \leq \sum_{i=a}^b f(i) \leq \int_{a-1}^b f(x)dx.
    \end{equation*}
    We now distinguish the two cases $t \in (0,1)$ and $t=0$.

    \paragraph{Case $t \in (0,1)$.}
    We have
    \begin{equation*}
        T_k - T_i
        = \sum_{j=i+1}^k \tau_j
        \leq \int_{i}^{k} \tau(x+1)^{t-1}dx = \frac{\tau}{t}((k+1)^t - (i+1)^t),
    \end{equation*}
    which implies that
    \begin{equation*}
        e^{2 \gamma (T_k - T_i)} \leq \exp \left( \frac{2 \gamma \tau}{t} ((k+1)^t - (i+1)^t) \right),
    \end{equation*}
    and thus
    \begin{equation*}
        e^{2 \gamma T_{k}} \sum_{i=k-s_k}^{\infty} e^{-2 \gamma T_{i}} \leq \exp \left( \frac{2 \gamma \tau}{t} (k+1)^t \right) \int_{k-s_k-1}^{\infty} \exp \left( -\frac{2 \gamma \tau}{t} (x+1)^t \right)dx.
    \end{equation*}
    Using \Cref{lemma_gamma_int}, we have a closed form for the integral in the right-hand side of the above equation, hence, dividing further by $s_k$ and using $s_k \geq s(k+1)^{1-t}$, we have that
    \begin{equation*}
        \frac{e^{2 \gamma T_{k}}}{s_k} \sum_{i=k-s_k}^{\infty} e^{-2 \gamma T_{i}}
        \leq \frac1s \left( \frac{t}{2\gamma \tau} \right)^{\frac{1}{t}} (k+1)^{t-1}\exp \left( \frac{2 \gamma \tau}{t} (k+1)^t \right) \frac{\Gamma \left( \frac{1}{t},\left( \frac{2 \gamma \tau}{t} \right)(k-s_k)^t \right)}{t}.
    \end{equation*}
    Using the fact that $\Gamma(r,x) \sim x^{r-1} e^{-x}$ as $x \to \infty$ and $(k-s_k)^t \sim (k+1)^t - s t$ as $k \to \infty$, one can show that the right-hand side converges to $\frac{e^{2 \gamma \tau s}}{2\gamma \tau s}$ as $k$ goes to infinity and, therefore, there exists $\tilde{k} = \tilde{k}(\tau,s,t, \gamma) \in \N$ such that, for all $k \geq \tilde{k}$,
    \begin{equation*}
        \frac{1}{s_k}\sum_{i=1}^{s_k} \prod_{j=k+1-i}^k x_j^{-1} \leq \frac{e^{2 \gamma \tau s}}{\gamma \tau s}.
    \end{equation*}
    Assuming $\bar{k} \geq \tilde{k}$, equation \eqref{eq:variable_proof_to_prove_simplified} is satisfied for all $k \geq \bar{k}$ if
    \begin{equation*}
        x_{k+1} \geq a_k + s_k b_k \frac{e^{2 \gamma \tau s}}{\gamma \tau s}.
    \end{equation*}
    Replacing the expressions of $a_k$ and $b_k$, the above inequality can be simplified to
    \begin{equation*}
        \gamma \tau_{k+1} \leq 2\alpha \tau_k -4  \left( 1 + 6 \frac{e^{2 \gamma \tau s}}{\gamma \tau s} \right) L^2 q_k \tau_k^2.
    \end{equation*}
    Since $\tau_{k+1} \leq \tau_k$, the above equation is satisfied if
    \begin{equation*}
        \begin{split}
             & \gamma \tau_{k} \leq 2\alpha \tau_k -4 \left( 1 + 6 \frac{e^{2 \gamma \tau s}}{\gamma \tau s} \right) L^2 q_k \tau_k^2, \\
        \end{split}
    \end{equation*}
    which can be further recast into
    \begin{equation*}
        2\alpha-\gamma \geq 4 \left( 1 + 6 \frac{e^{2 \gamma \tau s}}{\gamma \tau s} \right) L^2 q_k \tau_k.
    \end{equation*}
    Since the above needs to be true for any $k \geq \bar{k}$, using the assumption $q_k \tau_k \to 0$ as $k \to \infty$ and the choice $\gamma = \alpha$, the inequality is satisfied taking $\bar{k}$ large enough.

    \paragraph{Case $t = 0$.}
    Let us assume $2\gamma\tau \neq 1$.
    We have
    \begin{equation*}
        T_k - T_i
        \leq \int_{i}^k \tau (x+1)^{-1} d x
        = \tau (\log(k+1) - \log(i+1)),
    \end{equation*}
    hence,
    \begin{equation*}
        e^{2 \gamma T_{k}} \sum_{i=k-s_k}^{k-1} e^{-2 \gamma T_{i}}
        \leq e^{2\gamma \tau \log(k+1)} \int_{k-s_k-1}^{k} e^{-2\gamma \tau \log(x+1)}dx = \frac{(k+1)^{2 \gamma \tau}}{1-2\gamma \tau}\left[ (k+1)^{1-2\gamma \tau}-(k-s_k)^{1-2\gamma \tau} \right],
    \end{equation*}
    which, in turn, implies
    \begin{equation*}
        \sum_{i=1}^{s_k} \prod_{j=k-i+1}^k x_j^{-1}
        \leq \frac{(k+1)^{2 \gamma \tau}}{1-2\gamma \tau} \left[ (k+1)^{1-2\gamma \tau}-(k-s_k)^{1-2\gamma \tau} \right].
    \end{equation*}
    {\color{black}
    We divide now by $s_k = \ceil{s(k+1)} \geq s(k+1)$
    and set $R_k := \frac{k+1}{k-s_k} > 1$
    to obtain
    \begin{equation*}
        \begin{split}
            \frac{1}{s_k}\sum_{i=1}^{s_k} \prod_{j=k+1-i}^k x_j^{-1}
             & \leq
            \frac{1}{s(1-2\gamma \tau)}
            \left[1-R_k^{2\gamma\tau-1}\right] \\
             & =
            \frac{\abs{1-R_k^{2\gamma\tau-1}}}{s\abs{1-2\gamma\tau}}.
        \end{split}
    \end{equation*}
    Moreover, choosing $\bar{k}$ such that $(1-s)(\bar{k}+1) \geq 4$, we have, for all $k\geq \bar{k}$,
    \begin{equation*}
        k-s_k \geq \frac{1-s}{2}(k+1),
        \qquad
        R_k \leq \frac{2}{1-s}.
    \end{equation*}
    Therefore, if $2\gamma\tau<1$, then $R_k^{2\gamma\tau-1}\leq 1$, while if $2\gamma\tau>1$, then
    \begin{equation*}
        R_k^{2\gamma\tau-1}
        \leq
        \left(\frac{2}{1-s}\right)^{2\gamma\tau-1}.
    \end{equation*}
    In both cases,
    \begin{equation*}
        \abs{1-R_k^{2\gamma\tau-1}}
        \leq
        1+\left(\frac{2}{1-s}\right)^{2\gamma\tau-1}
        =
        1+\left(\frac{1-s}{2}\right)^{1-2\gamma\tau}.
    \end{equation*}
    Thus,
    \begin{equation*}
        \frac{1}{s_k}\sum_{i=1}^{s_k} \prod_{j=k+1-i}^k x_j^{-1}
        \leq
        \frac{1}{s\abs{1-2\gamma \tau}}
        \left[
            1+\left(\frac{1-s}{2}\right)^{1-2\gamma \tau}
            \right].
    \end{equation*}
    }
    Equation \eqref{eq:variable_proof_to_prove_simplified} is hence satisfied if
    \begin{equation*}
        x_{k+1} \geq a_k + s_k b_k \frac{1}{s\abs{1-2\gamma \tau}}\left[  1 + \left( \frac{1-s}{2} \right)^{1-2 \gamma \tau} \right]
    \end{equation*}
    The above can be simplified into
    \begin{equation}
        \gamma \tau_{k+1} \leq 2 \alpha \tau_k - 4L^2 \left( 1 + 6 \frac{1}{s\abs{1-2\gamma \tau}}\left[  1 + \left( \frac{1-s}{2} \right)^{1-2 \gamma \tau} \right] \right) q_k \tau_k^2.
        \label{eq:proof_t=0:gamma_tau_cond}
    \end{equation}
    Since $\tau_k > \tau_{k+1}$, \eqref{eq:proof_t=0:gamma_tau_cond} is satisfied if
    \begin{equation*}
        \begin{split}
             & \gamma \tau_{k} \leq 2 \alpha \tau_k - 4L^2 \left( 1 + 6 \frac{1}{s\abs{1-2\gamma \tau}}\left[  1 + \left( \frac{1-s}{2} \right)^{1-2 \gamma \tau} \right] \right) q_k \tau_k^2. \\
        \end{split}
    \end{equation*}
    We then obtain the condition
    \begin{equation*}
        2\alpha - \gamma \geq 4L^2 \left( 1 + 6 \frac{1}{s\abs{1-2\gamma \tau}}\left[  1 + \left( \frac{1-s}{2} \right)^{1-2 \gamma \tau} \right] \right) q_k \tau_k
    \end{equation*}
    for $k \geq \bar{k}$, which is satisfied taking $\bar{k}$ large enough with the choice $\gamma = \alpha$ since we assumed $q_k \tau_k \to 0$ as $k \to \infty$.
\end{proof}

\begin{thrm} \label{thm:conv-var}
    Let $\{ \bbV_m \}_{m \geq 0}$ be a sequence of nested approximation spaces with $\mathrm{dim}(\bbV_m) = m$ and $\overline{\bbV}_{\infty} = L^2_\rho(\Gamma;U)$.
        {\color{black} Let further \Cref{asmp:str-cvx,asmp:lip-cont} be satisfied, $e_k$ be the squared error defined in \eqref{eq:error-notation}, $\tau_k > 0$ be the step size, $m_k$ be the approximation space size, $s_k$ be the memory size satisfying the sampling inequality \eqref{eq:samp-ineq}, and $q_k = \sup_{y\in\Gamma} w_k(y)$.}
    Assume further that the measure $\mu_k$ is quasi-optimal for $\bbV_{m_k}$ in the sense of \eqref{eq:samp-meas-quasi-opt} with constant $C_\mu \geq 1$ independent of $k$ ($\mu_k$ is optimal if $C_\mu=1$), and that $q_k \leq q m_k^\xi$ for constants $q>0$ and $\xi \in [0,1]$.
    \begin{enumerate}
        \item{(Exponential projection error decay)} \label{case:exp-dec}
              If there exist constants $C,\eta>0$ such that
              \begin{equation} \label{eq:exp-decay-asmp}
                  \min_{v \in \bbV_m} \norm{\nabla g(u^*,\cdot) - v}_{L^2_\rho(\Gamma;U)} \leq C e^{-\eta m},
              \end{equation}
              choose $s>0$, $t = \frac{1}{3}$, $\tau > \frac{\sqrt{\eta (1-\log(2)) s}}{2^{\frac13} 3 \alpha}$, and
              \begin{equation}\label{eq:mi-exp-dec}
                  m_k = \floor{\frac{1}{2}\sqrt{\frac{(1-\log(2))}{\eta C_\mu}} s_k^{\frac12}}.
              \end{equation}
              Then, there exist $\bar{k} \geq 0$ and a constant $C'>0$ such that
              \begin{equation*}
                  e_k \leq C' (k+1)^{\frac13\xi} e^{- \hat{\eta} (k+1)^{\frac{1}{3}}} \quad \forall \, k \geq \bar{k},
              \end{equation*}
              with $\hat{\eta} = \frac{1}{\sqrt[3]{2}}\sqrt{\frac{\eta (1-\log(2)) s}{C_\mu}}$.
        \item{(Algebraic projection error decay)} \label{case:alg-dec} If there exist constants $C,\eta>0$ such that
              \begin{equation}\label{eq:alg-decay-asmp}
                  \min_{v \in \bbV_m} \norm{\nabla g(u^*,\cdot) - v}_{L^2_\rho(\Gamma;U)} \leq C m^{-\eta},
              \end{equation}
              choose $s \in (0,1)$, $t=0$, $\tau > \frac{2 \eta + 1 - \xi}{\alpha }$, and $m_i \in \N$ such that
              \begin{equation}\label{eq:mi-alg-dec}
                  m_k = \floor{\frac{\kappa}{C_\mu} \frac{s_k}{\log(s_k)}},
                  \qquad \kappa = \frac{1-\log(2)}{2(1+2 \eta)}.
              \end{equation}
              Then, there exist $\bar{k} \geq 0$ and a constant $C'>0$ such that
              \begin{equation*}
                  e_k \leq C' \log(k+1)^{2\eta+\xi} (k+1)^{- 2 \eta - 1 + \xi} \quad \forall \, k \geq \bar{k}.
              \end{equation*}
    \end{enumerate}
\end{thrm}

\begin{proof}
    Notice first that in both cases \ref{case:exp-dec} and \ref{case:alg-dec} $q_k \leq q m_k^\xi$ implies $q_k = o(\tau_k^{-1})$, so that \Cref{thm:conv_var1} applies.
    Let $k \geq \bar{k}$, then,
    using the definition of $X_k$ and similar computations to the ones in \eqref{eq:log-X-bound}, we can bound
    \begin{equation*}
        X_k \leq e^{\alpha (T_{\bar{k}-1} - T_k)},
        \quad \frac{X_k}{X_{i+1}} \leq e^{\alpha (T_{i} - T_k)},
    \end{equation*}
    which, combined with the error bound \eqref{eq:error_bound_thm_var}, gives
    \begin{equation}
        e_k \leq \tilde{e}_{\bar{k}} e^{\alpha (T_{\bar{k}-1} - T_k)} + 2 \sum_{i=\bar{k}}^{k-1}  q_i \tau_i^2 D_{m_i}  e^{\alpha (T_{i} - T_k)}.
        \label{eq:result_from_thm_var}
    \end{equation}
    Furthermore, from \Cref{thm:theorem-approx-lstsq}, if $C_\mu m_i \leq \frac{1-\log(2)}{2(1+r_i)}\frac{s_i}{\log(s_i)}$,
    \begin{equation*}
        D_{m_i} \leq \left( 1 + 2\epsilon(s_i) \right)e_{m_i}(\nabla g(u^*, \cdot))^2+ 2\norm{\nabla g(u^*, \cdot)}_{L^2_\rho(\Gamma;U)}^2s_i^{-r_i},
        \qquad i \geq \bar{k},
    \end{equation*}
    where $e_{m_i}(\nabla g(u^*, \cdot))$ denotes the best approximation error and $\epsilon(s_i)$ is defined as in \Cref{thm:theorem-approx-lstsq}.

    \paragraph{Exponential error decay.}
    Let us consider first the case of exponential decay of the gradient error approximation \eqref{eq:exp-decay-asmp}, in which case the assumption implies $e_{m_i}(\nabla g(u^*, \cdot)) \leq C e^{-\eta m_i}$, so that for $i \geq \bar{k}$
    \begin{equation}\label{eq:exp-err-bound-mi}
        \begin{split}
            D_{m_i} & \leq C \left( 1 + \epsilon(s_i) \right) e^{-2\eta m_i}+ 2\norm{\nabla g(u^*, \cdot)}_{L^2_\rho(\Gamma;U)}^2s_i^{-r_i}                    \\
                    & = e^{2 \eta} C \left( 1 + \epsilon(s_i) \right) e^{-2\eta (m_i+1)}+ 2s_i\norm{\nabla g(u^*, \cdot)}_{L^2_\rho(\Gamma;U)}^2s_i^{-(r_i+1)} \\
        \end{split}
    \end{equation}
    Let us choose $m_i \in \N$ such that
    \begin{equation*}
        m_i \leq \frac{1-\log(2)}{2 C_\mu (1+r_i)} \frac{s_i}{\log(s_i)} \leq m_i + 1,
    \end{equation*}
    which ensures \eqref{eq:exp-err-bound-mi} holds, and, in order to balance the terms in the right-hand side of \eqref{eq:exp-err-bound-mi}, let $r_i$ satisfy
    \begin{equation*}
        \frac{\eta (1-\log(2))}{C_\mu (1+r_i)} \frac{s_i}{\log(s_i)} = (1+r_i) \log(s_i),
    \end{equation*}
    that is $1+r_i = \tilde{\eta} \frac{s_i^{\frac{1}{2}}}{\log(s_i)}$ with $\tilde{\eta} = \sqrt{\frac{\eta (1-\log(2))}{C_\mu}}$.
    Note that, since $i \mapsto \frac{s_i^{\frac{1}{2}}}{\log(s_i)}$ is increasing for $i \geq \bar{i}$ with $s_{\bar{i}} \geq 2$, letting $\bar{k}$ be large enough ensures $r_i > 0$ for all $i \geq \bar{k}$.
    These choices lead to \eqref{eq:mi-exp-dec}.
    Then, using
    \begin{equation}\label{eq:Ti-Tk-bound-exp}
        T_i - T_k
        = - \sum_{j=i+1}^k \tau_j \leq - \int_{i+1}^{k+1} \tau (x+1)^{t-1} dx = \frac{\tau}{t}((i+2)^t - (k+2)^t),
    \end{equation}
    we can bound the sum in the right-hand side of \eqref{eq:result_from_thm_var} as
    \begin{equation*}
        \begin{split}
                 & \sum_{i=\bar{k}}^{k-1} q_i \tau_i^2 D_{m_i} e^{\alpha (T_{i} - T_k)}                                                                                                                                                                                                                      \\
            \leq & \sum_{i=\bar{k}}^{k-1} q_i \tau_i^2 \left( e^{2\eta} C \left( 1 + 2\epsilon(s_i) \right)e^{-2\eta (m_i+1)}+ 2 s_i \norm{\nabla g(u^*, \cdot)}_{L^2_\rho(\Gamma;U)}^2  s_i^{-(1+r_i)} \right) \exp\left(\frac{\alpha \tau}{t}((i+2)^t - (k+2)^t)\right)                                    \\
            =    & \sum_{i=\bar{k}}^{k-1} q_i \tau_i^2  \left( e^{2\eta}C(1 + 2\epsilon(s_i)) + 2 s_i \norm{\nabla g(u^*, \cdot)}_{L^2_\rho(\Gamma;U)}^2  \right) \exp\left(-\tilde{\eta} s_i^{\frac{1}{2}} + \frac{\alpha \tau}{t}((i+2)^t - (k+2)^t)\right)                                                \\
            \leq & \sum_{i=\bar{k}}^{k-1} q_i \tau_i^2  \left( e^{2\eta}C(1 + 2\epsilon(s_i)) + 2 s_i \norm{\nabla g(u^*, \cdot)}_{L^2_\rho(\Gamma;U)}^2  \right) \exp\left(-\frac{\tilde{\eta} s^{\frac{1}{2}}}{2^{\frac{1-t}{2}}} (i+2)^{\frac{1-t}{2}} + \frac{\alpha \tau}{t}((i+2)^t - (k+2)^t)\right), \\
            \leq & \left( e^{2\eta} C(1 + 2\epsilon(s_{\bar{k}})) + 2 \norm{\nabla g(u^*, \cdot)}_{L^2_\rho(\Gamma;U)}^2  \right) \exp\left(-\frac{\alpha \tau}{t} (k+2)^t\right)                                                                                                                            \\
                 & \times \sum_{i=\bar{k}}^{k-1} q_i \tau_i^2 (1+s_i) \exp\left(-\frac{\tilde{\eta} s^{\frac{1}{2}}}{2^{\frac{1-t}{2}}} (i+2)^{\frac{1-t}{2}} + \frac{\alpha \tau}{t}(i+2)^t\right),                                                                                                         \\
        \end{split}
    \end{equation*}
    where we used the fact that $s_i \geq s(i+1)^{1-t}$ and $(i+2)\leq 2(i+1)$ in the second inequality, and the fact that $\epsilon(\cdot)$ is a decreasing function in the last inequality.
    We now balance the exponents in the exponential terms, that is we set $\frac{1-t}{2} = t$, which yields $t = \frac{1}{3}$.
    Then, we can write
    \begin{equation}\label{eq:exp-dec-interm-res}
        \begin{split}
                 & \sum_{i=\bar{k}}^{k-1} q_i \tau_i^2 D_{m_i} e^{\alpha (T_{i} - T_k)}                                                                                                                               \\
            \leq & \tilde{q} \tau^2 \left( C(1 + 2\epsilon(s_{\bar{k}})) + 2 \norm{\nabla g(u^*, \cdot)}_{L^2_\rho(\Gamma;U)}^2 \right) \exp\left(- 3 \alpha \tau (k+2)^{\frac{1}{3}}\right)                          \\
                 & \times \sum_{i=\bar{k}}^{k-1} \left(1+s_i\right)(i+1)^{\frac13\xi-\frac{4}{3}}  \exp\left(\left(3 \alpha \tau - \frac{\tilde{\eta} s^{\frac{1}{2}}}{2^{\frac13}}\right)(i+2)^{\frac{1}{3}}\right), \\
        \end{split}
    \end{equation}
    where we used the fact that there exists a constant $\tilde{q} > 0$ such that $q_k \leq \tilde{q} (k+1)^{\frac13\xi}$.
    Let us now choose
    \begin{equation}\label{eq:expdec-tau-s-ineq}
        \tau > \frac{\tilde{\eta} s^{\frac{1}{2}}}{2^{\frac13} 3 \alpha},
    \end{equation}
    and notice that, using \Cref{lem:end-asymp-int}, we have the asymptotic expansion
    \begin{equation*}
        \int_{\bar{k}}^{k+1} x^a e^{b x^{\frac{1}{3}}} \sim \frac{3}{b} (k+1)^{a+\frac{2}{3}} e^{b (k+1)^{\frac{1}{3}}},
    \end{equation*}
    where $a\in \R$ and $b>0$.
    Noting that $1+s_i \sim s (i+1)^{\frac23}$, upper bounding the sum in \eqref{eq:exp-dec-interm-res} with the integral from $\bar{k}$ to $k+1$, and applying the result above with $a = \frac13\xi - \frac{2}{3} < 0$ and $b = 3 \alpha \tau - 2^{-\frac13} s^{\frac{1}{2}} \tilde{\eta} >0$ we can conclude that there exists a constant $C'>0$ independent of $k$ such that
    \begin{equation*}
        \begin{split}
            \sum_{i=\bar{k}}^{k-1} q_i \tau_i^2 D_{m_i} e^{\alpha (T_{i} - T_k)}
            \leq & C' (k+1)^{\frac13\xi} e^{- 2^{-\frac13} s^{\frac{1}{2}} \tilde{\eta} (k+1)^{\frac{1}{3}}}.
        \end{split}
    \end{equation*}
    Finally, owing to \eqref{eq:Ti-Tk-bound-exp} and the choice \eqref{eq:expdec-tau-s-ineq}, we have that
    \begin{equation*}
        \tilde{e}_{\bar{k}} e^{\alpha (T_{\bar{k}-1} - T_k)}
        \leq \tilde{e}_{\bar{k}} e^{3 \alpha \tau ((\bar{k}+1)^{\frac{1}{3}} - (k+2)^{\frac{1}{3}})}
        \leq \tilde{e}_{\bar{k}} e^{3 \alpha \tau (\bar{k}+1)^{\frac{1}{3}} - 2^{-\frac13}s^{\frac{1}{2}} \tilde{\eta} (k+2)^{\frac{1}{3}}},
    \end{equation*}
    hence thesis follows from \eqref{eq:result_from_thm_var}.

    \paragraph{Algebraic error decay.}
    Let us now consider the case \eqref{eq:alg-decay-asmp} where the approximation error decays algebraically with rate $\eta$.
    With the choice \eqref{eq:mi-alg-dec}, we have
    \begin{equation*}
        D_{m_i} \leq \left( 1 + 2\epsilon(s_i) \right)e_{m_i}(\nabla g(u^*, \cdot))^2+ 2\norm{\nabla g(u^*, \cdot)}_{L^2_\rho(\Gamma;U)}^2s_i^{-2 \eta}.
    \end{equation*}
    In order to have an algebraic decay of the contraction terms in the error bound \eqref{eq:result_from_thm_var}, so that it matches the algebraic decay of the approximation error, we consider $t = 0$.
    Then, using
    \begin{equation}\label{eq:Ti-Tk-bound-alg}
        T_i - T_k = - \sum_{j=i+1}^k \tau_j \leq - \int_{i+1}^{k+1} \tau (x+1)^{-1} dx = \tau (\log(i+2) - \log(k+2))
    \end{equation}
    and $m_i+1\geq \frac{\kappa}{C_\mu} \frac{s_i}{\log(s_i)}$,
    we can bound the sum in the right-hand side of \eqref{eq:result_from_thm_var} as
    \begin{equation*}
        \begin{split}
                 & \sum_{i=\bar{k}}^{k-1} q_i \tau_i^2 D_{m_i} e^{\alpha (T_{i} - T_k)}                                                                                                                                                                                                                  \\
            \leq & (k+2)^{-\alpha \tau} \sum_{i=\bar{k}}^{k-1} q_i \tau_i^2 (i+2)^{\alpha \tau} \left( C \left( 1 + 2\epsilon(s_i) \right) 2^{2\eta} (m_i+1)^{-2\eta}+ 2\norm{\nabla g(u^*, \cdot)}_{L^2_\rho(\Gamma;U)}^2 s_i^{- 2 \eta} \right)                                                        \\
            \leq & (k+2)^{-\alpha \tau} \sum_{i=\bar{k}}^{k-1} q_i \tau_i^2 (i+2)^{\alpha \tau} \left( C \left( 1 + 2\epsilon(s_i) \right) \left(\frac{2 C_\mu}{\kappa}\right)^{2\eta} \frac{s_i^{-2\eta}}{\log(s_i)^{-2\eta}}+ 2\norm{\nabla g(u^*, \cdot)}_{L^2_\rho(\Gamma;U)}^2 s_i^{-2\eta} \right) \\
            \leq & \tilde{q} \max\left\{1, 2^{\alpha \tau - 2 \eta - 2 + \xi}\right\} s^{-2\eta} \left( C \left( 1 + 2\epsilon(s_{\bar{k}}) \right) \left(\frac{2 C_\mu}{\kappa}\right)^{2\eta} + 2\norm{\nabla g(u^*, \cdot)}_{L^2_\rho(\Gamma;U)}^2 \right)                                            \\
                 & \times (k+2)^{-\alpha \tau} \sum_{i=\bar{k}}^{k-1} \log(s_i)^{2 \eta+\xi} (i+2)^{\alpha \tau - 2 \eta - 2 + \xi},
        \end{split}
    \end{equation*}
    where we assumed $\bar{k}$ large enough so that $\frac{m_i+1}{m_i} \leq 2$ in the second inequality, and the fact that there exists a constant $\tilde{q}$ such that $q_i \leq \tilde{q} \frac{(i+1)^\xi}{\log(s_i)^{\xi}}$, $s_i \geq s(i+1)$, and
    $\left(\frac{i+2}{i+1}\right)^{\alpha \tau - 2 \eta - 2 + \xi} \leq \max\left\{1, 2^{\alpha \tau - 2 \eta - 2 + \xi}\right\}$ in the third inequality.
    Setting now
    \begin{equation}\label{eq:algdec-tau-s-ineq}
        \alpha \tau > {2 \eta + 1 - \xi}
    \end{equation}
    and using the fact that, owing to \Cref{lem:end-asymp-int} used after the change of variable $u=\log(x)$, for $a>0$ and $b>-1$ it holds $\int_{\bar{k}}^k \log(x)^a x^b dx \sim \frac{\log(k+1)^a}{b+1} (k+1)^{b+1}$ as $k \to \infty$, we obtain
    \begin{equation*}
        \begin{split}
            (k+2)^{-\alpha \tau} \sum_{i=\bar{k}}^{k-1} \log(s_i)^{2 \eta+\xi} (i+2)^{\alpha \tau - 2 \eta - 2 + \xi} \sim \frac{\log(k+1)^{2 \eta+\xi}}{\alpha \tau - 2 \eta - 1 + \xi} (k+1)^{- 2 \eta - 1 + \xi}.
        \end{split}
    \end{equation*}
    The above is sufficient to conclude that there exists a constant $C'$ independent of $k$ such that
    \begin{equation*}
        \sum_{i=\bar{k}}^{k-1} q_i \tau_i^2 D_{m_i} e^{\alpha (T_{i} - T_k)}
        \leq C' \log(k+1)^{2 \eta+\xi} (k+1)^{- 2 \eta - 1 + \xi}.
    \end{equation*}
    Finally, owing to \eqref{eq:Ti-Tk-bound-alg} and the choice \eqref{eq:algdec-tau-s-ineq}, we have that
    \begin{equation*}
        \tilde{e}_{\bar{k}} e^{\alpha (T_{\bar{k}-1} - T_k)}
        \leq \tilde{e}_{\bar{k}} e^{\alpha \tau (\log(\bar{k}+1) - \log(k+2))}
        = \tilde{e}_{\bar{k}} (\bar{k}+1)^{\alpha \tau} (k+2)^{-\alpha \tau}
        \leq \tilde{e}_{\bar{k}} (\bar{k}+1)^{\alpha \tau} (k+2)^{-2\eta - 1 + \xi},
    \end{equation*}
    hence the thesis follows from \eqref{eq:result_from_thm_var}.
\end{proof}

After analyzing the convergence properties of the proposed method, it is natural to address their computational cost.
Suppose that each evaluation of the function $g$ at a control-data pair $(u, y)$ incurs a unit cost denoted by $C_g > 0$.
When the sampling measure $\mu$ is fixed, both \Cref{alg:algo-fixed} and \Cref{algo:algo-variable} prescribe exactly one gradient evaluation per iteration.
Therefore, the total computational cost incurred after $k$ iterations, which we call $\mathrm{Work}(k)$, is precisely $C_g k$.
The following proposition shows that, on average, this favorable cost structure is preserved asymptotically even in the presence of the resampling strategy.

\begin{prpstn}\label{prop:work}
    Under the same assumptions and definitions of \Cref{thm:conv-var}, the expected cost at iteration $k$ satisfies
    \begin{equation*}
        \frac{\Exp{\mathrm{Work}(k)}}{k} \xrightarrow[k \to \infty]{} \gamma C_g,
    \end{equation*}
    where $\gamma = 1$ in Case \ref{case:exp-dec} and $\gamma = 1+s < 2$ in Case \ref{case:alg-dec}.
\end{prpstn}
\begin{proof}
    Let us consider the case where $\mu_k$ changes with the iterations and is the optimal sampling measure in the sense of \eqref{eq:opt-samp-meas} associated to the approximation space $\bbV_{m_k}$, since the case where the sampling measure is fixed is trivial (no resampling is necessary).
    Let $k \in \N$ be the iteration counter.
    Assume $k\geq \tilde{k}$, where $\tilde{k} := \inf\{i \in \N \,:\,m_i \geq 1\}$ as no resampling is needed if $m_k=0$.
    Denote $b_{k,j}$ the $j$-th Bernoulli variable sampled as a result of running \Cref{algo:resampling}.
    Then, the number of gradient evaluations at iteration $k$ is
    \begin{equation*}
        1 + \sum_{j=1}^{s_k} (1-b_{k,j}),
    \end{equation*}
    which, in expectation, gives
    \begin{equation*}
        1 + \sum_{j=1}^{s_k} \frac{m_k-m_{k-1}}{m_k} = 1 + \frac{s_k}{m_k} (m_k-m_{k-1}).
    \end{equation*}
    Letting $C_g=1$ to ease notation, the expected work at iteration $k$ coincides with the sum of the expected number of gradient evaluations up to iteration $k$, that is
    \begin{equation*}
        \Exp{\mathrm{Work}(k)}
        = k + \sum_{i=\tilde{k}}^k \frac{s_i}{m_i} (m_i-m_{i-1}).
    \end{equation*}
    Since $\frac{s_i}{m_i}$ is an increasing sequence it holds $\frac{s_i}{m_i} (m_i-m_{i-1}) \leq s_i - s_{i-1}$, hence we can bound the sum in the above by the telescoping sum $\sum_{i=\tilde{k}}^k (s_i-s_{i-1}) = s_{k} - s_{\tilde{k}-1}$.
    The thesis then follows from the definitions of $s_k$ for the two cases of \Cref{thm:conv-var}.
\end{proof}

Finally, we can combine the convergence results of \Cref{thm:conv-var} and the asymptotic behavior of the expected cost from \Cref{prop:work} to obtain the asymptotic complexity of \Cref{algo:algo-variable}.
\begin{crllr}
    Under the same assumptions and definitions of \Cref{thm:conv-var}, consider a tolerance $0<\epsilon \lesssim 1$ and define $\calW(\epsilon)$ the total expected computational work that \Cref{algo:algo-variable} incurs to achieve $J(u_k) - J(u^\ast) \leq \epsilon$.
    \begin{enumerate}
        \item In Case \ref{case:exp-dec}, for all $\delta>0$ there exist $C'>0$ such that
              \begin{equation*}
                  \calW(\epsilon)
                  \leq \left( \frac{1+\delta}{\hat{\eta}} \log\left(\frac{C'}{\epsilon}\right) \right)^3,
              \end{equation*}
              if $\epsilon$ is sufficiently small.
        \item In Case \ref{case:alg-dec}, for all $\delta>0$ there exist $C'>0$ such that
              \begin{equation*}
                  \calW(\epsilon)
                  \leq \left( \frac{C'}{\epsilon} \right)^{\frac{1}{2\eta + 1 - \xi - \delta}},
              \end{equation*}
              if $\epsilon$ is sufficiently small.
    \end{enumerate}
\end{crllr}
\begin{proof}
    Owing to \Cref{asmp:lip-cont}, we have $J(u_k) - J(u^\ast) \leq \bbE_{Y\sim \rho}[\ell(Y)] \norm{u_k-u^\ast}_{U}^2$.
    The rest follows from \Cref{thm:conv-var} and \Cref{prop:work}.
\end{proof}

\section{Numerical experiments}
\label{sec:numerics}
In this section we present numerical experiments performed in the context of PDE-constrained optimization.
Links to the code used are provided at the end of the paper.

\subsection{Diffusion equation with one parametric dimension}
\label{sec:num-one-dim}

In this subsection we study a simple family of optimal control problems governed by a
diffusion-type PDE with exactly one stochastic (parametric) degree of
freedom, while the spatial domain is the fixed square $D=[0,1]^{2}$.
More precisely, given a control $u \in L^2(D)$ and a realization $y$ of the stochastic parameter $Y \sim \mathcal{U}([-1,1])$, the state $z(\cdot; u, y)$ satisfies the equation
\begin{equation}\label{eq:pde-one-dim}
    \begin{cases}
        -\mathrm{div}(\tilde{y} \cdot  \nabla z(x;u,y)) = u(x) \qquad & \text{in } [0,1]^2          \\
        z(x;u,y) = 0 \qquad                                           & \text{on } \partial[0,1]^2,
    \end{cases}
\end{equation}
where
\begin{equation*}
    \tilde{y} = a \exp\left( \frac{(y+1)\log(b/a)}{2} \right)
\end{equation*}
for some $b > a > 0$.
One can show that the map $(u,y) \mapsto z(u,y) \in H_0^1(D)$ is well-defined for all $y\in\Gamma$ and $u\in L^2(D)$, $z(\cdot;u,\cdot) \in L^{\infty}(\Gamma;H^1_0(D))$ for all $u\in L^2(D)$, and $u \mapsto z(\cdot; u, y) \in H^1_0(D)$ is Fréchet differentiable for all $y \in \Gamma = [-1,1]$ (see, e.g., \cite[\S 6]{kouri2018optimization}).
We aim to solve the minimization problem
\begin{equation}\label{eq:min-prob-one-dim}
    \min_{u \in L^2(D)} J(u) := \mathbb{E}\left[ \frac{1}{2}\norm{z(u,Y)-z_d}_{L^2(D)}^2 + \frac{\beta}{2}\norm{u}_{L^2(D)}^2\right] ,
\end{equation}
where $z(u,y)$ is short for $z(\cdot;u,y)$ which solves the PDE \eqref{eq:pde-one-dim} and $z_d$ is a target state in $L^2(D)$.
Note that, for this simple problem, we can get rid of the uncertainty and compute analytically the minimum of the above quantity when $z_d$ is an eigenfunction of the Laplacian.
Indeed, the following proposition holds.
\begin{prpstn}\label{prop:one-dim}
    Let $z_d$ be an eigenfunction of the Laplacian $- \Delta$ associated to eigenvalue $\lambda > 0$ and define
    \begin{equation*}
        g(u,y) := \frac{1}{2}\norm{z(u,y)-z_d}_{L^2(D)}^2 + \frac{\beta}{2}\norm{u}_{L^2(D)}^2,
    \end{equation*}
    so that $J(u) = \Exp{g(u,Y)}$.
    Then, the minimizer of \eqref{eq:min-prob-one-dim} is proportional to $z_d$ and is given by
    \begin{equation}\label{eq:u-ast-prop}
        u^\ast = \left[ \frac{\alpha \lambda}{1+\tilde{\beta}\lambda^{2}} \right] z_d,
    \end{equation}
    where
    \begin{equation}\label{eq:tildey-one-dim}
        \alpha = \frac{\Exp{\tilde{Y}^{-1}}}{\Exp{\tilde{Y}^{-2}}},
        \quad \tilde{\beta}= \frac{\beta}{\Exp{\tilde{Y}^{-2}}}.
    \end{equation}
    Moreover, the parameter-to-gradient map $y \mapsto \nabla g(u,y) \in L^2(D)$ is analytic for all $u \in L^2(D)$.
    \label{prop:proportional_prop}
\end{prpstn}
\begin{proof}
    See \Cref{sec:app-prop-one-dim}.
\end{proof}
Furthermore, the quantities $\Exp{\tilde{Y}^{-1}}$ and $\Exp{\tilde{Y}^{-2}}$ can be explicitly computed via direct integration, indeed
\begin{equation*}
    \begin{split}
        \Exp{\tilde{Y}^{-1}} = \frac{b-a}{ab\log(b/a)}\quad \text{and} \quad \Exp{\tilde{Y}^{-2}}=\frac{b^2-a^2}{2a^2b^2\log(b/a)}.
    \end{split}
\end{equation*}
This allows us to easily assess the performance of our algorithm.

In practice, we chose $z_d = \sin(\pi x_1) \sin(\pi x_2)$ and discretized the PDE on a finite element space $V_h$ ($\bbP_1$ finite elements on a regular mesh with $81$ nodes) and obtained gradients of the discretized functional via the adjoint method, {\color{black} which requires two PDE solves per (stochastic) gradient evaluation}.
Furthermore, we considered the approximation space $\bbV_m = \calV_m \otimes V_h$ where $\calV_m$ is the $m$-dimensional space of polynomials of degree less than or equal to $m-1$.
Owing to the smoothness of $y \mapsto \nabla g(u,y)$ given by \Cref{prop:one-dim}, the spaces $\bbV_m$ achieve exponential decay of the gradient approximation.
Since the distribution on the parameter space is uniform, we use Legendre polynomials to form an orthonormal basis of $\calV_m$.
The sampling measure we employ is approximation space-independent and is fixed to the arcsine distribution on $\Gamma$, that is $d\mu(y) = \frac{2}{\pi \sqrt{1-y^2}} \frac{dy}{2}$, which entails the weight supremum
\begin{equation*}
    q = \sup_{y \in \Gamma} \frac{\pi \sqrt{1-y^2}}{2}
    = \frac{\pi}{2} \approx 1.5.
\end{equation*}

\begin{figure}[ht]
    \centering
    \begin{subfigure}[t]{0.5\textwidth}
        \includegraphics[width=\linewidth]{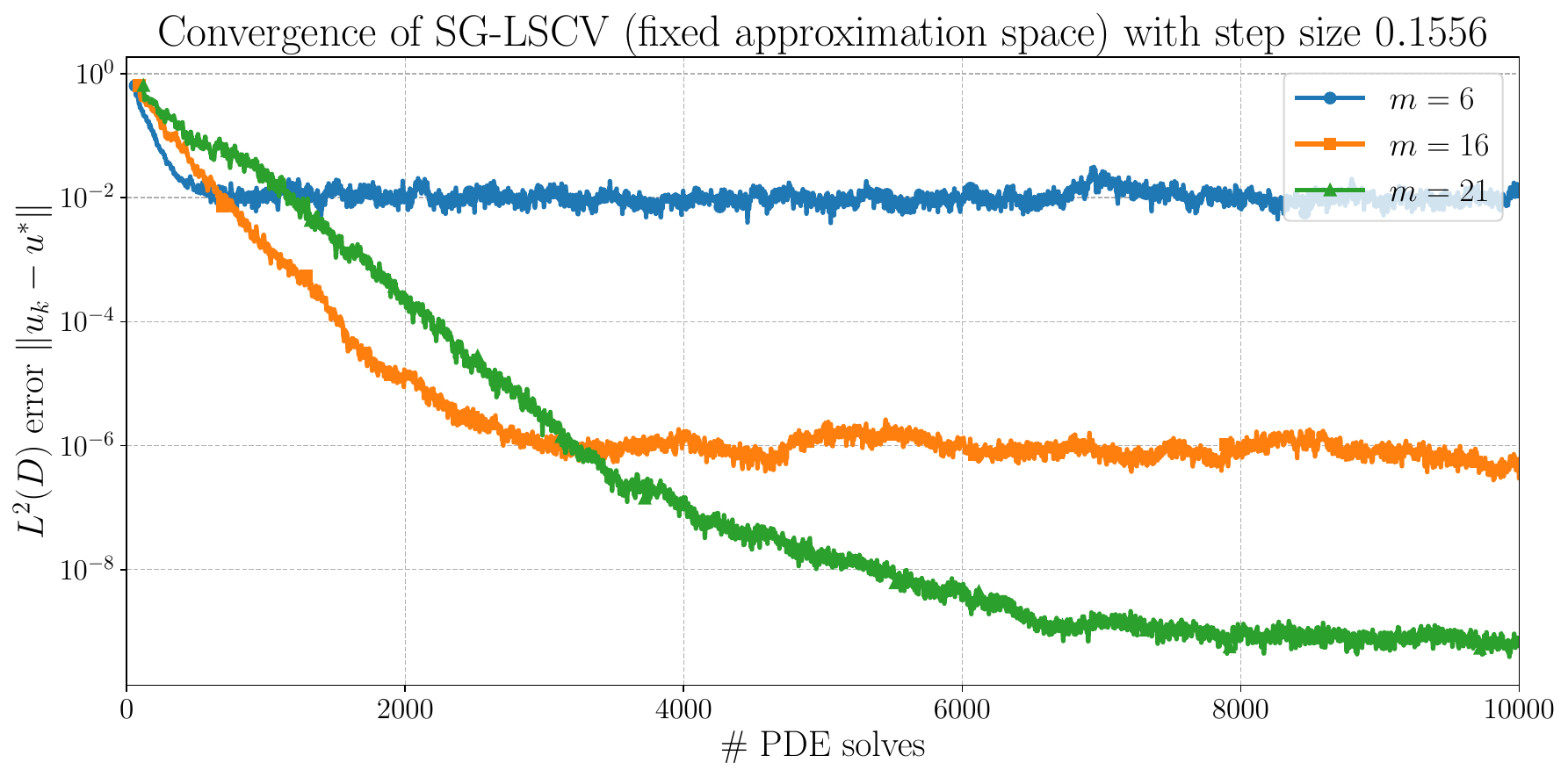}
        \caption{Comparison for $\tau=0.1556$ (optimized for $m=6$).}
        \label{subfig1:0.1556}
    \end{subfigure}%
    \hfill
    \begin{subfigure}[t]{0.5\textwidth}
        \includegraphics[width=\linewidth]{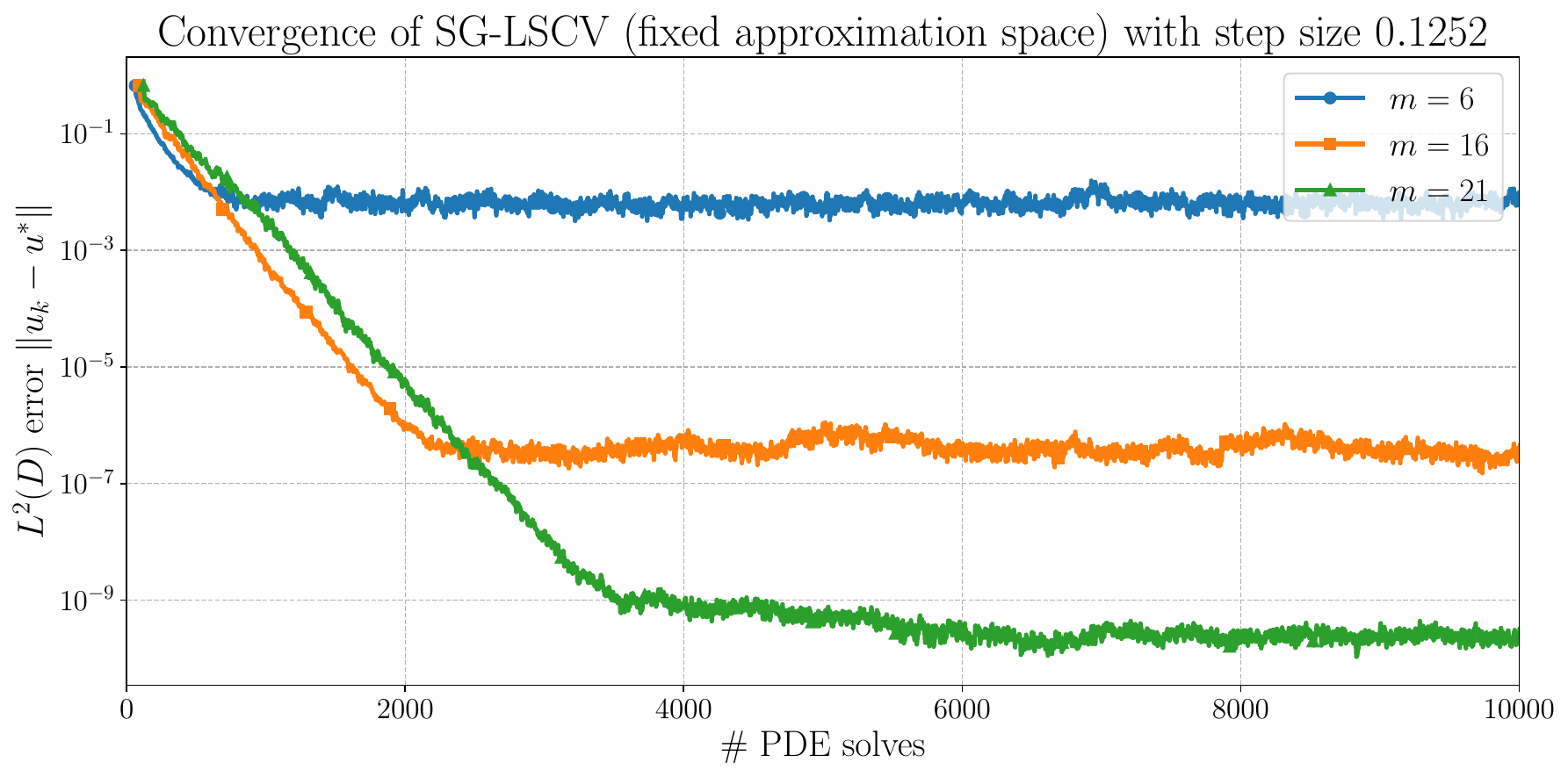}
        \caption{Comparison $\tau=0.1252$ (optimized for $m=16$).}
        \label{subfig1:0.1252}
    \end{subfigure}
    \\
    \begin{subfigure}[t]{0.5\textwidth}
        \includegraphics[width=\linewidth]{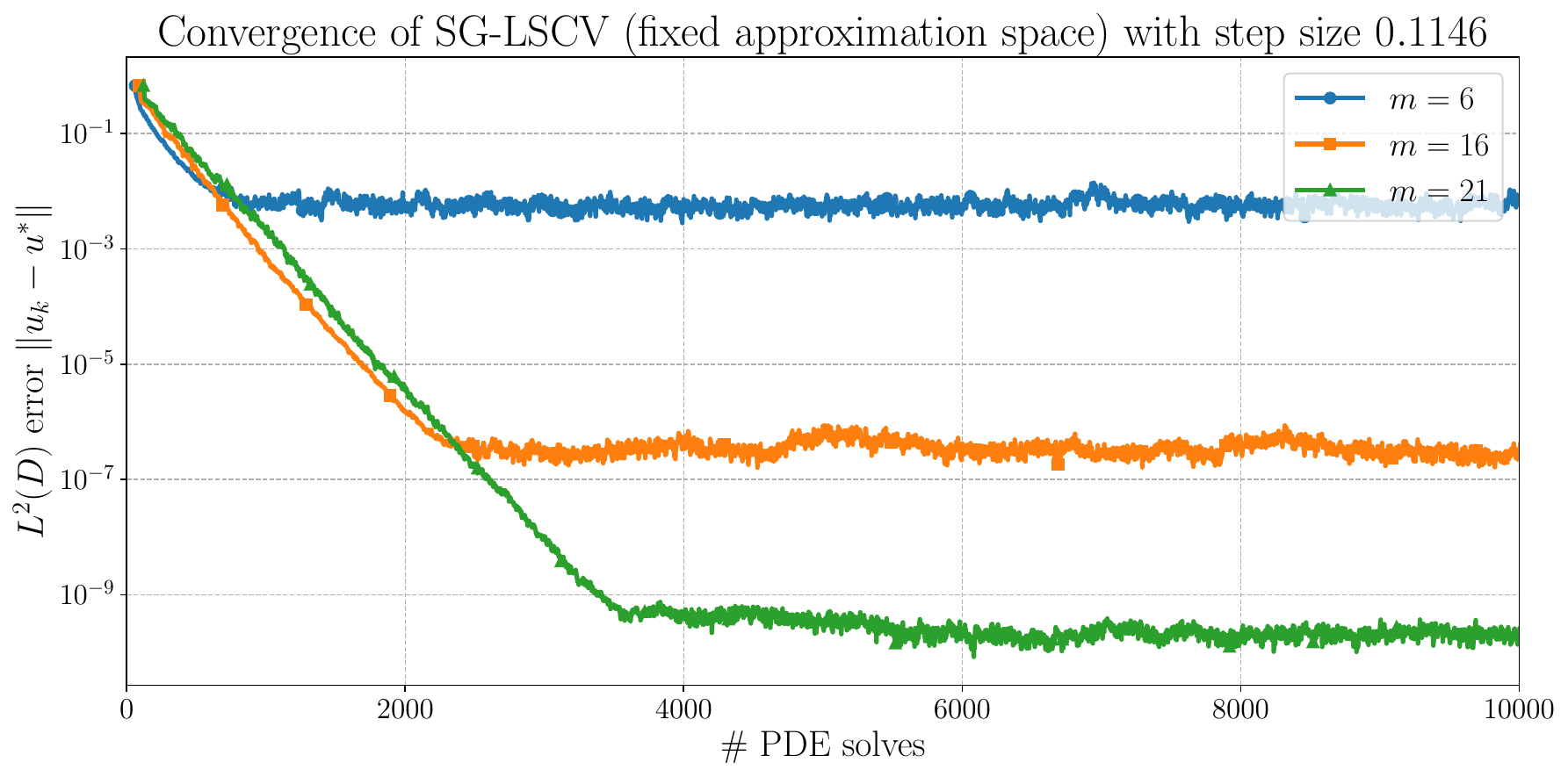}
        \caption{Comparison $\tau=0.1146$ (optimized for $m=21$).}
        \label{subfig1:0.1146}
    \end{subfigure}%
    \hfill
    \begin{subfigure}[t]{0.5\textwidth}
        \includegraphics[width=\linewidth]{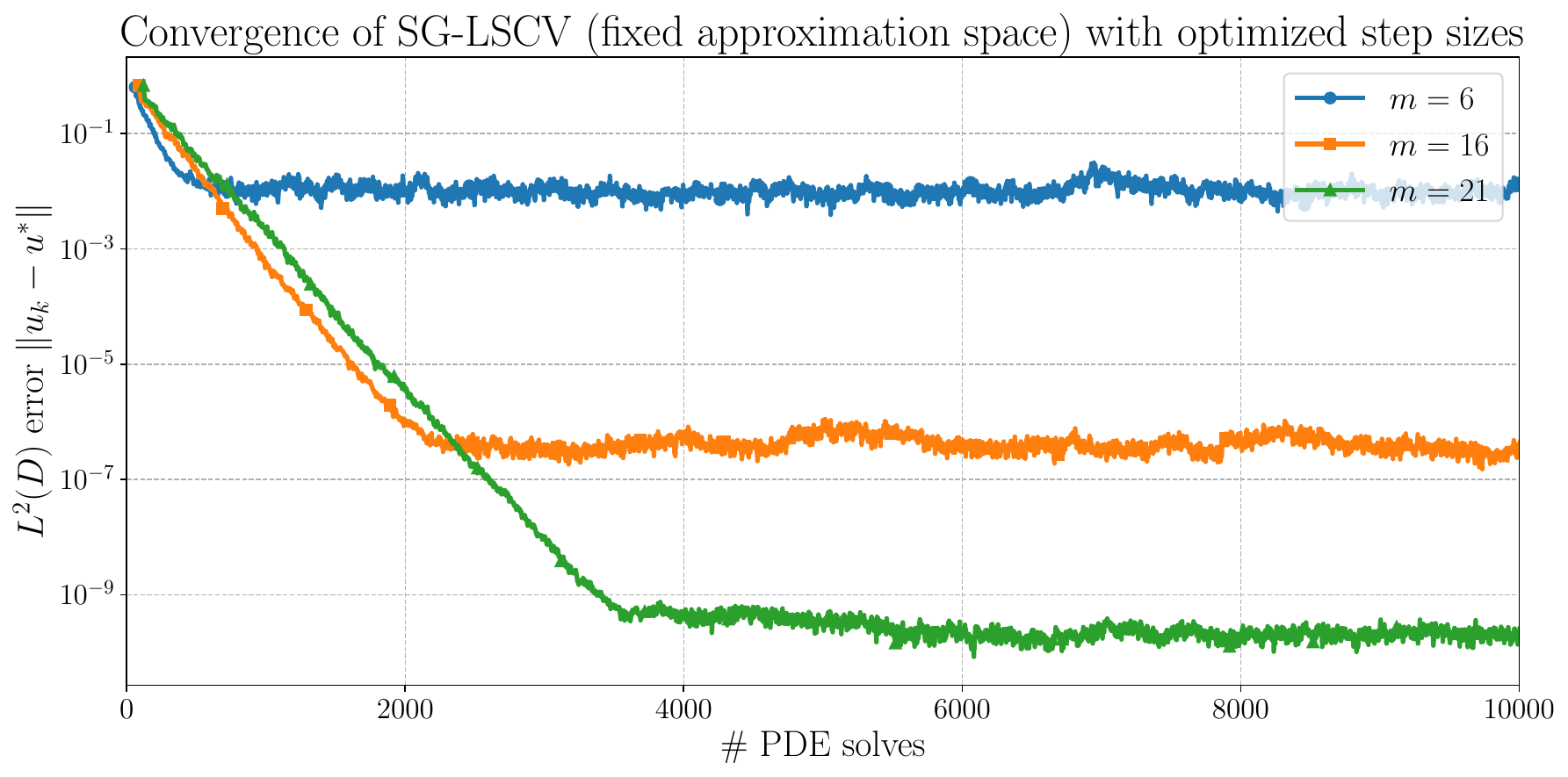}
        \caption{Comparison with optimized step sizes.}
        \label{subfig1:beststep}
    \end{subfigure}
    \caption{Convergence behavior of the \algname algorithm with fixed polynomial approximation space of dimension $m \in \{6,16,21\}$, for different choices of the step size.
            {\color{black} The algorithm uses one gradient evaluation per iteration, which corresponds to two PDE solves.}
    }
    \label{fig:fixed_sglscv_comparisons}
\end{figure}

{\color{black} We show in Figure \ref{fig:fixed_sglscv_comparisons} the convergence plots of $\norm{u_k-u^*}_{L^2(D)}$ versus the number of PDE solves for the \algname algorithm with fixed polynomial approximation space of dimension $m=6,16,21$.}
We report the geometric average over $40$ runs of each experiment.
In Figure \ref{fig:fixed_sglscv_comparisons}, we emphasize the effect of adapting the step size in the fixed \algname algorithm.
Specifically, we optimize the step size for $m=6,16,21$ and obtain $\tau=0.1556,0.1252,0.1146$, respectively.
This confirms the theoretical result emphasized in \Cref{rem:rate-large-s} that the optimal step size should shrink as $m$ increases.
This observation is further confirmed by the comparisons for fixed step size in \Cref{sub@subfig1:0.1556,sub@subfig1:0.1252,sub@subfig1:0.1146}, which show that the algorithms with smaller $m$ benefit from a larger step size and vice-versa.
Furthermore, it is clear that all algorithms feature two regimes, one of exponential convergence to the minimum and one of stagnation (no convergence), as predicted by \Cref{thm:conv-fixed}.
The exponential convergence rate is faster for smaller $m$'s at the price of a higher stagnation error, and vice-versa.
The latter observation still holds true even when we plot the algorithms with the respective optimized step size in \Cref{sub@subfig1:beststep}, meaning that this behavior is not an artifact stemming from a potentially non-optimal step size.

\begin{figure}[ht]
    \centering
    \includegraphics[width=\textwidth]{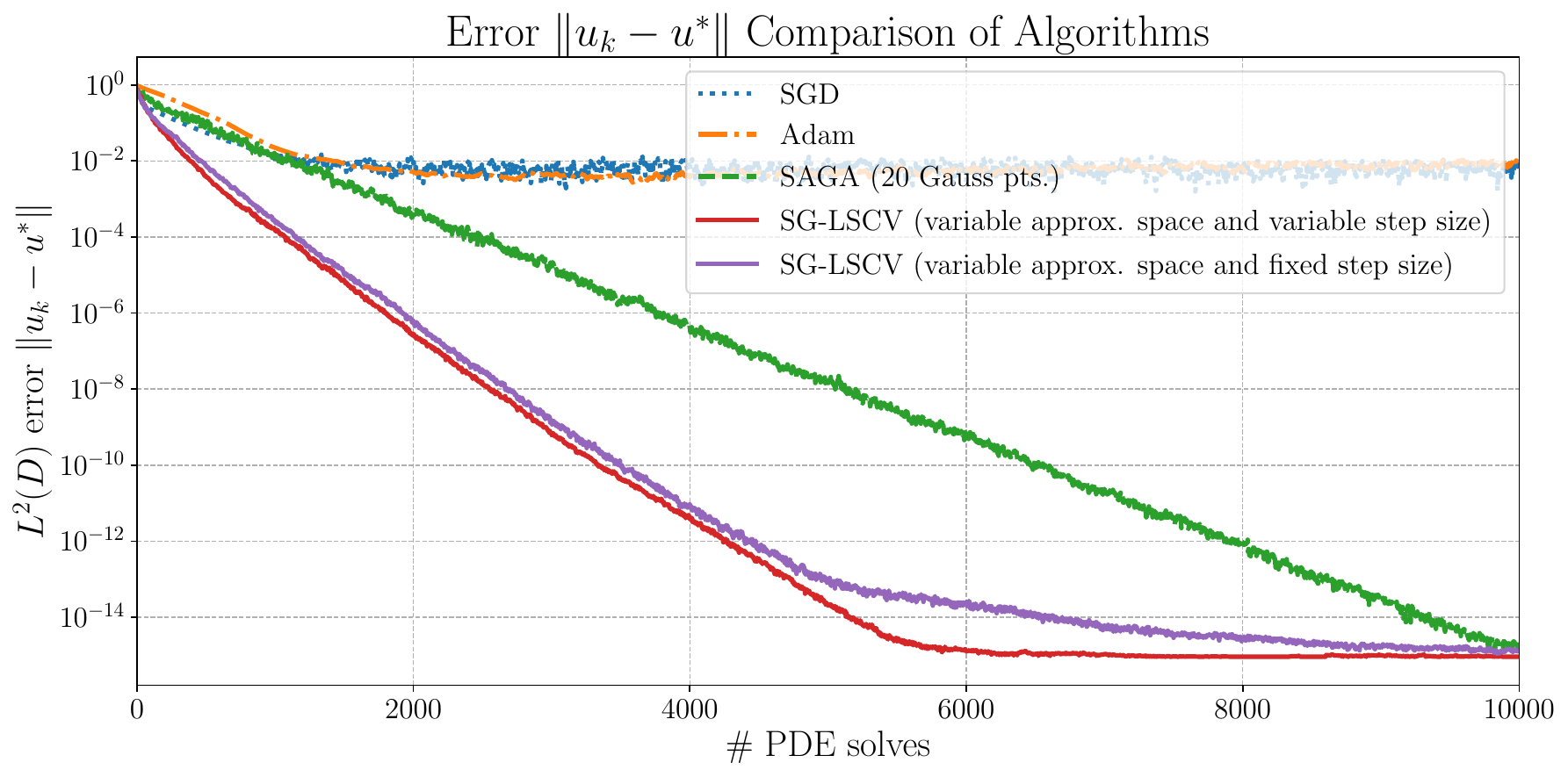}
    \caption{Comparison of SGD, Adam, SAGA and two versions of variable \algname.
    {\color{black} Each algorithm uses one gradient evaluation per iteration, which corresponds to two PDE solves.}}
    \label{fig:diff_alg_comparison}
\end{figure}

\begin{rmrk}
    Concerning the SAGA algorithm, we use importance sampling by taking a Gauss-Legendre quadrature in order to approximate the
    expected value.
\end{rmrk}

We show in \Cref{fig:diff_alg_comparison} the convergence plots of $\norm{u_k-u^*}_{L^2(D)}$ versus the number of {\color{black} PDE solves} for SGD, Adam, SAGA, and \algname with variable approximation space, both with a fixed step size and a memory-dependent step size (as \Cref{thm:conv-var} prescribes).
Note that SAGA was used on a discretization of the functional via a $20$-point Gaussian quadrature, which makes the comparison fair since it exploits the parameter regularity of the solution map.
The step size used for SAGA is tuned and equal to $\tau = 0.0769$,  while SGD was run with a step size $\tau_k = \frac{1}{19 + 10^{-3} k}$.
We chose the Adam parameters to be $\tau=0.01$ (step size), $\beta_1= 0.9$, and $\beta_2 = 0.99$, via grid search.
The fixed step size used for \algname with variable approximation space is $0.125$, while the memory-dependent one is given by the rule $\tau_k = \frac{1}{0.1165 m_k + 6.3433}$, where the numerical constants were tuned fitting the optimal step sizes found for the fixed approximation space algorithm.
The latter also justify the use of the above rule,
which is slightly different from the one used in \Cref{thm:conv-var,thm:conv_var1}, as it fits the data better.
\Cref{fig:diff_alg_comparison} shows that adapting the step size for the variable
\algname algorithm can result in a better convergence rate (notice the reduced rate once the error ~$10^{-13}$ is reached), even though the difference is not striking in this case.
Moreover, we also notice that the variable \algname outperforms SAGA, taking almost half the iterations to reach convergence at the same cost (one gradient evaluation per iteration).

\subsection{Advection-diffusion equation with five-dimensional random parameter}

We mimic the problem setting from \cite{doi:10.1137/18M1224076}.
We model the transport of a passive contaminant in the spatial domain $D\subset [0,1]^2$.
The control $u\in L^2(D)$ represents the rate at which a neutralizing chemical is injected; its purpose is to reduce the contaminant concentration $z$.
Uncertainty is introduced through the independent and identically distributed random variables $Y_1,\dots,Y_5\sim\mathcal U(0,1)$, collected in the random vector $Y:=(Y_1,\dots,Y_5)$, so that a realization of it is denoted $y:=(y_1,\dots,y_5)\in\Gamma:= [0,1]^5$: the wind field $V$, the diffusion coefficient $\varepsilon$, and the source term $f$ all depend on $y$.
Given a control $u\in L^2(D)$, the corresponding state $z(\cdot;u,y)\in H_{\Gamma_d}^1(D)$ is, for every $y\in\Gamma$, the solution of the following advection-diffusion equation given in weak form
\begin{equation}\label{eq:weak-state-5-dim}
    \int_D \bigl( \varepsilon(y)\,\nabla z\cdot\nabla v + V(y)\cdot\nabla z\,v \bigr)\,dx
    =\int_D \bigl( f(y)-u \bigr)v\,dx ,\qquad \forall v\in H_{\Gamma_d}^1(D) ,
\end{equation}
where homogeneous Dirichlet conditions are imposed on $\Gamma_d=\partial D$.
The data in \eqref{eq:weak-state-5-dim} are chosen as
\begin{equation*}
    \begin{aligned}
        f(x,y)         & =\exp\!\Bigl(-\|x-h(y)\|^2_2\Bigr), & h(y) & =(y_1,y_2), \\
        \varepsilon(y) & =0.5+\exp(y_3-1),                   &      &             \\
        V(x,y)         & =(y_4-y_5 x_1, y_5 x_2)^T .         &      &
    \end{aligned}
\end{equation*}
Thus, the pair $(y_1,y_2)$ parametrizes the source location, $y_3$ scales the diffusion coefficient, and $(y_4,y_5)$ set a linear advection field.
For a fixed regularization parameter $\beta>0$ we seek a control that minimizes the expected quadratic cost
\begin{equation}\label{eq:objective}
    J(u):= \Exp{g(u,Y)},
    \quad g(u,y) := \frac12 \norm{z(x,y)}_{L^2(D)}^2 + \frac{\beta}{2}\,\|u\|_{L^2(D)}^2,
\end{equation}
that is the optimal-control problem seeks $u$ that minimizes the expected total mass of contaminant while penalizing large injections.
To write the gradient of $J$ we introduce, for each $y\in\Gamma$, the adjoint variable $p(\cdot;u,y) \in H_{\Gamma_d}^1(D)$ which solves the equation
\begin{equation}\label{eq:adjoint-5-dim}
    \int_D \bigl( \varepsilon(y)\,\nabla v\cdot\nabla p + V(y)\cdot\nabla v p \bigr) dx
    =\int_D z v dx ,\qquad \forall v\in H_{\Gamma_d}^1(D) .
\end{equation}
The $L^2$-gradient then reads
\begin{equation}\label{eq:gradient-5-dim}
    \nabla J(u)=\Exp{\nabla g (u,Y)},
    \quad \nabla g (u,y) = p(\cdot,y) + \beta u.
\end{equation}
In the numerical experiments that follow we set $\beta=10^{-4}$ and discretize \eqref{eq:weak-state-5-dim} and \eqref{eq:adjoint-5-dim} using a piecewise linear finite element space $V_h$ on a uniform triangular mesh of size $h$ such that the number of total elements is $81$.

\begin{figure}[ht]
    \centering
    \includegraphics[width=\textwidth]{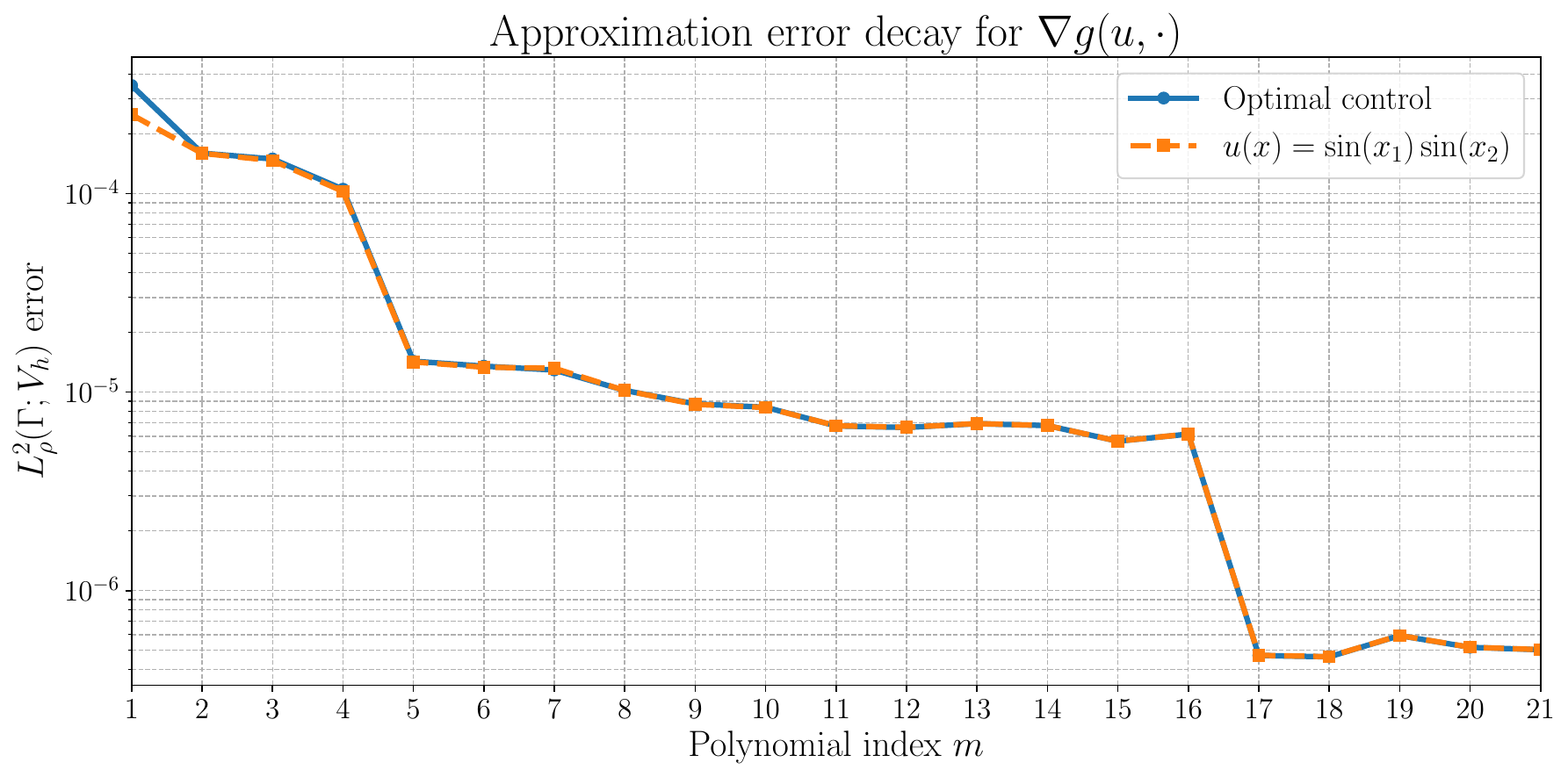}
    \caption{Estimate of the error of the polynomial approximation of the gradient $\norm{(I-\Pi^{\bbV_m}) [\nabla g(u,\cdot)]}_{L^2_\rho(\Gamma;L^2(D))}$ for $u$ the optimal control and $u(x)=\sin(x_1)\sin(x_2)$.}
    \label{fig:approx-error-decay}
\end{figure}

The parametric regularity of PDEs of the type \eqref{eq:weak-state-5-dim} has been studied in recent years (see, e.g., \cite{herrmann2019multilevel}).
It turns out that the maps $y \mapsto z(\cdot;u,y) \in H_{\Gamma_d}^1(D)$ and $y \mapsto p(\cdot;u,y) \in H_{\Gamma_d}^1(D)$ and, owing to \eqref{eq:gradient-5-dim}, $y \mapsto \nabla g(u,y) \in L^2(D)$ depend smoothly on $y$ for all $u \in L^2(D)$.
This justifies the use of polynomial approximation of $y \mapsto \nabla g(u,y)$ in this setting.
Indeed, the approximation spaces we employ in \Cref{alg:algo-fixed,algo:algo-variable} are polynomial spaces defined on the parametric domain $\Gamma = [0,1]^5$, as detailed in \Cref{sec:tensor-product-spaces}.
Specifically, we use multivariate polynomial spaces with an isotropic hyperbolic cross multi-index set defined by
\begin{equation*}
    \Lambda_m := \left\{\nu \in \mathbb{N}^5 \,:\, \sum_{j=1}^{5} \log(\nu_j+1) \leq \log(m+1) \right\},
\end{equation*}
so that $\calV_m = \mathrm{span}\{\phi_\nu \,:\, \nu \in \Lambda_m\}$, where each $\phi_\nu$ is the tensorized Legendre polynomial associated with the multi-index $\nu$ (note that, in this case, $\mathrm{dim}(\calV_m) = \abs{\Lambda_m} \neq m$).
However, this choice is not necessarily optimal.
Indeed, Figure \ref{fig:approx-error-decay} illustrates the decay of the approximation error when approximating $\nabla g(u, \cdot)$, both for the optimal control $u=u^\ast$ and for an arbitrary smooth control $u(x)=\sin(x_1)\sin(x_2)$. The observed stagnation of the approximation error between $m=5$ (with $|\Lambda_5|=56$) and $m=16$ (with $|\Lambda_{16}|=346$), and again after $m=17$ (with $|\Lambda_{17}|=421$), suggests that alternative multi-index designs could yield better performance.
On the other hand, the similarity between the two curves suggests that the polynomial space design to approximate $\nabla g(u,\cdot)$ features some robustness with respect to $u$, so that one could attempt building better index sets, e.g., employing greedy adaptive methods \cite{cohen2018multivariate}, without the knowledge of the optimal control (recall, indeed, that the theoretical results in \Cref{sec:conv-analyis} rely on the approximability of $\nabla g(u^\ast,\cdot)$).
Nevertheless, the isotropic hyperbolic cross structure was deliberately chosen to keep the polynomial space design simple and its dimension tractable.
We emphasize that our proposed algorithms provide considerable flexibility in designing these approximation spaces.
The sampling measure we employ is independent of the approximation space, and it is fixed to the tensorized arcsine distribution on $\Gamma$, that is $d\mu(y) = \prod_{i=1}^{5} (\pi \sqrt{y_i (1-y_i)})^{-1} dy_i$, which entails the weight supremum
\begin{equation*}
    q = \sup_{y \in \Gamma} \prod_{i=1}^{5} \pi \sqrt{y_i (1-y_i)}
    = \left( \frac{\pi}{2} \right)^5 \approx 9.5.
\end{equation*}
In contrast, the SAGA algorithm requires discretization of the functional through quadrature rules.
To leverage gradient regularity, one typically employs tensorized Gauss-Legendre quadratures; however, such quadratures become prohibitive in higher dimensions without sparsification. Sparsification itself introduces negative quadrature weights, compromising the convexity of the functional. Therefore, in high dimensions, alternative quadrature methods such as Monte Carlo or Quasi Monte Carlo become necessary, sacrificing the favorable convergence properties of Gauss-Legendre quadratures.
Following \cite{doi:10.1137/18M1224076}, our numerical experiments use SAGA with tensorized Gauss-Legendre quadratures employing either $5$ or $8$ points per dimension.
    {\color{black} Additionally, we include the "full gradient" solution to show that the method is computationally prohibitive in our setting. Indeed, each iteration of the full gradient algorithm requires solving $2\cdot 5^5$ PDEs, while each stochastic step solves only two PDEs.}
Since the exact minimizer is unavailable, the reference solution is computed using full-gradient descent with a $5$-point tensorized Gauss-Legendre quadrature.

\begin{figure}[ht]
    \centering
    \includegraphics[width=\textwidth]{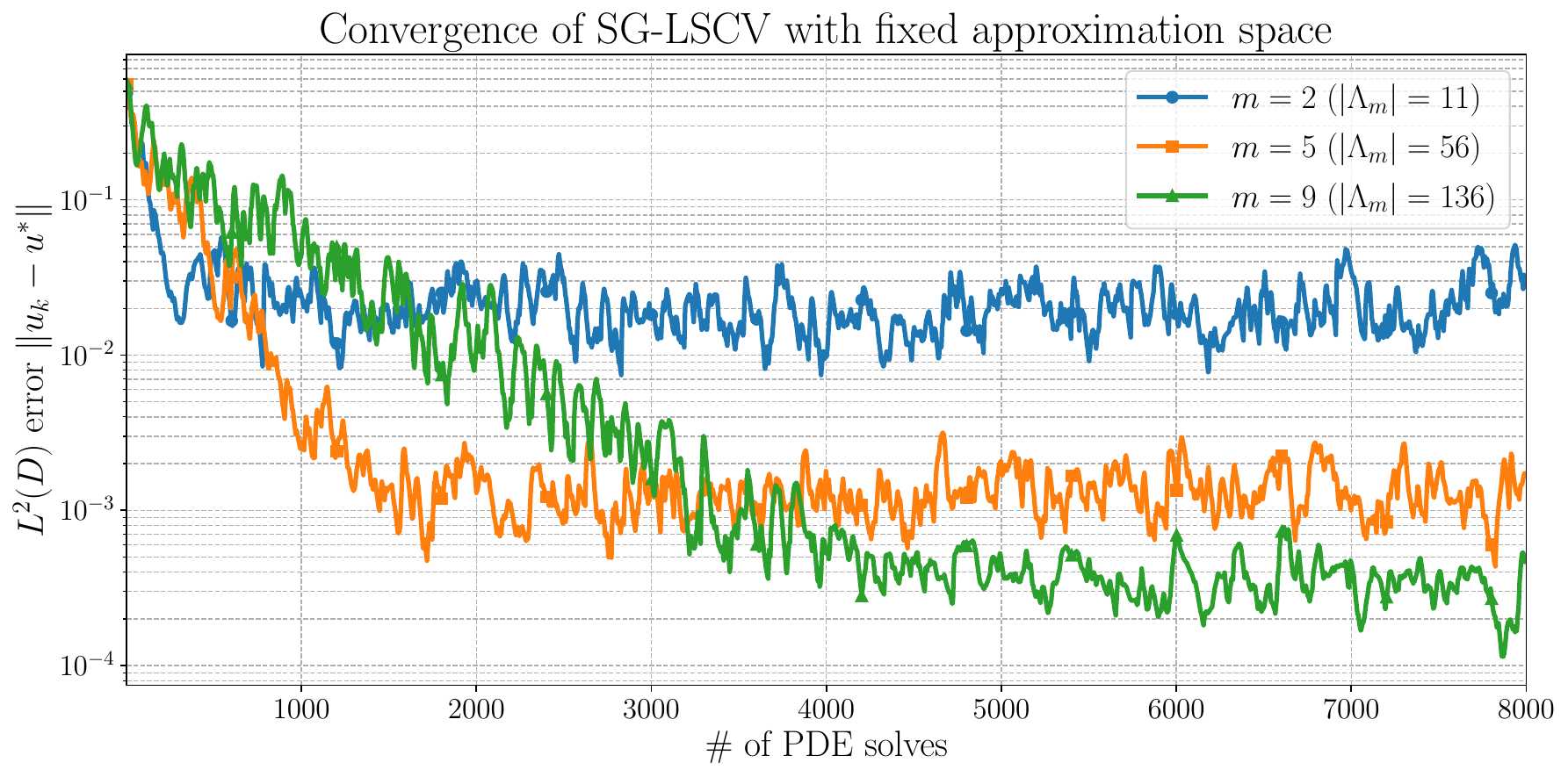}
    \caption{Convergence of the \algname algorithm with fixed approximation spaces $\bbV_m$ with $m\in\{2,5,9\}$ on the $5$-dimensional problem. The plot shows the exponential moving average of the error $\norm{u_k-u^\ast}_{L^2(D)}$.
            {\color{black} The algorithm uses one gradient evaluation per iteration, which corresponds to two PDE solves.}
    }
    \label{fig:fixed_sglscv_comparison}
\end{figure}

\Cref{fig:fixed_sglscv_comparison} displays the convergence of the \algname algorithm with fixed approximation spaces of increasing parameter {\color{black} $m \in \{2,5,9\}$}, corresponding respectively to polynomial spaces of cardinality $|\Lambda_m| = \{6, 56, 136\}$.
The methods employ step sizes tuned via grid-search, specifically $\tau=100,50,40$ for $m=2,5,9$, respectively.
As expected, higher-dimensional approximation spaces enable improved long-term accuracy, with the case $m=9$ achieving the lowest error plateau.
However, the benefit of increasing $m$ is not monotonic in early iterations: the method with $m=5$ initially outperforms $m=9$, consistently with \Cref{eq:thm-conv-fixed} and the observations in \Cref{sec:num-one-dim}.

\begin{figure}[ht]
    \centering
    \includegraphics[width=\textwidth]{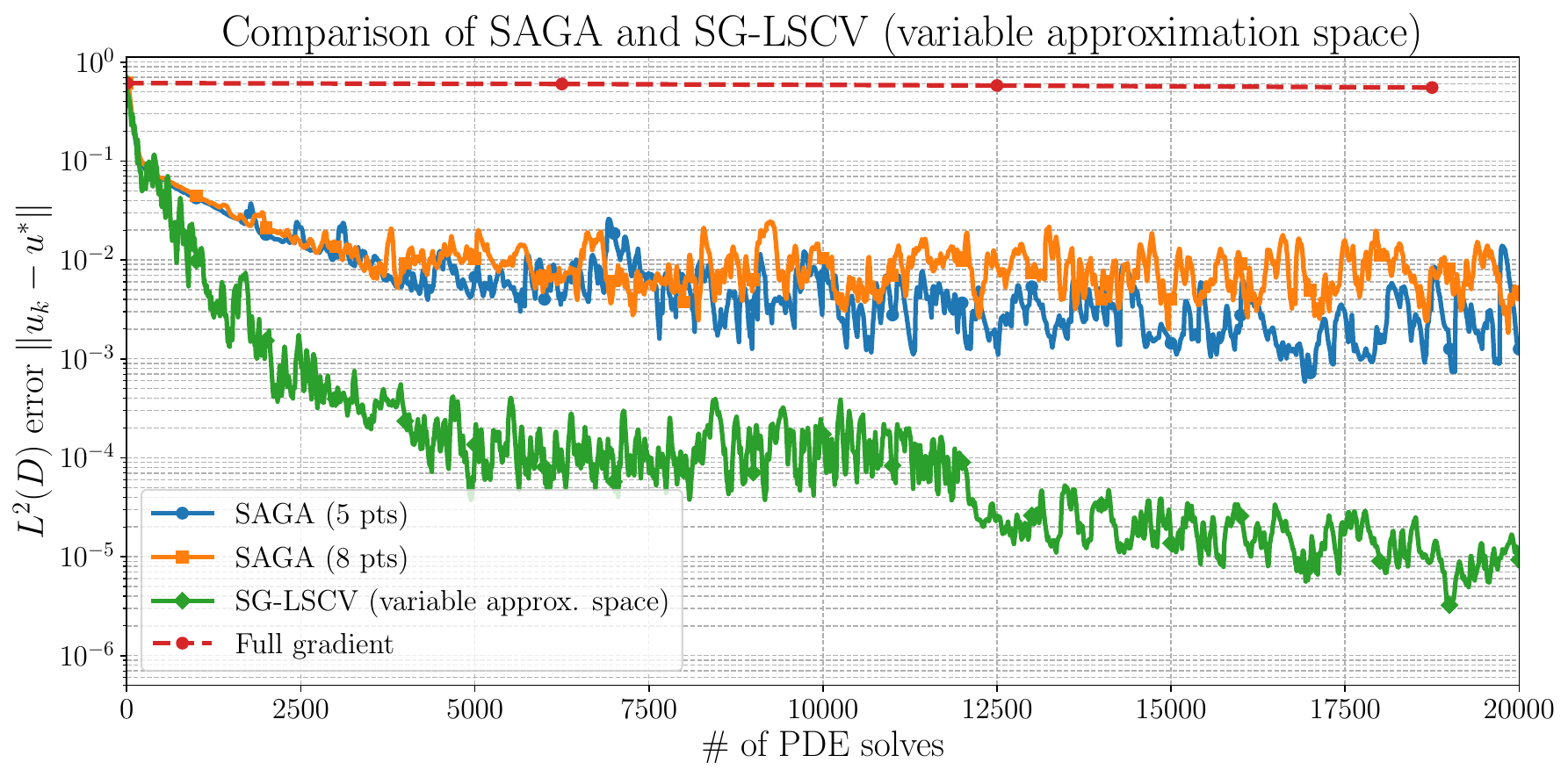}
    \caption{Comparison of stochastic optimization methods in 5 dimensions. The plot shows the exponential moving average of the error $\norm{u_k-u^\ast}_{L^2(D)}$ over $10^5$ iterations for SAGA with $5$ and $8$ quadrature points, and \algname with variable approximation spaces, {\color{black} which use one gradient evaluation per iteration, namely two PDE solves}. The full gradient method was run with step size $\tau=100$ {\color{black} and instead uses $5^5$ gradient evaluations, which corresponds to $2 \cdot 5^5$ PDE solves}.
    }
    \label{fig:5dim_alg_comparison}
\end{figure}

In \Cref{fig:5dim_alg_comparison}, we report convergence plots for the error $\|u_k - u^*\|_{L^2(D)}$ comparing two SAGA algorithms and the \algname algorithm with variable approximation space and memory-dependent step size.
The two SAGA variants differ in the number of quadrature points used in the discretization of the gradient: the first uses a $5$-point tensorized Gauss-Legendre quadrature, while the second uses $8$ points.
The \algname algorithm uses the step size rule $\tau_k = \frac{1}{2.3 \times 10^{-4} \abs{\Lambda_k} + 4.7 \times 10^{-3}}$, where the constants were found by fitting the optimal step sizes for the fixed approximation space algorithm, and incrementally enriches the polynomial approximation space by increasing the parameter $m$ in steps of three, skipping intermediate $m$'s and thus accelerating access to higher-fidelity surrogates.
This is important because of the stagnation phenomenon in $m$ of the approximation error observed in \Cref{fig:approx-error-decay}.

Indeed, we observe that the \algname algorithm significantly outperforms both SAGA variants in terms of convergence rate and final error.
Notably, around {\color{black} $12000$ PDE solves (iteration $6000$)}, \algname exhibits a marked drop in error, consistent with the activation of higher-degree approximations.
This behavior aligns with the sharp decrease in approximation error shown in Figure~\ref{fig:approx-error-decay} at $m=17$.
Since the update strategy for $m$ reaches $m=17$ around iteration $6000$, this explains the renewed decay observed in the convergence curve.

These results confirm that aggressively enriching the approximation space -- while controlling step size and stability of the least-squares approximation -- can yield faster and more robust convergence than methods relying on fixed discretizations, such as SAGA.

\section{Conclusions}
In this paper, we introduced and analyzed a stochastic optimization algorithm designed explicitly to address optimization problems involving expectations over continuous random variables.
Our method leverages memory-based control variates fitted on linear spaces through optimal weighted least-squares approximations, significantly reducing variance in gradient estimations and improving computational efficiency.
We proposed two versions, one with fixed approximation space and one where the approximation space size grows with the iterations.
We provided rigorous theoretical analyses, proving convergence for both fixed and variable approximation spaces under standard assumptions of strong convexity and Lipschitz continuity of the gradient.
In particular, we showed that, under the assumption that there exist linear spaces such that the projection error of $y \mapsto \nabla g(u^\ast,y)$ decays exponentially in the approximation space size, the error on the design converges to zero subexponentially fast, while, if the projection error decays algebraically, the convergence in the design is also algebraic.

We assessed the practical performance of our algorithm on PDE-constrained optimization problems with uncertain parameters, demonstrating its robust convergence behavior and computational advantages compared to existing variance reduction algorithms such as SAGA.
Our numerical experiments indicated substantial improvements in convergence speed, highlighting the efficacy of exploiting gradient regularity through optimal weighted polynomial approximation often present in these settings.

Our technique holds potential for broader applications, notably in machine learning contexts characterized by atomic probability measures.
In such discrete-data scenarios, one can construct polynomial bases orthogonal with respect to the empirical measure, thus naturally exploiting the smoothness of the gradient of the objective with respect to the data.
Unlike finite-sum algorithms such as SAGA, which inherently do not assume or exploit gradient regularity, our methodology explicitly leverages this regularity, thereby offering potential computational benefits and enhanced convergence properties in high-dimensional, data-intensive optimization problems.

\medskip
\noindent\textbf{Data availability statement} \\
The code used to produce the numerical results is freely available on Zenodo: \url{https://doi.org/10.5281/zenodo.16964631} \cite{raviola2025stochasticcode}.

\clearpage
\begin{appendices}

    \crefname{appendix}{appendix}{appendices}
    \Crefname{appendix}{Appendix}{Appendices}
    \crefalias{section}{appendix}      

    \section{\texorpdfstring{Proof of \Cref{prop:sup-weight}}{Proof of}}
    \label{sec:appendix-prop-sup-weight}

    The function $w_\Lambda(y)$ depends on $\calV_\Lambda$ and the measure $\rho$, but not on the chosen basis $\{\phi_\nu\}_{\nu \in \Lambda}$ \cite{cohen2016optimal}, hence, if the constant function is in $\calV_\Lambda$, we can assume without loss of generality that $\phi_{\bar{\nu}} = 1$ for some $\bar{\nu}\in \Lambda$.
    Then,
    \begin{equation*}
        w_\Lambda(y)
        = \frac{m}{1+\sum_{\nu \in \Lambda\setminus \bar{\nu}} \phi_\nu^2(y)}
        \leq m.
    \end{equation*}

    Assume now $\rho$ is Gaussian and $d=1$, so that we can use $\phi_j = H_{j-1}$, where $\{H_j\}_{j\geq 0}$ denotes the Hermite basis introduced in \Cref{exmp:hermite}, and replace the subscript $\Lambda$ with $m$.
    Define $S_m := \frac{m}{w_m} = \sum_{j=0}^{m-1} H_j^2$ and let us first show that $S_m$ has a minimum in $0$.
    Using the formulas
    \begin{equation*}
        H_{m+1} = \frac{1}{\sqrt{m+1}}\left(y H_m - \sqrt{m} H_{m-1} \right),
        \quad H_m' = \sqrt{m} H_{m-1},
    \end{equation*}
    we compute the derivative
    \begin{equation*}
        S'_m = 2 \sum_{j=1}^{m-1} \sqrt{j} H_{j-1} H_j
    \end{equation*}
    and notice that $S_m'(0)=0$.
    Furthermore, $S'_m$ is an odd function.
    Note then
    \begin{equation*}
        \begin{split}
            S'_m
             & = \sum_{j=1}^{m-1} \sqrt{j} H_{j-1} H_j + H_0H_1 + \sum_{j=2}^{m-1} H_{j-1} \left(y H_{j-1} - \sqrt{j-1} H_{j-2} \right) \\
             & = \sum_{j=1}^{m-1} \sqrt{j} H_{j-1} H_j + H_0H_1 + \sum_{j=1}^{m-2} y H_{j}^2 - \sum_{j=1}^{m-2} \sqrt{j} H_{j-1} H_{j}  \\
             & = \sqrt{m-1} H_{m-2} H_{m-1} + H_0H_1 + \sum_{j=1}^{m-2} y H_{j}^2                                                       \\
        \end{split}
    \end{equation*}
    Now, when $m \geq 2$
    \begin{equation*}
        \sqrt{m} H_{m-1} H_m =
        \sum_{j=1}^{m-1} (-1)^{m+j+1} y H_j^2 + (-1)^{m-1} H_0 H_1.
    \end{equation*}
    Indeed, $\sqrt{2}H_1H_2 = H_1(yH_{1} - H_0) = y H_1^2 - H_0 H_1$ while, assuming $m>2$ and the property to be true for $m-1$, we have
    \begin{equation*}
        \begin{split}
            \sqrt{m} H_{m-1} H_m
            = & H_{m-1} (y H_{m-1} - \sqrt{m-1} H_{m-2})                               \\
            = & y H_{m-1}^2 - \sqrt{m-1} H_{m-2} H_{m-1}                               \\
            = & y H_{m-1}^2 - \sum_{j=1}^{m-2} (-1)^{m+j} y H_j^2 - (-1)^{m-2} H_0 H_1 \\
            = & \sum_{j=1}^{m-1} (-1)^{m+j+1} y H_j^2 + (-1)^{m-1} H_0 H_1.
        \end{split}
    \end{equation*}
    Hence,
    \begin{equation*}
        S'_m
        =
        \begin{cases}
            2 \sum_{j=0}^{\frac{m-3}{2}} y H_{2j+1}^2          & \text{if $m$ is odd},  \\
            2 H_0H_1 + 2 \sum_{j=1}^{\frac{m-2}{2}} y H_{2j}^2 & \text{if $m$ is even}.
        \end{cases}
    \end{equation*}
    Since $H_0H_1 = y$, it holds that $S'_m(y) > 0$ if $y>0$ and $S'_m(y) < 0$ if $y<0$.
    Finally, observing that $S_m(y) \to +\infty$ as $\abs{y} \to \infty$ we conclude that $S_m$ has a minimum in $y=0$.

    We now aim to show that $S_m(0) \gtrsim \sqrt{m}$ using the fact that
    \begin{equation*}
        H_{2j}(0) = (-1)^j \frac{\sqrt{(2j)!}}{j!2^j}, \qquad H_{2j+1}(0) = 0,
    \end{equation*}
    which implies
    \begin{equation*}
        H_{2j}(0)^2 = \frac{(2j)!}{(j!)^2 \cdot 4^j} = \frac{1}{4^j} \binom{2j}{j},
        \qquad H_{2j+1}(0)^2 = 0.
    \end{equation*}
    Hence, we have
    \begin{equation*}
        S_m(0) = \sum_{j=0}^{\lfloor \frac{m-1}{2} \rfloor} \frac{1}{4^j} \binom{2j}{j}.
    \end{equation*}
    Using the bounds
    \begin{equation*}
        \sqrt{2 \pi j} \left( \frac{j}{e}\right)^j \leq j! \leq e \sqrt{j} \left( \frac{j}{e}\right)^j,
        \quad j\geq 1,
    \end{equation*}
    we obtain
    \begin{equation*}
        H_{2j}(0)^2 = \frac{1}{4^j} \binom{2j}{j} \geq \frac{2 \sqrt{\pi}}{e^2 \sqrt{j}}.
    \end{equation*}
    Hence,
    \begin{equation*}
        S_m(0) \geq 1+\sum_{j=2}^{\floor{\frac{m-1}{2}}} \frac{2 \sqrt{\pi}}{e^2 \sqrt{j}}
        \geq 1+ \frac{2 \sqrt{\pi}}{e^2} \int_2^{\frac{m-1}{2}} \frac{1}{\sqrt{t}} dt
        = 1+\frac{4 \sqrt{\pi}}{e^2} (\sqrt{\frac{m-1}{2}} - \sqrt{2})
        = 1+\frac{2 \sqrt{2 \pi}}{e^2} (\sqrt{m-1} - 2).
    \end{equation*}
    Therefore, for $m$ sufficiently large, the final result for the one-dimensional case is
    \begin{equation*}
        \max_y w_m(y) = \frac{m}{S_m(0)} \leq \frac{m}{1+\frac{2 \sqrt{2 \pi}}{e^2} (\sqrt{m-1} - 2)} \lesssim \sqrt{m}.
    \end{equation*}
    The result for all $m \geq 1$ follows by increasing the implicit constant.

    \section{\texorpdfstring{Complementary proofs from \Cref{sec:conv-analyis}}{Complementary proofs from}}

    \begin{lmm} \label{lemma_gamma_int}
        Assuming $\gamma >0$, $\tau>0$, $a\in\R$, and $0<t<1$, we have the following identity
        \begin{equation*}
            \int_{a}^{\infty} \exp\left( -2 \gamma \frac{\tau}{t} (x+1)^t \right)dx = \left( \frac{t}{2\gamma \tau} \right)^{\frac{1}{t}} \frac{\Gamma \left( \frac{1}{t}, \frac{2 \gamma \tau}{t} (a+1)^t \right)}{t},
        \end{equation*}
        where $\Gamma(r,x) := \int_{x}^{\infty} y^{r-1} e^{-y} dy$ is the upper incomplete $\Gamma$ function.
    \end{lmm}

    \begin{proof}
        Let us apply the following change of variable
        \begin{equation*}
            u = \left( \frac{2 \gamma \tau}{t} \right)^{\frac1t} (x+1).
        \end{equation*}
        Then, we have
        \begin{equation*}
            du = \left( \frac{2 \gamma \tau}{t} \right)^{\frac1t} dx.
        \end{equation*}
        We can therefore rewrite the integral $I:=\int_{a}^\infty \exp\left(-2\gamma \frac{\tau}{t} (x+1)^t\right) dx$ as
        \begin{equation*}
            \begin{split}
                I
                = \left( \frac{t}{2 \gamma \tau} \right)^{\frac{1}{t}} \int_{\left( \frac{2 \gamma \tau}{t} \right)^{\frac{1}{t}}(a+1)}^{\infty} e^{-u^t}du.
            \end{split}
        \end{equation*}
        Applying now the change of variable $u = y^{\frac{1}{t}}$, which implies $du = \frac{1}{t} y^{\frac{1}{t}-1}dy$, and using $y^{\frac{1}{t}-1} = u^{1-t}$, we have
        \begin{equation*}
            \int_{\left( \frac{2 \gamma \tau}{t} \right)^{\frac{1}{t}}(a+1)}^{\infty} e^{-u^t}du
            = \frac1t \int_{\left( \frac{2 \gamma \tau}{t} \right)(a+1)^t}^{\infty} e^{-y} y^{\frac1t -1} dy
            = \frac{1}{t} \Gamma\left(\frac1t, \frac{2\gamma\tau}{t}(a+1)^t\right).
        \end{equation*}
        Finally, we obtain the desired equality
        \begin{equation*}
            I = \int_{a}^{\infty} \exp\left( -2 \gamma \frac{\tau}{t} (x+1)^t \right)dx =  \left( \frac{t}{2 \gamma \tau} \right)^{\frac{1}{t}}   \frac{\Gamma \left( \frac{1}{t}, \frac{2 \gamma \tau}{t} (a+1)^t \right)}{t}.
        \end{equation*}
    \end{proof}

    \begin{lmm}[Endpoint Asymptotics for Exponential Integrals]
        \label{lem:end-asymp-int}
        Let $f \colon [a,\infty) \to (0, \infty)$ with $a \in \R$ be a continuously differentiable function and let $\phi \colon [a,\infty) \to \mathbb{R}$ be a twice continuously differentiable function such that $\phi'(x) > 0$ and $\phi''(x)$ exist for all $x \geq a$.
        Assume further that there exist $C, C', \xi, \eta>0$, $\delta >-1$, such that
        \begin{equation*}
            \frac{1}{C} x^{-\xi} \leq \frac{f(x)}{\phi'(x)} \leq C (1+x^\eta)
        \end{equation*}
        for all $x \geq a$, $\phi'(x) \geq C' x^{\delta}$, and $\abs{R(x)} \to 0$ as $x \to \infty$ with
        \begin{equation*}
            R(x) := \frac{f'(x)\phi'(x) - f(x)\phi''(x)}{f(x) \phi'(x)^2}.
        \end{equation*}
        Then, as $k \to \infty$,
        \begin{equation}\label{eq:endp-asympt}
            \int_a^k f(x)\, e^{\phi(x)}\, dx \sim \frac{f(k)}{\phi'(k)}\, e^{\phi(k)}.
        \end{equation}
        In particular, this is verified for $f(x) = x^\alpha$ and $\phi(x) = x^\beta$ for $\alpha \in \R$ and $\beta > 0$.
    \end{lmm}

    \begin{proof}
        Define for $c>b\geq a$
        \begin{equation*}
            I(b,c) := \int_b^c f(x) e^{\phi(x)} dx.
        \end{equation*}
        Then,
        {\color{black}
                \begin{equation*}
                    \begin{split}
                        I(b,c)
                        =  \int_b^c \frac{f(x)}{\phi'(x)} \phi'(x) e^{\phi(x)} dx
                        \leq  \max_{x\in [b,c]} \frac{f(x)}{\phi'(x)} (e^{\phi(c)}-e^{\phi(b)})
                        \leq  C (1 + c^\eta) (e^{\phi(c)}-e^{\phi(b)})
                    \end{split}
                \end{equation*}}
        and
        \begin{equation*}
            I(b,c) \geq \min_{x\in [b,c]} \frac{f(x)}{\phi'(x)} (e^{\phi(c)}-e^{\phi(b)})
            \geq \frac{1}{C} c^{-\xi} (e^{\phi(c)}-e^{\phi(b)}),
        \end{equation*}
        hence, up to a redefinition of $C$,
        \begin{equation*}
            \frac{I(a,\frac{k}{2})}{I(\frac{k}{2},k)}
            \leq C (1+k^{\eta+\xi}) \frac{e^{\phi\left(\frac{k}{2}\right)}-e^{\phi(a)}}{e^{\phi(k)}-e^{\phi\left(\frac{k}{2}\right)}}.
        \end{equation*}
        Now we have that there exists $\chi \in \left(\frac{k}{2},k\right)$ such that, up to redefinition of $C'$, $\phi(k)-\phi\left(\frac{k}{2}\right) = \phi'(\chi)\frac{k}{2} \geq C'k^{\delta+1}$.
        Since we assumed $\delta > -1$, we can conclude
        \begin{equation*}
            \frac{I(a,\frac{k}{2})}{I(\frac{k}{2},k)} \xrightarrow[k\to\infty]{}0.
        \end{equation*}

        We now rewrite $I\left(\frac{k}{2},k\right)$ using integration by parts, which gives
        \begin{equation*}
            I\left(\frac{k}{2},k\right)
            = \left[ \frac{f(x)}{\phi'(x)} e^{\phi(x)} \right]_{\frac{k}{2}}^k
            - \int_{\frac{k}{2}}^k \left( \frac{f'(x)\phi'(x) - f(x)\phi''(x)}{(\phi'(x))^2} \right) e^{\phi(x)} dx
            = \left[ \frac{f(x)}{\phi'(x)} e^{\phi(x)} \right]_{\frac{k}{2}}^k
            - \int_{\frac{k}{2}}^k R(x) f(x) e^{\phi(x)} dx.
        \end{equation*}
        Then, we have that
        \begin{equation*}
            \abs{\int_{\frac{k}{2}}^k R(x) f(x) e^{\phi(x)} dx} \leq I\left( \frac{k}{2}, k \right) \max_{x\in\left[\frac{k}{2},k\right]} \abs{R(x)} = o\left(I\left( \frac{k}{2}, k \right)\right),
        \end{equation*}
        which implies \eqref{eq:endp-asympt}.

        Finally, if $f(x) = x^\alpha$ and $\phi(x) = x^\beta$, the assumptions trivially hold with $\xi=\eta=\alpha-\beta+1$ and $\delta = \beta-1>-1$.
        This concludes the proof.
    \end{proof}

    \section{\texorpdfstring{Proof of \Cref{prop:one-dim}}{Proof of}}
    \label{sec:app-prop-one-dim}

    Consider the Laplacian operator $\Delta : H^1_0 \to L^2(D)$, which is invertible, and note that, dividing \eqref{eq:pde-one-dim} by $\tilde{y}>0$, we can write the solution as $z(u,y)= - \frac{1}{\tilde{y}}  \Delta^{-1}(u)$.
    Plugging this into \eqref{eq:min-prob-one-dim} gives
    \begin{equation*}
        \begin{split}
            J(u) & =\frac{1}{2} \mathbb{E}\left[ \norm{-\frac{1}{\tilde{Y}}\Delta^{-1}(u)-z_d}_{L^2(D)}^2 \right] + \frac{\beta}{2}\norm{u}_{L^2(D)}^2                                                                                                           \\
                 & = \frac{1}{2} \mathbb{E}\left[ \frac{1}{\tilde{Y}^2}\norm{\Delta^{-1}u}_{L^2(D)}^2+\frac{2}{\tilde{Y}} \langle \Delta^{-1}u,z_d  \rangle +\norm{z_d}_{L^2(D)}^2 \right] + \frac{\beta}{2}\norm{u}_{L^2(D)}^2                                  \\
                 & = \frac{1}{2} \left(  \mathbb{E}\left[ \frac{1}{\tilde{Y}^2}\right]\norm{\Delta^{-1}u}_{L^2(D)}^2 + 2\mathbb{E}\left[ \frac{1}{\tilde{Y}}\right] \langle \Delta^{-1}u,z_d\rangle + \norm{z_d}_{L^2(D)}^2 + \beta \norm{u}_{L^2(D)}^2  \right)
        \end{split}
    \end{equation*}
    Let us define the modified functional $\tilde{J}(u)$ by completing the square in the right-hand side of the last equality above -- which can be done adding a $u$-independent term -- and dividing everything by $\mathbb{E}[\tilde{Y}^{-2}]$, that is
    \begin{equation*}
        \begin{split}
            \tilde{J}(u) := \frac{1}{2} \left( \norm{\Delta^{-1}(u) + \alpha z_d}^2 + \tilde{\beta} \norm{u}^2\right),
        \end{split}
    \end{equation*}
    so that it holds $\mathrm{argmin}_u \tilde{J}(u) = \mathrm{argmin}_u J(u)$.

    Let us now compute $\mathrm{argmin}_u \tilde{J}(u)$.
    Recalling that the Fréchet derivative of $\norm{Au}_{L^2(D)}^2$, where $A$ is an operator from $L^2(D)$ to $L^2(D)$ with adjoint $A^\ast$, is $2A^*Au$, and using the fact that $z_d$ is an eigenfunction of the Laplacian associated to eigenvalue $\lambda$, namely
    \begin{equation}
        - \Delta^{-1}z_d = \lambda^{-1} z_d,
        \label{eq:eignfunc-one-dim}
    \end{equation}
    we can write the Fréchet derivative of $\tilde{J}$ as
    \begin{equation*}
        \begin{split}
            \nabla \tilde{J}(u)
             & = \nabla \frac{1}{2} \left[ \norm{\Delta^{-1}\!\left(u -  \alpha\lambda  z_d\right)}_{L^2(D)}^2 + \tilde{\beta}\norm{u}_{L^2(D)}^2 \right] \\
             & = \Delta^{-\ast}\Delta^{-1}\!\left[ u-  \alpha\lambda  z_d \right] + \tilde{\beta}u.
        \end{split}
    \end{equation*}
    Then, imposing $\nabla \tilde{J}(u)=0$ and using the fact that $\Delta^{-1}$ is self-adjoint, yields
    \begin{equation*}
        u = \alpha\lambda \Big[ \Delta^{-2} + \tilde{\beta} I  \Big]^{-1}\Delta^{-2} z_d,
    \end{equation*}
    which, owing to \eqref{eq:eignfunc-one-dim}, shows \eqref{eq:u-ast-prop}.

    Finally, fixing $u \in L^2(D)$, it is easy to see that
    \begin{equation*}
        \nabla g(u,y)
        = \frac{1}{\tilde{y}^{2}}\Delta^{-2} u + \frac{1}{\tilde{y}} \Delta^{-1}z_d + \beta u,
    \end{equation*}
    which is analytic in $\tilde{y}$ on the interval $[a,b]$ with $0<a<b<\infty$.
    Furthermore, owing to \eqref{eq:tildey-one-dim} $\tilde{y}$ is analytic in $y$ on the interval $[-1,1]$, whence $y\mapsto \nabla g(u,y)$ is analytic as a composition of analytic functions.
    This concludes the proof.

\end{appendices}
\bibliographystyle{plain}
\bibliography{bibliography}
\end{document}